\documentclass[leqno,11pt]{article}  %putting formula numbers on the left side
\usepackage{amsmath}
\usepackage{amssymb}

\usepackage[german,english]{babel}
\usepackage{amsthm}
\usepackage{xypic}  % commutative diagrams

\title{On the $p$-adic cohomology of some $p$-adically uniformized varieties}

\author{\textsc{Elmar Grosse-Kl\"onne}}
\date{}

\theoremstyle{plain} 
\newtheorem{satz}{Theorem}[section]  %@@
\newtheorem{kor}[satz]{Corollary}  %@@
\newtheorem{lem}[satz]{Lemma}  %@@
\newtheorem{pro}[satz]{Proposition}  %@@
\newcommand{\ho}{\mbox{\rm Hom}}  %@@
 
\newcommand{\spec}{\mbox{\rm Spec}}  %@@

\newcommand{\quot}{\mbox{\rm Quot}}  %@@
\newcommand{\spf}{\mbox{\rm Spf}}  %@@

\newcommand{\bi}{\mbox{\rm im}}  %@@
\newcommand{\ke}{\mbox{\rm Ker}}  %@@
\newcommand{\gr}{\mbox{\rm gr}}

\newcommand{\id}{\mbox{\rm id}}

\newcommand{\dlog}{\mbox{\rm dlog}}

\newcommand{\diag}{\mbox{\rm diag}}

\theoremstyle{remark}

\theoremstyle{definition}

\DeclareMathOperator{\Hom}{Hom}

\begin{document}
\maketitle
\footnote[0]
    {2000 \textit{Mathematics Subject Classification}.
    Primary 14F30}                               %2000MSC 
\footnote[0]{\textit{Key words and phrases}. Drinfel'd's symmetric space, Hodge decomposition, filtered $(\Phi,N)$-module}

\begin{abstract}
Let $K$ be a finite extension of ${\mathbb Q}_p$ and let $X$ be Drinfel'd's symmetric space of dimension $d$ over $K$. Let $\Gamma\subset {\rm SL}\sb {d+1}(K)$ be a cocompact discrete (torsionfree) subgroup and let ${{X}}_{\Gamma}=\Gamma\backslash {X}$, a smooth projective ${{K}}$-variety. In this paper we investigate the de Rham and log crystalline (log convergent) cohomology of local systems on $X_{\Gamma}$ arising from $K[\Gamma]$-modules. (I) We prove the monodromy weight conjecture in this context. To do so we work out, for a general strictly semistable proper scheme of pure relative dimension $d$ over a cdvr of mixed characteristic, a rigid analytic description of the $d$-fold iterate of the monodromy operator acting on de Rham cohomology. (II) In cases of arithmetical interest we prove the (weak) admissibility of this cohomology (as a filtered $(\phi,N)$-module) and the degeneration of the relevant Hodge spectral sequence.    
\end{abstract}

% \addtolength{\textwidth}{1.0in}
% \setlength{\hoffset}{-.5in}
%
%\addtolength{\topmargin}{-34pt}
%\addtolength{\textheight}{68pt}
%\baselineskip13,6pt
%\setlength{\parindent}{0pt}
%\setlength{\parskip}{0.8ex plus 0.2ex minus 0.2ex}
%\setlength{\oddsidemargin}{3cm}
%\frenchspacing
%\sloppy
%\pagestyle{myheadings}
%\markboth{}{}

%\textwidth14cm
%\textheight9in
%\hoffset+1in
%\voffset-0,8cm

\begin{center} {\bf Introduction} 
\end{center}

Recall that $X$ is the $K$-rigid analytic space complementary in ${\mathbb P}^d_K$ to all $K$-rational hyperplanes. The group $G={\rm GL}\sb {d+1}(K)$ acts naturally on $X$. The ${\ell}$-adic ($\ell\ne p$) cohomology of $X$ and its quotients $X_{\Gamma}$ are well understood. For example, by now we know that the cohomology of certain \'{e}tale coverings of $X$ realizes the local (smooth) Langlands correspondence for $G$, cf. e.g. \cite{dat}. Here instead we consider the de Rham and log crystalline (for our purposes equivalently: log convergent) cohomology (with coefficients) of the $X_{\Gamma}$ with all their additional structure elements provided by $p$-adic Hodge theory. The quite advanced knowledge in the case $G={\rm GL}\sb {2}({\mathbb Q}_p)$ has recently found impressive applications (e.g. \cite{iosp} and, in connection with a hoped for $p$-adically continuous local Langlands correspondence, \cite{bre}). This seems to be a good reason for trying to also deal with general $G={\rm GL}\sb {d+1}(K)$. 

{\bf (a)} Let ${\mathfrak X}$ be the natural strictly semistable $G$-equivariant formal ${\mathcal O}_K$-scheme underlying $X$, let ${\mathfrak{X}}_{\Gamma}=\Gamma\backslash {\mathfrak X}$. Let $M$ be a $K[\Gamma]$-module with $\dim_K(M)<\infty$. The associated $\Gamma$-equivariant constant sheaf $M$ on the rigid space $X$, resp. the formal scheme ${\mathfrak X}$ (in its Zariski topology), descends to a locally constant sheaf $M^{\Gamma}$ on $X_{\Gamma}$ (or even a local system, in the ordinary topological sense, on the Berkovich analytic space associated with $X_{\Gamma}$), resp. on the formal scheme ${\mathfrak{X}}_{\Gamma}$. (Thus $M^{\Gamma}$ does {\it not} mean the subspace of $\Gamma$-invariants, or its associated constant sheaf, of the {\it abstract} $K[\Gamma]$-module $M$.) Similarly, the sheaf complex $M\otimes_K\Omega_X^{\bullet}$, resp. $M\otimes_{{\mathcal O}_K}\Omega_{\mathfrak X}^{\bullet}$, with differential ${\rm id}_M\otimes d$, descends to a sheaf complex $M^{\Gamma}\otimes_K\Omega_{X_{\Gamma}}^{\bullet}$, resp. $M^{\Gamma}\otimes_{{\mathcal O}_K}\Omega_{{\mathfrak X}_{\Gamma}}^{\bullet}$. As these complexes have interesting global cohomology only in middle degree $d$, our central object of study is the de Rham cohomology group$$H_{dR}^d(X_{\Gamma},M)=H^d(X_{\Gamma},M^{\Gamma}\otimes_K\Omega_{X^{\Gamma}}^{\bullet})=H^d({\mathfrak X}_{\Gamma}, M^{\Gamma}\otimes_{{\mathcal O}_K}\Omega_{{\mathfrak X}_{\Gamma}}^{\bullet})$$endowed with additional structure elements as follows.\\(i) There is a covering spectral sequence\begin{gather}E_2^{r,s}=H^r(\Gamma,M\otimes_KH^s_{dR}(X))\Rightarrow H^{r+s}_{dR}(X_{\Gamma},M)\label{einlcovss}\end{gather}providing $H_{dR}^d(X_{\Gamma},M)$ with the covering filtration $H^{d}_{dR}(X_{\Gamma},{{M}})=F^0_{\Gamma}\supset F^1_{\Gamma}\supset\ldots\supset F^{d+1}_{\Gamma}=0$.\\(ii) The 'stupid' filtration $(M\otimes_K \Omega^{\bullet}_{X,\ge i})_{i\ge 0}$ gives rise to the spectral sequence\begin{gather}E_1^{r,s}=H^s(X_{\Gamma},{{M}}^{\Gamma}\otimes_K\Omega^{r}_{X_{\Gamma}} )\Rightarrow H^{r+s}_{dR}(X_{\Gamma},{{M}})\label{einlhodss}\end{gather}with corresponding filtration $H^{d}_{dR}(X_{\Gamma},{{M}})=F^0_{H}\supset F^1_{H}\supset\ldots\supset F^{d+1}_{H}=0.$\\(iii) If more specifically $M$ is (the restriction to $\Gamma$) of an (irreducible) $K$-rational representation of $G$ then there is another $G$-stable filtration $({\mathcal F}^{r,\bullet})_r$ on $(M\otimes_K \Omega^{\bullet}_{X})$, given by filtrations of the $(M\otimes_K \Omega^{i}_{X})$ by certain $G$-stable coherent ${\mathcal O}_X$-submodules, hence another spectral sequence\begin{gather}E_1^{st}=H^{s+t}(X_{\Gamma},({\mathcal F}^{s,\bullet}/{\mathcal F}^{s+1,\bullet})^{\Gamma})\Rightarrow H^{s+t}(X_{\Gamma},{{M}}^{\Gamma}\otimes_K\Omega^{\bullet}_{{X}_{\Gamma}} )\label{einltilhob}\end{gather} with corresponding (suitably renumbered) filtration $H^{d}_{dR}(X_{\Gamma},M)=F^0_{red}\supset F^1_{red}\supset\ldots\supset F^d_{red}\supset F^{d+1}_{red}=(0).$\\(iv) General facts on log crystalline, log convergent and log rigid cohomology and how they are related can be learned from \cite{faalet}, \cite{ogusII}, \cite{ogusc}, \cite{shiho1}, \cite{shiho2}, \cite{shihoArX}. Let $T={\rm Spf}(W(k),({\mathbb N}\to W(k), 1\mapsto0))$, the 'Hyodo-Kato' formal base log scheme. The sheaf $M^{\Gamma}$ on the reduction ${\mathfrak X}_{\Gamma,0}$ of ${\mathfrak X}_{\Gamma}$ gives rise to a locally constant isocrystal $M^{\Gamma}\otimes{\mathcal O}$ on ${\mathfrak X}_{\Gamma}$ and we may consider its log convergent cohomology $H_{conv}^d( {\mathfrak X}_{\Gamma,0}/T,M)$ with respect to $T$. (We might just as well work with log crystalline cohomology, for our purposes here this does not make a difference.) We have a Hyodo-Kato-type isomorphism\begin{gather} H_{conv}^d( {\mathfrak X}_{\Gamma,0}/T,M)\cong H^{d}_{dR}(X_{\Gamma},M)\label{einlhk}\end{gather}(depending on the choice of a uniformizer in ${\mathcal O}_K$). Since the general theory provides $H_{conv}^d( {\mathfrak X}_{\Gamma,0}/T,M)$ with a nilpotent linear monodromy operator $N$ we obtain, by transport of structure, a nilpotent linear monodromy operator $N$ on $H^{d}_{dR}(X_{\Gamma},M)$ (independent on the choice of a uniformizer in ${\mathcal O}_K$). We may consider the corresponding monodromy filtration, the convolution of the image- and the kernal filtration for $N$.\\(v) Suppose that in addition we are given a $K[\Gamma]$-linear automorphism $\Phi$ of $M$, a power of which acts as a scalar $\alpha\in K$. Tensoring with the $f$-fold iterate --- $f$ the degree of the residue field of $K$ over its prime field --- of the Frobenius endomorphism of the constant $F$-isocrystal on ${\mathfrak X}_{\Gamma,0}$ we obtain a $K$-linear Frobenius endomorphism $\Phi$ on $H_{conv}^d( {\mathfrak X}_{\Gamma,0}/T,M)$, hence on $H^{d}_{dR}(X_{\Gamma},M)$ (via the isomorphism (\ref{einlhk})). It satisfies $N\Phi=p^f\Phi N$ and provides $H^{d}_{dR}(X_{\Gamma},M)$ with a ($p$-adic) slope filtration, and if $\alpha$ is a Weil number then this is also a weight filtration (in the archimedean sense). 

{\bf (b)} The filtrations $F^{\bullet}_{\Gamma}$, $F^{\bullet}_{H}$ and $F^{\bullet}_{red}$ have been defined and studied by Schneider \cite{schn}. He proved the degeneration of (\ref{einlcovss}) in $E_2$ (see also \cite{hk} for another proof) and he computed the vector space dimensions of the graded pieces of the filtration $F^{\bullet}_{\Gamma}$. He conjectured that for any $M$ we have\begin{gather}H_{dR}^{d}(X_{\Gamma},{{M}})=F_H^{i+1}\bigoplus F_{\Gamma}^{d-i}\quad\quad(0\le i\le d-1).\label{einlspliteins}\end{gather}He conjectured that if $M$ is an (irreducible) $K$-rational representation of $G$ then the spectral sequence (\ref{einltilhob}) degenerates in $E_1$ and we have $F_{red}^j=F_H^{j}$ for all $j$, in particular he conjectured \begin{gather}H_{dR}^{d}(X_{\Gamma},{{M}})=F_{red}^{i+1}\bigoplus F_{\Gamma}^{d-i}\quad\quad(0\le i\le d-1).\label{einlspliteinsred}\end{gather}Iovita and Spiess \cite{iovspi} first proved all these conjectures for the trivial $G$-representation $M=K$ (and there is another proof by Alon and de Shalit), later we proved it furthermore for the standard representation $M=K^{d+1}$ of $G$ and its dual \cite{latt}. In \cite{hk} we showed that if $\Phi$ is as in (a)(v) then the filtration $F^{\bullet}_{\Gamma}$ is just the slope (resp. weight) filtration. In 2003, de Shalit \cite{desh} and Ito \cite{ito} independently proved the purity of the monodromy filtration ('monodromy weight conjecture'), i.e. the coincidence of the monodromy filtration with the weight filtration, in the case $M=K$ (with trivial $\Phi$). For $d=1$ the validity of all these conjectures is known for $K$-rational representations of $G$.

Finally we mention that in the ${\ell}$-adic setting ($\ell\ne p$), Dat \cite{dat} proved (for any $d$ and arbitrary coefficients) the analog of the purity of the monodromy filtration.

{\bf (c)} The first purpose of the present text is to prove the purity of the monodromy filtration for general $M$ and $\Phi$ as in (a)(v). Since, as just recalled, the slope (resp. weight) filtration coincides with $F^{\bullet}_{\Gamma}$ we need to show that $F^{\bullet}_{\Gamma}$ is the monodromy filtration. The computation of $F^{\bullet}_{\Gamma}$ in \cite{schn} and some elementary linear algebra show that it is enough to prove that $N^d$, the $d$-fold iterate of $N$, acts as non-trivially as possible on $H_{dR}^{d}(X_{\Gamma},{{M}})$. So we must understand $N^d$.

In general, for a strictly semistable (formal) ${\mathcal O}_K$-scheme ${\mathfrak Y}$ and a locally constant sheaf $E$ of finite dimensional $K$-vector spaces on the special fibre ${\mathfrak Y}_0$, the monodromy operator $N$ on $H_{conv}^*({\mathfrak Y}_0/T,E)$ provides a monodromy operator $N$ on $H_{dR}^*({\mathfrak Y}\otimes K,E)$ via the isomorphism $H_{conv}^*({\mathfrak Y}_0/T,E)\cong H_{dR}^*({\mathfrak Y}\otimes K,E)$, and it is an open problem to what extent $N$ can be described in terms of rigid analysis on ${\mathfrak Y}\otimes K$ alone, i.e. avoiding the passage to $H_{conv}^*({\mathfrak Y}_0/T,E)$. Assume that ${\mathfrak Y}\otimes K$ is of pure dimension $d$. One of the main results of \cite{coliov} is such a description of $N$ in the case $d=1$ and trivial $E$. Here we generalize this in that we desribe $N^d$ on $H^d({\mathfrak Y}\otimes K,E)$ for general $d$ and $E$. Namely, we show $N^d=\alpha\circ{\rm res}\circ\beta$, where $\alpha$ and $\beta$ are connecting homomorphisms in the canonical (corresponding to the decomposition of ${\mathfrak Y}_0$ into irreducible components) {C}ech spectral sequence converging to $H^*({\mathfrak Y}\otimes K,E)$, and where ${\rm res}$ is essentially given by residue maps. That is, ${\rm res}$ is given by expanding $d$-forms on generalized annuli (the preimages, under the specialization map ${\mathfrak Y}\otimes K\to{\mathfrak Y}_0$, of $(d+1)$-fold intersections of irreducible components of ${\mathfrak Y}_0$) into convergent Laurent series in uniformizing variables $\{T_i\}_i$ and extracting the coefficients of $\wedge_i{\rm dlog}(T_i)$.

Coming back to our varieties $X_{\Gamma}$ we then see that $N^d$ is essentially given by a residue map\begin{gather}H^0(\Gamma,M\otimes H_{dR}^d(X))\longrightarrow H^d(\Gamma,M)\label{einlres}\end{gather}and that our task is to prove its bijectivity. Let ${\mathcal C}_{har}^d(M)$ denote the space of $M$-valued harmonic $d$-cochains on the Bruhat-Tits building of ${\rm SL}_{d+1}(K)$. The residue map (\ref{einlres}) factors naturally as$$H^0(\Gamma,M\otimes H_{dR}^d(X))\longrightarrow H^0(\Gamma,{\mathcal C}_{har}^d(M) )\longrightarrow H^d(\Gamma,M).$$That the first arrow (even $M\otimes H_{dR}^d(X))\to{\mathcal C}_{har}^d(M)$) is an isomorphism follows from one of the main results in \cite{schtei}. That the second arrow is an isomorphism has been proved in \cite{gar}.

{\bf (d)} Let $\breve{K}$ be the completion of the maximal unramified extension of $K$. It is well known, cf. the book \cite{rz} by Rapoport and Zink, that $\breve{\mathfrak X}={\mathfrak X}\otimes_{{\mathcal O}_K}{\mathcal O}_{\breve{K}}$ is a moduli space for certain formal $p$-divisible groups endowed with some additional structures. The corresponding universal $p$-divisible group ${\mathcal G}^u$ over $\breve{\mathfrak X}$ then gives rise, by Dieudonn\'{e} theory, to a $G$-equivariant filtered (convergent) $F$-isocrystal $\breve{E}$ on $\breve{\mathfrak X}$. This is essentially the datum of a filtered ${\mathcal O}_{\breve{\mathfrak X}}\otimes_{{\mathcal O}_{\breve{K}}}\breve{K}$-module $\breve{E}(\breve{\mathfrak X})$ with connection, together with Frobenius structures on evaluations of $\breve{E}$ at Frobenius thickenings of the special fibre $\breve{\mathfrak X}_0$ of $\breve{\mathfrak X}$. It is certainly of great interest (e.g. with an eye towards applications for a hoped for $p$-adically continuous local Langlands correspondence for $G$) to understand the de Rham/log crystalline/log convergent cohomology, as a filtered $(\phi,N)$-module, of all tensor powers $\breve{E}^{\otimes r}$ of $\breve{E}$. These split up into certain direct summands $\breve{E}_M$ associated with $K$-rational representations $M$ of $G$. (We have $\breve{E}=\breve{E}_{{\mathcal S}}$ for the standard representation ${\mathcal S}=K^{d+1}$ of $G$.) Serving as prominent example material the various structural features of $\breve{E}$ are discussed in \cite{rz} as the development of the general theory in that book proceeds. We gathered from \cite{rz} once more in full detail the explicit description of $\breve{E}$, thereby deriving one for all the $\breve{E}_M$. The relevance with respect to the first part comes from the following fact. For $\delta\in{\rm Gal}(K/{\mathbb Q}_p)$ let $M_{\delta}$ denote the $\delta$-twist of the $K[G]$-module $M$. Then the $G$-equivariant ${\mathcal O}_{\breve{\mathfrak X}}\otimes_{{\mathcal O}_{\breve{K}}}\breve{K}$-module with connection $\breve{E}_M(\breve{\mathfrak X})$ is a direct sums of copies of $M_{\delta}\otimes_K({\mathcal O}_{\breve{\mathfrak X}}\otimes_{{\mathcal O}_{\breve{K}}}\breve{K})$ (all ${\delta}$) such that\\--- the connection is just ${\rm id}_{M_{\delta}}\otimes d$,\\--- the filtration on $M_{\delta}\otimes_K({\mathcal O}_{\breve{\mathfrak X}}\otimes_{{\mathcal O}_{\breve{K}}}\breve{K})$ gives rise to a filtration on the de Rham complex $M_{\delta}\otimes_K(\Omega^{\bullet}_{\breve{\mathfrak X}}\otimes_{{\mathcal O}_{\breve{K}}}\breve{K})$ which is of the sort considered in (a)(iii) if $\delta={\rm id}_K$, and of the sort considered in (a)(ii) if $\delta\ne {\rm id}_K$, and\\--- a suitable power of the Frobenius action on $\breve{E}_M$ stabilizes each of its direct summands and acts on each of them like in (a)(v).\\ In other words, $\breve{E}_M$ (with a suitable power of its Frobenius structure) arises by base field extension $K\to \breve{K}$ from a direct sum of filtered convergent $F$-isocrystals on ${\mathfrak X}$ as considered in the first part.

For suitable $\Gamma$ it happens that $\breve{\mathfrak X}_{\Gamma}=\Gamma\backslash\breve{\mathfrak X}$ arises by base field extension to $\breve{K}$ of some Shimura variety ($p$-adic uniformization). In that case ${\mathcal G}^u$ descends to a $p$-divisible group ${\mathcal G}^u_{\Gamma}$ over $\breve{\mathfrak X}_{\Gamma}$, and the product of ${\mathcal G}^u_{\Gamma}$ with its dual is the $p$-divisible group of a universal abelian scheme ${\mathcal A}$ over $\breve{\mathfrak X}_{\Gamma}$. Moreover, $\breve{E}_M$ descends to a convergent filtered $F$-isocrystal on $\breve{\mathfrak X}_{\Gamma}$ which is a direct summand in the relative de Rham/crystalline cohomology of some power ${\mathcal A}^r$ of the $\breve{\mathfrak X}_{\Gamma}$-scheme ${\mathcal A}$. One deduces that the absolute (i.e. over ${\mathcal O}_{\breve{K}}$) de Rham/log crystalline cohomology groups (with constant coefficients) of the strictly semistable ${\mathcal O}_{\breve{K}}$-scheme ${\mathcal A}^r$ are direct sums, as filtered $(\phi,N)$-modules, of the cohomology groups of the various $\breve{E}_M$ --- hence of the cohomology groups of the (base field extensions of) filtered convergent $F$-isocrystals as considered in the first part. Similarly, the Hodge spectral sequence for ${\mathcal A}^r\otimes_{{\mathcal O}_{\breve{K}}}\breve{K}\to\breve{K}$ is the direct sum of various spectral sequences of the form (\ref{einlhodss}) and (\ref{einltilhob}) (base field extended). In one direction, this implies that the purity of the monodromy filtration proven in the first part implies the purity of the monodromy filtration on the (absolute) log crystalline cohomology (with constant coefficients) of ${\mathcal A}^r$ (any $r$). In the reverse direction we have two applications. Firstly, the (well known) degeneration of the Hodge spectral sequence for ${\mathcal A}^r\otimes_{{\mathcal O}_{\breve{K}}}\breve{K}\to\breve{K}$ implies the degeneration of the spectral sequences (\ref{einlhodss}) for all {\it non-trivial} twists of $K$-rational $G$-representations, and of the spectral sequences (\ref{einltilhob}). The degeneration of the spectral sequences (\ref{einlhodss}) for non-trivial twists of $K$-rational $G$-representation came to us as a real surprise since for non-twisted $K$-rational $G$-representations the spectral sequences (\ref{einlhodss}) almost never degenerate. Secondly, the weak admissibility of the cohomology of ${\mathcal A}^r$, as proven by Tsuji, implies the weak admissibility of the cohomology of the $\breve{E}_M$. Since on the latter the slope filtration for Frobenius coincides with $F^{\bullet}_{\Gamma}$, this may be interpreted as a weak form for the decomposition conjectures (\ref{einlspliteins}) (for all non-trivial twists of $K$-rational $G$-representations) and (\ref{einlspliteinsred}).

{\bf (e)} The text is organized as follows. In section \ref{conj} we recall the various spectral sequences computing $H^{d}_{dR}(X_{\Gamma},M)$ as well as Schneider's conjectures. In section \ref{holdisrep} we briefly mention the link he observed with certain 'holomorphic discrete series representations' of $G$. In section \ref{monosection} we prove the purity of the monodromy filtration. As explained, this rests on the rigid analytic description of $N^d$ which in the appendix sections \ref{resisec} and \ref{steenbrink} we explain for general strictly semistable formal schemes of relative dimension $d$ over a complete discrete valuation ring of mixed characteristic. We also recall in section \ref{steenbrink} some facts from \cite{hk}. In section \ref{filconfiso} we describe $K$-versions of the filtered convergent $F$-isocrystals $\breve{E}_M$, in section \ref{appel} we explain their genesis from the universal $p$-divisible group ${\mathcal G}^u$ and prove the results indicated in (d) above. In the appendix section \ref{unifsect} we recall the statement on $p$-adic uniformization of Shimura varieties in the setting relevant for us. The appendix section \ref{klaussec} discusses the Hodge spectral sequence of an abelian scheme over a smooth projective variety over a field of characteristic zero.\\

{\bf Acknowledgments:} I wish to thank K. K\"unnemann, M. Rapoport, P. Schneider and M. Strauch for helpful conversations related to this work. These results have been presented at the workshop on rigid analysis and its applictions held at Regensburg in February 2008, and at the Number Theory Seminar of the Hebrew University at Jerusalem in March 2008. For the invitation to the former I am grateful to G. Kings and K. K\"unnemann, for the invitation to the latter I heartily thank E. de Shalit. I thank the referee for carefully reading the text and providing helpful comments.\\

{\it Notations:} We fix a prime number $p$ and a finite extension $K$ of ${\mathbb{Q}}_p$ of degreee $n$ with ring of integers ${\cal O}_K$ and residue field $k$. Let $\pi\in K$ be a prime element. Let $K^t$ be the maximal subfield of $K$ unramified over ${\mathbb Q}_p$, let $f=[K^t:{\mathbb Q}_p]$ and let $\sigma$ denote the Frobenius endomorphism of $K^t$. We fix $d\in \mathbb{N}$ and write $G={\rm GL}\sb {d+1}(K)$. 

\section{De Rham cohomology and spectral sequences}

\label{conj}

Drinfel'd's symmetric space of dimension $d$ over $K$ is the $K$-rigid space$$X={\mathbb{P}^{d}_K}-(\mbox{the union of all $K$-rational hyperplanes}).$$The group $G$ acts on $X$ as follows: $g\in G$ sends a line through the origin and $(z_0,\ldots,z_d)\ne 0$ to the line through the origin and $(z_0,\ldots,z_d)g^{-1}$. Let ${\mathfrak{X}}$ be the strictly semistable formal ${\mathcal O}_K$-scheme with generic fibre $X$ introduced in \cite{mus}. The action of $G$ on $X$ extends naturally to an action on ${\mathfrak{X}}$. For any $0\le j\le d$ the set $F^j$ of non-empty intersections of $(j+1)$-many pairwise distinct irreducible components of ${\mathfrak{X}}\otimes_{{\mathcal O}_K}k$ is in natural bijection with the set of $j$-simplices of the Bruhat Tits building of ${\rm PGL}\sb {d+1}/K$. 

In the sequel, for sheaves ${\mathcal G}$ on $X$ we write ${\mathcal G}$ also for the push forward sheaf on $\mathfrak{X}$ under the specialization map $sp:X\to \mathfrak{X}$; we use tacitly and repeatedly Kiehl's result that if ${\mathcal G}$ is coherent we have ${\mathbb R}^tsp_*{\mathcal G}=0$ for all $t>0$.  

Let $\Gamma\subset {\rm PGL}\sb {d+1}(K)$ be a discrete torsionfree and cocompact subgroup such that the quotient ${{\mathfrak X}}_{\Gamma}=\Gamma\backslash{{\mathfrak X}}$ is a projective ${\mathcal O}_{{K}}$-scheme with strictly semistable reduction. Its generic fibre is ${{X}}_{\Gamma}=\Gamma\backslash {X}$, a smooth projective ${{K}}$-variety.

Let $M$ be a $K[\Gamma]$-module with $\dim_KM<\infty$; we write ${M}$ also for the constant sheaf on $\mathfrak{X}$ generated by $M$. It descends to a local system $M^{\Gamma}$ on $\mathfrak{X}_{\Gamma}$ as follows: $M^{\Gamma}$ is the sheaf associated to the presheaf which assigns to an open subset $U\subset \mathfrak{X}_{\Gamma}$ the space of $\Gamma$-invariant in $\underline{M}(U\times_{\mathfrak{X}_{\Gamma}}\mathfrak{X})$, where $\underline{M}$ denotes the constant sheaf on $\mathfrak{X}$ generated by $M$. The de Rham complex ${{M}}\otimes_K\Omega^{\bullet}_{{X}}$ with differential $\id_M\otimes d$ descends to a de Rham complex ${{M}}^{\Gamma}\otimes_K\Omega^{\bullet}_{{X}_{\Gamma}}$ on $\mathfrak{X}_{\Gamma}$. (Thus, unless the $\Gamma$-action on $M$ is trivial, ${{M}}^{\Gamma}\otimes_K\Omega^{\bullet}_{{X}_{\Gamma}}$ is {\it not} meant to denote the de Rham complex $\Omega^{\bullet}_{{X}_{\Gamma}}$ tensored with the constant sheaf with value the $\Gamma$-invariants of the abstract $K[\Gamma]$-module $M$.) 
We have the covering spectral sequence\begin{gather}E_2^{r,s}=H^r(\Gamma,M\otimes_KH^s_{dR}(X))\Rightarrow H^{r+s}(X_{\Gamma},{{M}}^{\Gamma}\otimes_K\Omega^{\bullet}_{{X}_{\Gamma}}    )\label{covss}\end{gather}which degenerates in $E_2$, as is shown in \cite{schn} (or \cite{hk}). If $M$ contains no non-zero $\Gamma$-invariant vector we have\begin{gather}H^{j}(X_{\Gamma},{{M}}^{\Gamma}\otimes_K\Omega^{\bullet}_{{X}_{\Gamma}}    )=0\quad\quad\quad(j\ne d).\label{othdeg}\end{gather}Denote by \begin{gather}H^{d}(X_{\Gamma},{{M}}^{\Gamma}\otimes_K\Omega^{\bullet}_{{X}_{\Gamma}}    )=F^0_{\Gamma}\supset F^1_{\Gamma}\supset\ldots\supset F^{d+1}_{\Gamma}=0\label{gamfil}\end{gather}the filtration on $H^{d}(X_{\Gamma},{{M}}^{\Gamma}\otimes_K\Omega^{\bullet}_{{X}_{\Gamma}}    )$ induced by (\ref{covss}). By \cite{schn} Theorem 2 and Proposition 2, section 1, we have for $i=0,\ldots,d+1$:\begin{gather}\dim_KF^i_{\Gamma}=\left\{\begin{array}{l@{\quad:\quad}l}(d+1-i)\mu(\Gamma,M)&\mbox{}\,\,d\,\,\mbox{is odd or}\,\,2i>d\\(d+1-i)\mu(\Gamma,M)&\mbox{}\,\,d\,\,\mbox{is even and}\,\,2i\le d\end{array}\right.\label{pecomp}\\\mu(\Gamma,M)=\mu(\Gamma,M^*)\label{dimdu}\end{gather}Here $\mu(\Gamma,M)=\dim_KH^d(\Gamma,M)$ and $M^*=\Hom_K(M,K)$ and we must assume $d\ge2$ and that $M$ contains no non-zero $\Gamma$-invariant vector. On the other hand we have the Hodge spectral sequence\begin{gather}E_1^{r,s}=H^s(X_{\Gamma},{{M}}^{\Gamma}\otimes_K\Omega^{r}_{X_{\Gamma}} )\Rightarrow H^{r+s}(X_{\Gamma},{{M}}^{\Gamma}\otimes_K\Omega^{\bullet}_{{X}_{\Gamma}}    )
\label{hodss}\end{gather}which gives rise to the Hodge filtration$$H^{d}(X_{\Gamma},{{M}}^{\Gamma}\otimes_K\Omega^{\bullet}_{{X}_{\Gamma}}    )=F^0_{H}\supset F^1_{H}\supset\ldots\supset F^{d+1}_{H}=0.$$

If $\Xi_0,\ldots,\Xi_{d}$ denote the standard projective coordinate functions on $\mathbb{P}^{d}_K$, then $z_j=\Xi_j/\Xi_0$ for $j=1,\ldots,d$ are holomorphic functions on $X$. Let$$\overline{u}(z)=
\left(\begin{array}{cc}1&-z_1\quad\cdots\quad-z_d\\{0}&I_d\end{array}\right)\in{\rm SL}\sb {d+1}({\mathcal O}_X(X)).$$

 For $0\le i\le d$ define the obvious cocharacter $e_i:{\mathbb{G}}_m\to{\rm GL}\sb {d+1}$, i.e. the one which sends $t$ to the diagonal matrix $e_i(t)$ with $e_i(t)_{ii}=t$, $e_i(t)_{jj}=1$ for $i\ne j$ and $e_i(t)_{j_1j_2}=0$ for $j_1\ne j_2$.

Let now $M$ be an irreducible $K$-rational representation of $G={\rm GL}\sb {d+1}(K)$ of highest weight $(\lambda_0\ge\lambda_1\ge\ldots\ge\lambda_d)$. By this we mean that there exists a non zero vector $m\in M$ such that $K.m$ is stable under upper triangular matrices and generates $M$ as a $G$-representation, and such that $gm=\prod_{i=0}^da_i^{\lambda_i}m$ for all diagonal matrices $g=e_0(a_0)\cdots e_d(a_d)\in G$. Replacing $M$ by $M\otimes{\rm det}^{-\lambda_d}$ we may assume $\lambda_d=0$ --- in fact all the following constructions depend only on $M$ viewed as an ${\rm SL}\sb {d+1}(K)$-representation. We then set $$r_M=\lambda_0+\ldots+\lambda_d=\lambda_0+\ldots+\lambda_{d-1}.$$ We grade $M$ by setting$${\gr}^s M=\{m\in M\quad|\quad e_0(a_0)m=a_0^{r_M-s}m\,\,\mbox{for all}\,\,a_0\in K\}$$for $s\in\mathbb{Z}$, and we filter $M$ by setting$${f}^sM=\bigoplus_{s'\ge s}{\gr}^{s'}M.$$Thus ${\gr}^s M$ resp. ${f}^s M$ as defined here is ${\gr}^{\lambda_0-r_M+s} M$ resp. ${f}^{\lambda_0-r_M+s} M$ as defined in \cite{latt}. We have ${f}^{r_M+1}M=0$ and ${f}^{r_M-\lambda_0}M=M$. We filter $M\otimes_K\Omega^{j}_X$ by setting\begin{gather}f^r(M\otimes_K\Omega^{j}_X)={\mathcal O}_X.\overline{u}(z)(f^rM)\otimes_{{\mathcal O}_X}\Omega^j_X.\label{verkl}\end{gather}We let$${\mathcal F}^{r,\bullet}=[f^{r}(M\otimes_K\Omega^0_X)\longrightarrow f^{r-1}(M\otimes_K\Omega^1_X)\longrightarrow f^{r-2}(M\otimes_K\Omega^2_X)\longrightarrow\ldots].$$

It follows from \cite{schn} that this is a ${\rm SL}\sb {d+1}(K)$-stable filtration of $M\otimes_K\Omega^{\bullet}_X$ by subcomplexes and that the differentials on the graded pieces ${\mathcal F}^{r,\bullet}/{\mathcal F}^{r+1,\bullet}$ are ${\mathcal O}_X$-linear (notations and normalizations in loc. cit. are different). We obtain the spectral sequence\begin{gather}E_1^{r,s}=h^{r+s}({\mathcal F}^{r,\bullet}/{\mathcal F}^{r+1,\bullet})\Rightarrow h^{r+s}(M\otimes_K\Omega^{\bullet}_X)\label{ogusfil}.\end{gather}The following is \cite{schn} Lemma 9, section 3 (observe that $X$ is a Stein space).

\begin{pro}\label{peterco} (Schneider) The terms ${D}^j(M)=E_1^{r_M-\lambda_j+j,\lambda_j-r_M}$ for $0\le j\le d$ are the only non vanishing $E_1$-terms in (\ref{ogusfil}).
\end{pro}

In particular this means that the inclusions ${\mathcal F}^{r_M-\lambda_{j}+j}\to{\mathcal F}^{r_M-\lambda_{j-1}+j}$ are quasiisomorphisms. Moreover, as explained in \cite{schn} (and recalled in \cite{latt}) it follows from \ref{peterco} that there are differential operators ${D}^j(M)\to {D}^{j+1}(M)$ for $0\le j\le d-1$ such that$${D}^0(M)\longrightarrow {D}^1(M)\longrightarrow{D}^2(M)\longrightarrow\ldots \longrightarrow {D}^d(M)$$is a ${\rm SL}\sb {d+1}(K)$-equivariant complex and such that the following holds: in the derived category $D({\mathfrak{X}})$ of abelian sheaves on ${\mathfrak{X}}$ there are canonical and ${\rm SL}_{d+1}(K)$-equivariant isomorphisms ${\mathcal F}^{r_M-\lambda_{j}+j}\cong{D}^{\bullet}(M)_{\ge j}$ for any $j$, compatible when $j$ varies (more precisely: there is third filtered ${\rm SL}_{d+1}(K)$-equivariant sheaf complex on ${\mathfrak{X}}$ whose $j$-th filtration step maps ${\rm SL}\sb {d+1}(K)$-equivariantly and quasiisomorphically to both ${\mathcal F}^{r_M-\lambda_{j}+j}$ and ${D}^{\bullet}(M)_{\ge j}$). In particular, ${M}\otimes_K\Omega^{\bullet}_X\cong{D}^{\bullet}(M)$.

Let $\Gamma<{\rm SL}\sb {d+1}(K)$ be as before. From what we just said it follows that the spectral sequences\begin{gather}E_1^{st}=H^{s+t}(X_{\Gamma},({\mathcal F}^{s,\bullet}/{\mathcal F}^{s+1,\bullet})^{\Gamma})\Rightarrow H^{s+t}(X_{\Gamma},{{M}}^{\Gamma}\otimes_K\Omega^{\bullet}_{{X}_{\Gamma}} )\label{tilhob}\end{gather}
\begin{gather}E_1^{st}=H^t(X_{\Gamma},{D}^{s}(M)^{\Gamma})\Rightarrow H^{s+t}(X_{\Gamma},{D}^{\bullet}(M)^{\Gamma})=H^{s+t}(X_{\Gamma},M^{\Gamma}\otimes_K\Omega^{\bullet}_{{X}_{\Gamma}} )\label{redss}\end{gather}are isomorphic. Let$$H^{d}(X_{\Gamma},M^{\Gamma}\otimes_K\Omega^{\bullet}_{{X}_{\Gamma}} )={F}^{r_M-\lambda_0}_I\supset{F}^{r_M-\lambda_0+1}_{I}\supset\ldots\supset{F}^{r_M+d}_{I}\supset{F}^{r_M+d+1}_{I}=(0)$$be the filtration induced by (\ref{tilhob}), let $$H^{d}(X_{\Gamma},M^{\Gamma}\otimes_K\Omega^{\bullet}_{{X}_{\Gamma}} )=F^0_{red}\supset F^1_{red}\supset\ldots\supset F^d_{red}\supset F^{d+1}_{red}=(0)$$be the filtration induced by $(\ref{redss})$. From \ref{peterco} it follows that for all $d\ge j\ge 1$ we have\begin{gather}F_{red}^j={F}_I^{r_M-\lambda_{j-1}+j}={F}_I^{r_M-\lambda_{j-1}+j+1}=\ldots={F}_I^{r_M-\lambda_j+j}.\label{redhod}\end{gather}

{\bf Conjecture:} (Schneider \cite{schn}) (a) For any $K[\Gamma]$-module $M$ we have \begin{gather}H^{d}(X_{\Gamma},{{M}}^{\Gamma}\otimes_K\Omega^{\bullet}_{{X}_{\Gamma}}    )=F_H^{i+1}\bigoplus F_{\Gamma}^{d-i}\quad\quad(0\le i\le d-1).\label{spliteins}\end{gather} (b) For any irreducible $K$-rational representation $M$ we have \begin{gather}H^{d}(X_{\Gamma},{{M}}^{\Gamma}\otimes_K\Omega^{\bullet}_{{X}_{\Gamma}}    )=F_{red}^{i+1}\bigoplus F_{\Gamma}^{d-i}\quad\quad(0\le i\le d-1).\label{spliteinsred}\end{gather}(c) For all $M$ as in (b) we have $F_{red}^j=F_H^{j}$ for all $j$.\\(d) For all $M$ as in (b) the spectral sequence (\ref{redss}) degenerates in $E_1$.

\section{Holomorphic discrete series representations}
\label{holdisrep}

Let $L_1$ denote the algebraic subgroup of ${\rm GL}_{d+1}$ consisting of matrices $(a_{ij})_{0\le i, j\le d}\in{\rm GL}_{d+1}$ with $a_{00}=1$ and $a_{ij}=0$ if {\it either} $i=0$ {\it or} $j=0$. Thus $L_1\cong{\rm GL}_{d}$. The irreducible $K$-rational representations of $L_1(K)$ are again parametrized by their highest weight. The set of possible highest weights $\mu$ is in bijection with the set of integer valued $d$-tuples $(\mu_1,\ldots,\mu_d)$ satisfying $\mu_1\ge\ldots\ge\mu_d$. We simply write $\mu=(\mu_1,\ldots,\mu_d)$. If $\lambda=(\lambda_0,\ldots,\lambda_d)$ (with $\lambda_0\ge\ldots\ge\lambda_d$) is the highest weight of a $K$-rational $G$-representation we set $$\mu(\lambda(j))=(\lambda_0-\lambda_j+j+1,\ldots,\lambda_{j-1}-\lambda_j+j+1,\lambda_{j+1}-\lambda_j+j,\ldots,\lambda_{d}-\lambda_j+j)$$for $0\le j\le d$. Note that conversely, given a $d$-tuple $\mu=(\mu_1,\ldots,\mu_d)$ satisfying $\mu_1\ge\ldots\ge\mu_d$, we have $\mu=\mu(\lambda(j))$ for some $\lambda$ if and only if \begin{gather}\mu_s\ne s\quad\quad{\mbox{ for all }}\,\,1\le s\le d.\label{deraoc}\end{gather}Indeed, if this condition holds then the unique $j$ which works is the largest $j\in\{0,\ldots,d\}$ for which $\mu_j\ge j+1$ (where we set $\mu_0=+\infty$). Also note that in this case there exists exactly one $\lambda$ with $\lambda_d=0$ and $\mu=\mu(\lambda(j))$.
 
Schneider \cite{schn} (described there in different normalizations) associated to any irreducible $K$-rational representation $V$ of $L_1(K)$ of highest weight $\mu$ a ${\rm SL}_{d+1}(K)$-equivariant vector bundle $V(\mu)$ on $X$. The corresponding ${\rm SL}_{d+1}(K)$-representations on its space of global sections $\Gamma(X,V(\mu))$ have recently been studied by Orlik \cite{orlik}, and it is of great interest to understand how they contribute, via Proposition \ref{hoder} below, to the de Rham cochomology of $X$ and its quotients $X_{\Gamma}$. We want to explain why the conjecture (d) at the end of section \ref{conj} predicts that $\Gamma$-group cohomology of $\Gamma(X,V(\mu))$ occurs only in a single degree, except possibly for those $\mu$ not satisfying (\ref{deraoc}).

Let again $M$ be an irreducible $K$-rational representation of $G={\rm GL}\sb {d+1}(K)$ of highest weight $(\lambda_0\ge\lambda_1\ge\ldots\ge\lambda_d)$ (with $\lambda_d=0$).

\begin{pro}\label{hoder} (Schneider) For $0\le j\le d$ let ${D}^j(M)$ be as in \ref{peterco}. Then there exists an isomorphism of ${\rm SL}_{d+1}(K)$-equivariant vector bundles ${D}^j(M)\cong V(\mu(\lambda(j)))$.
\end{pro}

\begin{kor}\label{degcon} Suppose $M\ne K$ and let $0\le s\le d$. If the spectral sequence (\ref{redss}) degenerates in $E_1$ then we have $$H^t({\Gamma},\Gamma(X,V(\mu(\lambda(s)))))=0\quad\quad{\mbox{ whenever }}\,\,t\ne d-s.$$If moreover \begin{gather}H^{d}(X_{\Gamma},{{M}}^{\Gamma}\otimes_K\Omega^{\bullet}_{{X}_{\Gamma}}    )=F_{red}^{i+1}\oplus F_{\Gamma}^{d-i}\label{reddec}\end{gather} for all $0\le i\le d-1$ then we have$$\dim_K(H^{d-s}({\Gamma},\Gamma(X,V(\mu(\lambda(s))))))=\mu(\Gamma,M).$$
\end{kor}

{\sc Proof:} It follows from Proposition \ref{hoder} (and the fact that $X$ is a Stein space) that the spectral sequence (\ref{redss}) is isomorphic with\begin{gather}E_1^{st}=H^t({\Gamma},\Gamma(X,V(\mu(\lambda(s)))))\Rightarrow H^{s+t}(X_{\Gamma},{D}^{\bullet}(M)^{\Gamma})=H^{s+t}(X_{\Gamma},M^{\Gamma}\otimes_K\Omega^{\bullet}_{{X}_{\Gamma}} ).\label{redgamss}\end{gather}Now the first statement follows from formula (\ref{othdeg}) and the second statement then follows from formula (\ref{pecomp}).\hfill$\Box$\\

\section{The Monodromy Operator}

\label{monosection}

Let $M$ be a $K[\Gamma]$-module as before. For $0\le i\le d$ let $\widehat{X}^i$ denote the set of pointed $i$-simplices $(\sigma,v)$ of the Bruhat-Tits building of ${\rm PGL}_{d+1}/K$ (i.e. $\sigma$ is an $i$-simplex and $v$ is one of its vertices). We have a natural action of $G$ on $\widehat{X}^i$. Let ${\mathcal C}^i(M)$ denote the $K$-vector space of all functions $f$ on $\widehat{X}^i$ with values in $M$ satisfying the following condition: for any $(\sigma,v)\in\widehat{X}^i$ and any $g\in G$ with $g\sigma=\sigma$ we have $$f((\sigma,v))=(-1)^{d\omega({\rm det}(g))}f((\sigma,gv));$$ here $\omega:K^{\times}\to{\mathbb Z}$ denotes the valuation normalized by $\omega(\pi)=1$. Using alternating signs we may form the $\Gamma$-equivariant complex$${\mathcal C}^{\bullet}(M)=[{\mathcal C}^0(M)\stackrel{{\partial}_M^{0}}{\longrightarrow}{\mathcal C}^1(M)\stackrel{{\partial}_M^{1}}{\longrightarrow}\ldots\stackrel{{\partial}_M^{d-1}}{\longrightarrow}{\mathcal C}^d(M)\longrightarrow0].$$From \cite{gar} p.387 it follows that ${\mathcal C}^{\bullet}(M)$ is a $K[\Gamma]$-injective resolution of $M$. In particular this implies\begin{gather}H^d(\Gamma,M)=\frac{{\mathcal C}^{d}(M)^{\Gamma}}{\bi[{\mathcal C}^{d-1}(M)^{\Gamma}\stackrel{{\partial}_M^{d-1}}{\longrightarrow}{\mathcal C}^d(M)]}.\label{gammacoc}\end{gather}

Let ${\mathcal C}_{har}^{d}(M)\subset{\mathcal C}^{d}(M)$ denote the subspace of harmonic $d$-cochains with values in $M$. Thus, $f\in{\mathcal C}^{d}(M)$ belongs to ${\mathcal C}_{har}^{d}(M)$ if and only if for any $(\tau,v)\in\widehat{X}^{d-1}$ we have$$\sum_{\sigma\supset\tau}f((\sigma,v))=0$$where the sum runs over all $d$-simplices $\sigma$ containing $\tau$.    

\begin{lem}\label{garland} (Garland) Assume $d\ge2$. We have a direct sum decomposition$${\mathcal C}^{d}(M)^{\Gamma}={\mathcal C}^d_{har}(M)^{\Gamma}\bigoplus \bi[{\mathcal C}^{d-1}(M)^{\Gamma}\stackrel{{\partial}_M^{d-1}}{\longrightarrow}{\mathcal C}^d(M)].$$
\end{lem}

{\sc Proof:} Of course, to prove our Lemma we are allowed to extend scalars from $K$ to ${\mathbb C}$ (for some embedding $K\to{\mathbb C}$). By a result of Margulis \cite{mar} IX Theorem (5.12) (i) (here the hypothesis $d\ge2$ is used), $M\otimes_K{\mathbb C}$ carries a $\Gamma$-invariant positive definite hermitian form, i.e. $M\otimes_K{\mathbb C}$ is a unitary $\Gamma$-module. Therefore our claim is a special case of \cite{gar} Proposition 3.16.\hfill$\Box$\\

{\it Remark:} Strictly speaking, \cite{gar} Proposition 3.16 is formulated for a $\Gamma$-representation on an ${\mathbb R}$-vector space carrying a $\Gamma$-invariant positive-definite form, but the proof translates verbatim to the unitary case. Moreover, as pointed out in \cite{gar}, the existence of a $\Gamma$-invariant inner product on $M\otimes_K{\mathbb C}$ can be verified in many important examples circumventing Margulis' theorem. For example this applies in the context of section \ref{unifsect}: in that case $\Gamma\subset {S}'({\mathbb Q})$ (even $\Gamma\subset {S}'({\mathbb Z}[p^{-1}])$) for an inner form $S'$ of ${\rm SL}_{d+1}$ defined over ${\mathbb Q}$ such that ${S}'({\mathbb R})$ is compact and such that the embedding $\Gamma\to{\rm SL}_{d+1}(K)$ is induced by an isomorphism $S'(K)\cong{\rm SL}_{d+1}(K)$. This isomorphism allows us to view a $K$-rational representation $M$ of ${\rm SL}_{d+1}(K)$ as a $K$-rational representation of $S'(K)$. Thus by \cite{gar} p.414 there is a $\Gamma$-invariant inner product on $M\otimes_K{\mathbb C}$.\\

Consider the $\Gamma$-equivariant residue map$$\widetilde{res}: M\otimes_K H^d_{dR}(X)\longrightarrow {\mathcal C}^{d}(M)$$introduced in \cite{schtei}. It is constructed like the residue map of section \ref{resisec}, as follows. Let $(\sigma,v)\in \widehat{X}^d$. Restrict a class $\eta\in M\otimes_K H^d_{dR}(X)$ to $]\sigma^0[$, with $]\sigma^0[$ denoting the preimage of the interior of $\sigma$ under the residue map from $X$ to the Bruhat-Tits building of ${\rm PGL}_{d+1}/K$; equivalently, $\sigma$ corresponds to the (one-point) intersection of $d+1$ irreducible components of ${\mathfrak X}\otimes k$, and $]\sigma^0[$ is the preimage of this intersection under the specialization map $X\to{\mathfrak X}\otimes k$. Thus $]\sigma^0[$ is a generalized rigid analytic annullus and it makes sense to expand $\eta|_{]\sigma^0[}$ into a convergent Laurent series and to extract the coefficient of its non-integrable term, i.e. the residue of $\eta|_{]\sigma^0[}$. We define the value of $\widetilde{res}(\eta)\in{\mathcal C}^{d}(M)$ at $(\sigma,v)$ to be this residue of $\eta|_{]\sigma^0[}$, where the {\it sign} is fixed by the pointing $v$ of $(\sigma,v)$. Namely, as any simplex in the (oriented) Bruhat-Tits building of ${\rm PGL}_{d+1}/K$ carries a natural {\it cyclic} ordering of its vertices, to fix a pointing means fixing a total ordering of its vertices. Identifying $\sigma$ with the (one-point) intersection of $d+1$ irreducible components of ${\mathfrak X}\otimes k$ as indicated, this means fixing a total ordering of these $d+1$ irreducible components. This fixes the sign of the residue maps, see section \ref{resisec}.

By the main result of \cite{schtei} the map $\widetilde{res}$ takes in fact values in ${\mathcal C}^d_{har}(M)$ and induces an isomorphism $\widetilde{res}:M\otimes_K H^d_{dR}(X)\cong {\mathcal C}^d_{har}(M)$. In particular it induces an isomorphism$$\widetilde{res}^{\Gamma}:H^0(\Gamma,M\otimes_K H^d_{dR}(X))\cong H^0(\Gamma,{\mathcal C}^d_{har}(M)).$$Composing with the isomorphism $$H^0(\Gamma,{\mathcal C}^d_{har}(M))\cong H^d(\Gamma,M)$$ which we get by combining formula (\ref{gammacoc}) with Lemma \ref{garland}, we get an isomorphism $$res_{\Gamma}:H^0(\Gamma,M\otimes_K H^d_{dR}(X))\cong  H^d(\Gamma,M)=H^d(\Gamma,M\otimes_K H^0_{dR}(X)).$$Let$$\gamma:H^d(\Gamma,M\otimes_K H^0_{dR}(X))\longrightarrow H^{d}(X_{\Gamma},{{M}}^{\Gamma}\otimes_K\Omega^{\bullet}_{{X}_{\Gamma}}),$$$$\delta:H^{d}(X_{\Gamma},{{M}}^{\Gamma}\otimes_K\Omega^{\bullet}_{{X}_{\Gamma}})\longrightarrow H^0(\Gamma,M\otimes_K H^d_{dR}(X))$$denote the boundary maps in the spectral sequence (\ref{covss}).

By comparison isomorphisms with log convergent (log rigid, log crystalline) cohomology we get a nilpotent monodromy operator $N$ on $H^{d}(X_{\Gamma},{{M}}^{\Gamma}\otimes_K\Omega^{\bullet}_{{X}_{\Gamma}})$ (cf.\cite{hk} and section \ref{resisec}). Let $N^d$ denote its $d$-fold iterate.

\begin{satz}\label{resmono} At least up to sign we have, as endomorphisms of $H^{d}(X_{\Gamma},{{M}}^{\Gamma}\otimes_K\Omega^{\bullet}_{{X}_{\Gamma}})$:$$N^d=\gamma\circ res_{\Gamma}\circ \delta.$$
\end{satz}

{\sc Proof:} This follows from Theorem \ref{nviares}, applied to the proper semistable formal ${\mathcal O}_K$-scheme ${\mathfrak X}_{\Gamma}$ (i.e. what is denoted by ${\mathfrak X}$ in section \ref{resisec} is our ${\mathfrak X}_{\Gamma}$ here). The point is that the spectral sequence (\ref{drspec}), whose boundary maps enter into the statement of Theorem \ref{nviares}, is the $E_2$-version of the $E_1$-spectral sequence (\ref{covss}), whose boundary maps enter into the statement of Theorem \ref{resmono}. This can be deduced from the arguments from \cite{latt} section 5, but much less is needed for our purposes here, as follows. Let $Y={\mathfrak X}_{\Gamma}\otimes k$ and form $Y^{d+1}$ as in section \ref{resisec}. Let $F=M^{\Gamma}$, the locally constant sheaf on $Y$ obtained by descent of the constant sheaf with value $M$ on ${\mathfrak X}\otimes k$. Let $\alpha$, $\beta$ and $Res$ be as in Theorem \ref{nviares}, so we need to show $\gamma\circ res_{\Gamma}\circ \delta=\alpha\circ Res\circ \beta$. We have a natural map (up to once fixing a global sign)$$\kappa:{\mathcal C}^d(M)^{\Gamma}\longrightarrow H_{dR}^0(]Y^{d+1}[,F).$$The description of the residue map in either context immediately tells us $\kappa\circ res_{\Gamma}\circ \delta=Res\circ \beta$. On the other hand, both $\alpha$ and $\gamma$ are boundary morphisms induced by the inclusion of sheaf complexes $F\to\Omega^{\bullet}_{{X}_{\Gamma}}\otimes F$. Our claim follows.\hfill$\Box$\\    

Let $\Phi:M\to M$ be a $K[\Gamma]$-linear automorphism such that its $m$-th iterate $\Phi^m$ (some $m\in{\mathbb N}$) is multiplication by $\alpha\in K$. The induced automorphism $\Phi:M^{\Gamma}\to M^{\Gamma}$ of the descended sheaf $M^{\Gamma}$ on $X_{\Gamma}$ provides $H^{d}(X_{\Gamma},{{M}}^{\Gamma}\otimes_K\Omega^{\bullet}_{{X}_{\Gamma}})$ with a $K$-linear Frobenius automorphism which we again denote by $\Phi$ and which satisfies $N\Phi=p^f\Phi N$ (with $p^f=|k|$). This comes again by comparison with log convergent (log rigid, log crystalline) cohomology, see section \ref{resisec}. By \cite{hk} Theorem 7.3 this endomorphism $\Phi$ of $H^{d}(X_{\Gamma},{{M}}^{\Gamma}\otimes_K\Omega^{\bullet}_{{X}_{\Gamma}})$ respects the filtration $F_{\Gamma}^{\bullet}$ and $\Phi^m$ acts on the graded piece $F^r_{\Gamma}/F^{r+1}_{\Gamma}$ by multiplication with $p^{f(d-r)m}\alpha$. (In \cite{hk} Theorem 7.3 this is formulated in the case where $\Phi:M\to M$ is the identity, but the easy proof clearly applies in our more general situation). In particular, $F_{\Gamma}^{\bullet}$ is also the ($p$-adic) {\it slope-filtration}: it is characterized be the property that the eigenvalues of $\Phi$ on each of the graded pieces have the same $p$-adic size, and different $p$-adic sizes on different pieces. Similarly, $F_{\Gamma}^{\bullet}$ is characterized in the same way by the archimedean sizes (for all complex embeddings of $K$) of the eigenvalues of $\Phi$. In particular, if $\alpha$ is a Weil number then the graded pieces are pure of some weight and $F_{\Gamma}^{\bullet}$ is a {\it weight filtration}. 

{\it Remark:} In our later applications $\Phi$ indeed satisfies $\Phi^m=\alpha$ for some $m$ and $\alpha$, but here we could more generally just require that $\Phi$ be $p$-adically isoclinic, resp. be pure of some weight.

\begin{satz}\label{moweeas} Let $d\ge2$. The monodromy weight conjecture holds true for $H^d({X}_{\Gamma},M^{\Gamma}\otimes_K\Omega_{X^{\Gamma}}^{\bullet})$. More precisely, if $M$ contains no $\Gamma$-invariant vector, then$$F_{\Gamma}^r=\bi(N^{r})=\ke(N^{d+1-r})$$in $H^d({X}_{\Gamma},M^{\Gamma}\otimes_K\Omega_{X^{\Gamma}}^{\bullet})$ for all $r\ge0$.
\end{satz}

{\sc Proof:} For $M=K$ trivial this has been proved independently in \cite{desh} and \cite{ito}. Therefore we assume that $M$ contains no $\Gamma$-invariant vector (otherwise the argument must be slightly modified). The degeneration of the spectral sequence (\ref{covss}) implies that $\delta$ is surjective and $\gamma$ is injective. Moreover, by \cite{latt} formula (35), we know that the image of $\gamma$ is just $F_{\Gamma}^d$, whose dimension was recalled in formula (\ref{pecomp}): it is $\mu(\Gamma,M)$, where, recall,$$(d+1)\mu(\Gamma,M)=\dim_KH^d(X_{\Gamma},M^{\Gamma}\otimes_K\Omega_{X^{\Gamma}}^{\bullet}).$$As $res_{\Gamma}$ is bijective we get with Theorem \ref{resmono} that the image of $N^d$ has dimension $\mu(\Gamma,M)$. [In fact, it is enough to have Theorem \ref{resmono} (and Lemma \ref{garland}) only in the case where the $K[\Gamma]$-module $M$ is irreducible; the proof of $\dim_K(N^d)=\mu(\Gamma,M)$ for general $M$ follows from this case by considering a Jordan H\"older series of $M$.] But we also know $N^{d+1}=0$ [either from the geometry underlying $N$, or just (as we will see in a second) by $\bi(N^{d+1})\subset F_{\Gamma}^{d+1}$ which is $0$]. Together it follows from elementary linear algebra that more generally we have $$\dim_K(\bi(N^{r}))=\dim_K(\ke(N^{d+1-r}))=(d+1-r)\mu(\Gamma,M).$$ On the other hand, as $F_{\Gamma}^{\bullet}$ is the slope (and weight) filtration with $\Phi$-action on the graded pieces as described and as $N\Phi=p^f\Phi N$, we have $\bi(N^{r})\subset F_{\Gamma}^r\subset\ke(N^{d+1-r})$. Together the claim follows.\hfill$\Box$

\section{The convergent $F$-isocrystals ${\mathfrak E}_{M}$}
\label{filconfiso}

For a vector space $V$ over a field $E$ with structure morphism $\nu:E{\to}{\rm End}_{\mathbb Z}(V)$ and a field automorphism $\gamma$ of $E$ we write $V_{\gamma}$ to denote the same additive group $V$ but with the structure of $E$-vector space given by $\nu\circ\gamma$. If a group $H$ acts on $V$ through $E$-linear automorphisms then it also acts on $V_{\gamma}$ through $E$-linear automorphisms. We write $V_{\gamma}$ also for the constant sheaf generated by $V_{\gamma}$ on various schemes. \\ 

We recall some facts from \cite{hk}, cf. also section \ref{resisec}. Similarly to the spectral sequence (\ref{covss}) there is for every $K^t[\Gamma]$-module $N$ (with $\dim_{K^t}N<\infty$) with associated convergent isocrystal $N\otimes_{{\mathcal O}_{K^t}}{\mathcal O}^{conv}_{{\mathfrak {X}}_0}$ on ${\mathfrak {X}}_0$ a covering spectral sequence\begin{gather}E_2^{r,s}=H^r(\Gamma,H_{conv}^{s}({\mathfrak X}_{0},N\otimes_{{\mathcal O}_{K^t}}{\mathcal O}^{conv}_{{\mathfrak {X}}_0}))\Rightarrow H_{conv}^{r+s}({\mathfrak X}_{{\Gamma},0},(N\otimes_{{\mathcal O}_{K^t}}{\mathcal O}^{conv}_{{\mathfrak {X}}_0})^{\Gamma}).\label{covspss}\end{gather}which induces a filtration\begin{gather}H^d({\mathfrak X}_{{\Gamma},0},(N\otimes_{{\mathcal O}_{K^t}}{\mathcal O}^{conv}_{{\mathfrak {X}}_0})^{\Gamma})=S^0\supset S^1\supset\ldots\supset S^{d+1}=0.\label{gamffil}\end{gather} Here $H_{conv}$ denotes logarithmic convergent cohomology with respect to the logarithmic enlargement $T=(\spf({\mathcal O}_{K^t}),1\mapsto 0)$ (the 'Hyodo-Kato base') of the log scheme $(\spec(k),1\mapsto 0)$ over which ${\mathfrak X}_{{\Gamma},0}$ is a log smooth log scheme. Thus $H_{conv}^{*}({\mathfrak X}_{{\Gamma},0},(N\otimes_{{\mathcal O}_{K^t}}{\mathcal O}^{conv}_{{\mathfrak {X}}_0})^{\Gamma})$ is a finite dimensional $K^t$-vector space endowed with a nilpotent $K^t$-linear endomorphism $N$. If the convergent isocrystal $N\otimes_{{\mathcal O}_{K^t}}{\mathcal O}^{conv}_{{\mathfrak {X}}_0}$ carries a $\sigma$-linear Frobenius structure then also $H_{conv}^{*}({\mathfrak X}_{{\Gamma},0},(N\otimes_{{\mathcal O}_{K^t}}{\mathcal O}^{conv}_{{\mathfrak {X}}_0})^{\Gamma})$ receives a $\sigma$-linear Frobenius endomorphism $\varphi$ such that $N\varphi=p\varphi N$. The $f$-fold iterate of the ($p$-power) Frobenius action on the $s$-th convergent cohomology group with constant coefficients $H_{conv}^{*}({\mathfrak X}_{0},{\mathcal O}^{conv}_{{\mathfrak {X}}_0})$ of the reduction ${\mathfrak X}_{0}$ of ${\mathfrak X}$ is just multiplication with $p^{fs}$, see \cite{hk}. Hence if $N$ is endowed with a $\sigma$-linear endomorphism $\Phi_N$ such that its $mf$-fold iterate for some $m\in\mathbb{N}$ is multiplication with some $\alpha\in K^t$, and if $N\otimes_{{\mathcal O}_{K^t}}{\mathcal O}^{conv}_{{\mathfrak {X}}_0}$ is endowed with the Frobenius structure $\Phi_N\otimes\Phi$, then the $mf$-fold iterate of Frobenius on $H_{conv}^{s}({\mathfrak X}_{0},N\otimes_{{\mathcal O}_{K^t}}{\mathcal O}^{conv}_{{\mathfrak {X}}_0})$ is multiplication with $\alpha p^{mfs}$. If $\Phi_N$ commutes with the $\Gamma$-action then also $(N\otimes_{{\mathcal O}_{K^t}}{\mathcal O}^{conv}_{{\mathfrak {X}}_0})^{\Gamma}$ and hence $H^d({\mathfrak X}_{{\Gamma},0},(N\otimes_{{\mathcal O}_{K^t}}{\mathcal O}^{conv}_{{\mathfrak {X}}_0})^{\Gamma})$ receives a Frobenius action and by what we just said, $S^{\bullet}$ is Frobenius stable and the $mf$-fold iterate of Frobenius acts on $S^s/S^{s+1}$ as multiplication with $\alpha p^{mf(d-s)}$. Moreover we have a $K$-linear Hyodo-Kato isomorphism\begin{gather}H^{d}({X}_{{\Gamma}},N^{\Gamma}\otimes_{K^t}\Omega^{\bullet}_{{X}_{\Gamma}})\cong H_{conv}^{s}({\mathfrak X}_{\Gamma,0},(N\otimes_{{\mathcal O}_{K^t}}{\mathcal O}^{conv}_{{\mathfrak {X}}_0})^{\Gamma})\otimes_{K^t}{K}\label{modelhyk}\end{gather}such that the filtration (\ref{gamfil}) of the left hand side corresponds to the filtration $S^{\bullet}\otimes_{K^t}{K}$ of the right hand side. (This isomorphism depends on the choice of a uniformizer in ${K}$; we fix such a choice once and for all, e.g. our $\pi$.) In particular, the computations (\ref{pecomp}) give the $K^t$-vector space dimensions of the $S^s$.\\

{\it Remark:} If $N$ even carries a $K[\Gamma]$-structure (extending its $K^t[\Gamma]$-structure) then the decomposition $N\otimes_{K^t}K=\oplus_{\kappa}N_{\kappa}$, with $\kappa$ running through ${\rm Gal}(K/K^t)$, allows us to view the isomorphism (\ref{modelhyk}) as the direct sum, over all $\kappa$, of $K$-linear isomorphisms$$H^{d}({X}_{{\Gamma}},N_{\kappa}^{\Gamma}\otimes_{K}\Omega^{\bullet}_{{X}_{\Gamma}})\cong H_{conv}^{s}({\mathfrak X}_{\Gamma,0},(N_{\kappa}\otimes_{{\mathcal O}_{K^t}}{\mathcal O}^{conv}_{{\mathfrak {X}}_0})^{\Gamma}).$$Here we have no $K^t$-structure on the right hand side. (This point of view has often been adopted in \cite{hk}, as well as in section \ref{resisec}.)\\  

Let $M$ be an irreducible $K$-rational representation of $G$ of highest weight $(\lambda_0\ge\lambda_1\ge\ldots\ge\lambda_d)$ (with $\lambda_d=0$). Consider the $K^t$-vector space$${\bf M}_{M}=\bigoplus_{j\in{\mathbb Z}/f}\bigoplus_{i\in{\mathbb Z}/(d+1)}M_{\sigma^{j}}.$$We let $G$ act on ${\bf M}_{M}$ by the action on the direct summands. We use the ${\mathbb Q}_p$-linear identities $M_{\sigma^{j}}=M_{\sigma^{j-1}}$ to endow ${\bf M}_{M}$ with a $\sigma$-linear operator $\Phi_{{\bf M}_{M}}$ as follows: it sends an element $m$ from the $(j,i)$ direct summand\\--- to the element $p^{r_M}m$ in the $(j-1,i)$ direct summand if $j\ne0\in{\mathbb Z}/f$,\\--- to the element $p^{r_M}m$ in the $(j-1,i+1)$ direct summand if $j=0\in{\mathbb Z}/f$ and $i\ne d\in{\mathbb Z}/(d+1)$, and\\--- to the element $p^{r_M}\pi^{-r_M}m$ in the $(j-1,i+1)$ direct summand if $j=0\in{\mathbb Z}/f$ and $i=d\in{\mathbb Z}/(d+1)$. (Here we use the original $K$-vector space structure of the $(j-1,i+1)$ direct summand $M_{\sigma^{j}}=M$ to define its element $p^{r_M}\pi^{-r_M}m$.)\\Thus the $n(d+1)$-fold iterate of $\Phi_{{\bf M}_{M}}$ is multiplication with an element of $K^t$ of $p$-adic order $(n(d+1)-1)r_M$ (if $\pi^e=p$ for $e=n/f=[K:K^t]$ then this is the element $p^{(n(d+1)-1)r_M}$). In particular, $\Phi_{{\bf M}_{M}}$ on ${\bf M}_{M}$ is isoclinic of slope $\frac{(n(d+1)-1)r_M}{n(d+1)}$. 

We associate to $M$ a $G$-equivariant filtered convergent $F$-isocrystal ${\mathfrak E}_{M}$ on ${\mathfrak{X}}$ as follows. Let ${\mathcal O}^{conv}_{{\mathfrak {X}}_0}$ denote the structure sheaf of the convergent site of the special fibre ${\mathfrak X}_0$ of ${\mathfrak{X}}$: it assigns to each enlargement $Z$ of ${\mathfrak X}_0$ (i.e. $Z$ is a $p$-torsion free locally noetherian formal ${\mathcal O}_{K^t}$-scheme together with a $k$-morphism of a subscheme of definition of $Z$ to ${\mathfrak X}_0$, see \cite{and}, \cite{faalet}) the structure sheaf ${\mathcal O}_Z$ of $Z$. As a convergent isocrystal we set ${\mathfrak E}_{M}={\bf M}_{M}\otimes_{{\mathcal O}_{K^t}}{\mathcal O}^{conv}_{{\mathfrak {X}}_0}$. Thus the value ${\mathfrak E}_{M}(Z)$ of ${\mathfrak E}_{M}$ at an enlargement $Z$ of ${\mathfrak X}_0$ is given by ${\bf M}_{M}\otimes_{{\mathcal O}_{K^t}}{\mathcal O}_Z$. The Frobenius structure on ${\mathfrak E}_{M}$ is felt only at enlargements $Z$ of ${\mathfrak X}_0$ with a lifting $\Phi_Z:{\mathcal O}_Z\to {\mathcal O}_Z$ of the $p$-power Frobenius of the reduction: there it is given by $\Phi_{{\bf M}_{M}}\otimes\Phi_Z$. Finally, the filtration of ${\mathfrak E}_{M}$ is felt only at the evaluation ${\mathfrak E}_{M}({\mathfrak X})$ of ${\mathfrak E}_{M}$ at the enlargement ${\mathfrak X}$ of ${\mathfrak X}_0$. The $K$-vector space structure of $M$ provides a decomposition\begin{align}{\mathfrak E}_{M}({\mathfrak X})&=\bigoplus_{j\in{\mathbb Z}/f}\bigoplus_{i\in{\mathbb Z}/(d+1)}M_{\sigma^{j}}\otimes_{{\mathcal O}_{K^t}}{\mathcal O}_{{\mathfrak X}}\notag \\{}&=\bigoplus_{j\in{\mathbb Z}/f}\bigoplus_{i\in{\mathbb Z}/(d+1)}\bigoplus_{\delta\in{\rm Gal}(K/{\mathbb Q}_p)_{\sigma^{j}}}M_{\delta}\otimes_{{\mathcal O}_{K}}{\mathcal O}_{{\mathfrak X}}\notag\\{}&=\bigoplus_{i\in{\mathbb Z}/(d+1)}\bigoplus_{\delta\in{\rm Gal}(K/{\mathbb Q}_p)}M_{\delta}\otimes_{{\mathcal O}_{K}}{\mathcal O}_{{\mathfrak X}}\notag\end{align}Here we wrote ${\rm Gal}(K/{\mathbb Q}_p)_{\sigma^{j}}=\{\delta \in {\rm Gal}(K/{\mathbb Q}_p);\,\,\delta|_{K^t}=\sigma^{j}\}$. We filter the direct summands $M_{\delta}\otimes_{{\mathcal O}_{K}}{\mathcal O}_{{\mathfrak X}}$ as follows. If $\delta=\id_K$ we let$${f}^s(M\otimes_{{\mathcal O}_{K}}{\mathcal O}_{{\mathfrak X}})={\mathcal O}_{{\mathfrak X}}.\overline{u}(z)({f}^sM)$$for $s\in{\mathbb Z}$ (like in definition (\ref{verkl})).  If $\delta\ne\id_K$ and $s\le r_M$ we let$${f}^{s}(M_{\delta}\otimes_{{\mathcal O}_{K}}{\mathcal O}_{{\mathfrak X}})=M_{\delta}\otimes_{{\mathcal O}_{K}}{\mathcal O}_{{\mathfrak X}}$$and if $s\ge r_M+1$ we let ${f}^{s}(M_{\delta}\otimes_{{\mathcal O}_{K}}{\mathcal O}_{{\mathfrak X}})=0$. Together we get a filtration ${f}^{\bullet}{\mathfrak E}_{M}({\mathfrak X})$ of ${\mathfrak E}_{M}({\mathfrak X})$. (See the remark following Proposition \ref{iotas} below for why this at first sight strange looking filtration (in particular in its separation between the $\delta\ne{\rm id}_K$- and the $\delta={\rm id}_K$-direct summands) is in fact very natural.)

Now let $\Gamma<{\rm PGL}_{d+1}(K)$ be as before. The $G$-equivariant filtered convergent $F$-isocrystal ${\mathfrak E}_{M}$ descends to a filtered convergent $F$-isocrystal ${\mathfrak E}_{M}^{\Gamma}$ on ${{\mathfrak X}}_{\Gamma}$ (by taking $\Gamma$-invariants). We may view the evaluation ${\mathfrak E}_{M}^{\Gamma}({{\mathfrak X}}_{\Gamma})$ of ${\mathfrak E}_{M}^{\Gamma}$ at ${{\mathfrak X}}_{\Gamma}$ as a filtered ${\mathcal O}_{{{X}}_{\Gamma}}$-module with connection. Denoting the filtration of ${\mathfrak E}_{M}^{\Gamma}({{\mathfrak X}}_{\Gamma})$ by ${f}^{\bullet}{\mathfrak E}_{M}^{\Gamma}({{\mathfrak X}}_{\Gamma})$ we filter the de Rham complex ${\mathfrak E}_{M}^{\Gamma}({{\mathfrak X}}_{\Gamma})\otimes\Omega^{\bullet}_{{X}_{\Gamma}}$ by setting\begin{gather}{\mathcal F}^j_M=[({f}^{j}{\mathfrak E}_{M}^{\Gamma}({{\mathfrak X}}_{\Gamma}))\otimes\Omega^{0}_{{X}_{\Gamma}}\to({f}^{j-1}{\mathfrak E}_{M}^{\Gamma}({{\mathfrak X}}_{\Gamma}))\otimes\Omega^{1}_{{X}_{\Gamma}}\to({f}^{j-2}{\mathfrak E}_{M}^{\Gamma}({{\mathfrak X}}_{\Gamma}))\otimes\Omega^{2}_{{X}_{\Gamma}}\to\ldots].\label{drcof}\end{gather}This filtration gives rise to the spectral sequence\begin{gather}E_1^{st}=H^{s+t}({X}_{\Gamma},{\mathcal F}_M^{s}/{\mathcal F}_M^{s+1})\Rightarrow H^{s+t}({X}_{\Gamma},{\mathfrak E}_{M}^{\Gamma}({{\mathfrak X}}_{\Gamma})\otimes\Omega^{\bullet}_{{X}_{\Gamma}}).\label{tilhoa}\end{gather}

\begin{pro}\label{degeq} The following (i) and (ii) are equivalent:\\(i) the spectral sequence (\ref{tilhoa}) degenerates in $E_1$\\(ii) the spectral sequences (\ref{tilhob}) and (\ref{redss}) for $M$, and the spectral sequence (\ref{hodss}) for $M_{\delta}$, for any ${\delta\in{\rm Gal}(K/{\mathbb Q}_p)}$, $\delta\ne \id_K$, degenerate in $E_1$. 
\end{pro}

{\sc Proof:} By construction, ${\mathfrak E}_{M}({{\mathfrak X}})$ decomposes into a direct sum of $G$-equivariant filtered ${\mathcal O}_{{X}}$-modules with connection of the form $M_{\delta}\otimes_{K}{\mathcal O}_{{X}}$ with ${\delta\in{\rm Gal}(K/{\mathbb Q}_p)}$, each such copy for any fixed $\delta$ occuring exactly $(d+1)$ times. Taking $\Gamma$-invariants we get direct summands $(M_{\delta}\otimes_{K}{\mathcal O}_{{X}})^{\Gamma}$ of ${\mathfrak E}_{M}^{\Gamma}({{\mathfrak X}}_{\Gamma})$ as filtered modules with connection. It follows from our definitions that for $\delta=\id_K$ the associated spectral sequence computing its de Rham cohomology on $X_{\Gamma}$ is precicely the spectral sequence (\ref{tilhob}), which is isomorphic with spectral sequence (\ref{redss}). If $\delta\ne \id_K$ the associated spectral sequence computing its de Rham cohomology on $X_{\Gamma}$ is the spectral sequence (\ref{hodss}) with the filtration numbering shifted by $r_M$. Since these spectral sequences are direct summands of the spectral sequence (\ref{tilhoa}), the degeneration of the latter implies the degeneration of the former.\hfill$\Box$\\

Let $$D=H_{conv}^{d}({\mathfrak X}_{{\Gamma},0},{\mathfrak E}_{M}^{\Gamma}),$$$$D_{{K}}=H_{dR}^{d}({X}_{{\Gamma}},{\mathfrak E}_{M}^{\Gamma}({\mathfrak X}_{{\Gamma}}))=H^{d}({X}_{\Gamma},{\mathfrak E}_{M}^{\Gamma}({{\mathfrak X}}_{\Gamma})\otimes\Omega^{\bullet}_{{X}_{\Gamma}}).$$Thus $D$ is a finite dimensional $K^t$-vector space with $K^t$-linear nilpotent endomorphism $N$ and $\sigma$-linear endomorphism $\varphi$, while $D_K$ is a finite dimensional ${K}$-vector space endowed with a decreasing filtration $F^{\bullet}$. We will regard the ${K}$-vector space isomorphism\begin{gather}D_{{K}}\cong D\otimes_{K^t}{K}.\label{hkcoe}\end{gather} (depending on our chosen $\pi$, cf. section \ref{resisec}) henceforth as an identification and regard the pair $(D,D_{{K}})$ as a filtered $(\varphi,N)$-module. Let us gather what we know from our preceding discussion about $(D,D_{{K}})$. There is a direct sum decomposition as $K^t$-vector space$$D=\bigoplus_{j\in{\mathbb Z}/f}\bigoplus_{i\in{\mathbb Z}/(d+1)}D_{ij}.$$Each $D_{ij}$ is stable for the monodromy operator $N$ on $D$. The $\sigma$-linear Frobenius operator $\varphi$ on $D$ is bijective and restricts to bijections $\varphi:D_{ij}\to D_{i(j-1)}$ if $j\ne 0\in{\mathbb Z}/f$ and $\varphi:D_{ij}\to D_{(i+1)(j-1)}$ if $j=0\in{\mathbb Z}/f$. For any $i,j$ there is a filtration $$D_{ij}=S_{ij}^0\supset S_{ij}^1\supset\ldots\supset S_{ij}^{d+1}=0$$such that $$S^s=\bigoplus_{i,j}S_{ij}^s$$is $\varphi$-stable for all $s$ and such that the $n(d+1)$-fold iterate of $\varphi$ acting on $S^s/S^{s+1}$ is multiplication with an element of $K^t$ of $p$-adic order $n(d+1)(d-s)+(n(d+1)-1)r_M$ (if $\pi^e=p$ for $e=n/f=[K:K^t]$ then this is the element $p^{n(d+1)(d-s)+(n(d+1)-1)r_M}$). Moreover, $S_{ij}^s/S_{ij}^{s+1}$ has $K^t$-dimension equal to $\mu(\Gamma, M)$ for all $i, j, s$ (this follows from formula (\ref{pecomp}) first for $j=0$ and then for all $j$ in view of the described cyclic structure of the bijective $\sigma$-linear endomorphism $\varphi$). As a $K$-vector space, $D_K$ decomposes into a direct sum$$D_K=\bigoplus_{i\in{\mathbb Z}/(d+1)}\bigoplus_{\delta\in{\rm Gal}(K/{\mathbb Q}_p)}(D_K)_{i,\delta}$$such that the Hodge filtration on $D_K$ is the sum of filtrations on the summands $(D_K)_{i,\delta}$ of the following shape:$$(D_K)_{i,\delta}=F^{r_M}_{{i,\delta}}\supset\ldots\supset F^{r_M+d+1}_{{i,\delta}}=0$$if $\delta\ne\id_K$. For $\delta=\id_K$ the filtration is of the form$$(D_K)_{i,\delta}=F^{r_M-\lambda_0}_{{i,\delta}}\supset\ldots\supset F^{r_M+d+1}_{{i,\delta}}=0$$such that all graded pieces are zero except possibly $F^{r_M-\lambda_j+j}_{{i,\delta}}/F^{r_M-\lambda_j+j+1}_{{i,\delta}}$ for $j=0,\ldots,d$. The isomorphism (\ref{hkcoe}) is the direct sum, over all $i$ and $j$, of isomorphisms\begin{gather}\bigoplus_{\delta\in{\rm Gal}(K/{\mathbb Q}_p)_{\sigma^j}}(D_K)_{i,\delta}\cong D_{ij}\otimes_{K^t}K.\label{hyokako}\end{gather}For all $i,i'\in{\mathbb Z}/(d+1)$ and all $j\in{\mathbb Z}/f$ there are canonical $K^t$-linear isomorphims$$\rho_{i,i'}^j:D_{ij}\cong D_{i'j}$$ such that $\rho_{i',i''}^j\circ\rho_{i,i'}^j=\rho_{i,i''}^j$ and $\rho_{i,i}^j=\id_K$, compatible with the monodromy operator $N$ and strict for the filtrations $S_{ij}^{\bullet}$. Moreover, after applying $\otimes_{K^t}K$ they respect the direct sum decompositions of the left hand sides of (\ref{hyokako}) and are strict for the Hodge filtrations on the $(D_K)_{i,\delta}$. Consider the $K$-rational $G$-representation $$M^*=\ho_K(M,K)\otimes{\rm det}^{\lambda_0}.$$It has highest weight $(\lambda_0^*\ge\lambda_1^*\ge\ldots\ge\lambda_d^*)$ with $\lambda_i^*=\lambda_0-\lambda_{d-i}$. Note that $\lambda_d^*=0$. Let $(D^*,D^*_{{K}})$ denote the filtered $(\varphi,N)$-module constructed as the $d$-th cohomology of ${\mathfrak E}_{M^*}$. We define filtrations $F_{\bullet,\bullet}^{\bullet}$ and $S_{\bullet\bullet}^{\bullet}$ as before (we are using the same symbols where we should use distinct ones, e.g. decorated by a superscript $*$). Then Poincar\'{e} duality provides us for all $i\in{\mathbb Z}/(d+1)$ and all $\delta\in{\rm Gal}(K/{\mathbb Q}_p)$ with a perfect pairing\begin{gather}(D_K)_{i,\delta}\times(D_K^*)_{i,\delta}\longrightarrow K\label{dualm}\end{gather}and Serre duality tells us that it induces for all $0\le j\le d+1$ perfect pairings$$F_{i,\delta}^{r_M+j}\times(D_K^*)_{i,\delta}/F_{i,\delta}^{r_{M^*}+d+1-j}\longrightarrow K$$if $\delta\ne{\bf 1}=\id_K$, and perfect pairings$$F_{i,{\bf 1}}^{r_M+j-\lambda_j}\times(D_K^*)_{i,{\bf 1}}/F_{i,{\bf 1}}^{r_{M^*}+d+1-j-\lambda_{d-j+1}}\longrightarrow K$$(see \cite{latt}). Similarly, if $\delta|_{K^t}=\sigma^j$ then the pairing (\ref{dualm}) induces perfect pairings$$(D_K)_{i,\delta}\cap S_{ij}^{s}\otimes_{K^t}K\quad\times\quad(D_K^*)_{i,\delta}/((D_K^*)_{i,\delta}\cap S_{ij}^{d+1-s}\otimes_{K^t}K)\longrightarrow K$$for all $0\le s\le d+1$. This follows from Borel-Serre duality (similarly to \cite{schn} p. 626, 633) or alternatively from Poincar\'{e} duality in log crystalline cohomology (using the characterization of $S^{\bullet}$ as a slope filtration for the Frobenius structure). 

By Theorem \ref{moweeas} the monodromy filtration on $D$ (associated as usual to the nilpotent monodromy operator $N$) coincides (up to renumbering) with the slope filtration $S^{\bullet}$.

For $0\le j\le d$ let$$\widetilde{F}^j=\bigoplus_{i\in{\mathbb Z}/(d+1)}(F^{r_M-\lambda_j+j}_{i,{\bf 1}}\oplus    \bigoplus_{\delta\in{\rm Gal}(K/{\mathbb Q}_p)\atop \delta\ne{\bf 1}}F_{i,\delta}^{r_M+j})$$(${\bf 1}=\id_K$). This defines a Hodge filtration$$D_K=\widetilde{F}^0\supset\ldots\supset \widetilde{F}^{d+1}=0.$$The combined conjectures of Schneider stated at the end of section \ref{conj} predict it splits the slope filtration $S^{\bullet}$ of $D$, i.e. that for all $0\le j\le d-1$ we have$$D_K=\widetilde{F}^{j+1}\bigoplus S^{d-j}\otimes_{K^t}K.$$(In particular this predicts the Hodge numbers of the Hodge filtration of $D_K$.) It is a straightforward computation to show that if this splitting conjecture holds true then $(D,D_{{K}})$ is a (weakly) admissible filtered $(\phi,N)$-module, hence defines a semistable $p$-adic ${\rm Gal}(\overline{K}/K)$-representation. Conversely, the (weak) admissibility of $(D,D_{{K}})$ may be regarded as a weak form of the splitting conjecture. 

Let us say that $\Gamma$ is {\it of arithmetic type} if ${\mathfrak X}_{\Gamma}$ carries a universal abelian scheme as in section \ref{appel} (see also section \ref{unifsect}) below. There we will show:

\begin{satz}\label{wead} If $\Gamma$ is of arithmetic type then the spectral sequence (\ref{tilhoa}) degenerates in $E_1$ and $(D,D_{{K}})$ is a weakly admissible filtered $(\varphi,N)$-module (over $(K^t,K)$).
\end{satz}

In particular, if $\Gamma$ is of arithmetic type then the statements in Proposition \ref{degeq} (ii) hold true.\\

{\it Remark:} If $d\ge2$, the author does not know how to read off from a given discrete cocompact and torsion free subgroup $\Gamma\subset{\rm SL}_{d+1}(K)$ wether it is of arithmetic type, or how many such $\Gamma$ are of arithmetic type.

\section{The universal $p$-divisible group}
\label{appel}

Let ${\mathcal S}=K^{d+1}$ be the standard representation of $G={\rm GL}\sb {d+1}(K)$. We identify the projective $d$-space ${\mathbb P}_K^d$ (with its $G$-action), in which we defined $X$ as an open subspace, with the projective space ${\mathbb P}({\mathcal S}^*)$ over the dual ${\mathcal S}^*$ of ${\mathcal S}$. Let $K'$ be an extension field of $K$, complete with respect to an absolute value which extends the one on $K$. Let $X'$ be the $K'$-rigid space obtained from $X$ by scalar extension. Each $E$-valued point $x\in X'(E)$ of $X'$, for a finite field extension $E/K'$, corresponds to a line through the origin in ${\mathcal S}^*\otimes_KE$, or equivalently to a hyperplane $H_x$ through the origin in ${\mathcal S}\otimes_KE$.\\ 

{\bf Definition} The {\it tautological subbundle} ${\mathcal F}'$ of ${{\mathcal S}}\otimes_K{\mathcal O}_{X'}$ is the rank-$d$-sub vector bundle with the following property: for all $x\in {X'}(E)$ (any $E$) the isomorphism$${\mathcal S}\otimes_{K}E\longrightarrow ({{\mathcal S}}\otimes_K{\mathcal O}_{{X'}_E})\otimes_{{\mathcal O}_{{X'}_E}}\kappa(x),\quad m\mapsto m\otimes 1$$restricts to an isomorphism $H_x\to ({\mathcal F}'{\otimes}_{K'}E)\otimes_{{\mathcal O}_{{X'}_E}}\kappa(x)$. (Here ${X'}_E={X'}\otimes_{K'}E$ and ${\mathcal O}_{{X'}_E}\to\kappa(x)$ is the natural surjection to the residue class field $\kappa(x)$ at $x$ (which is $\cong E$).)\\

We describe a convergent filtered $F$-isocrystal on ${\mathfrak X}$ which is the (covariant) Dieudonn\'{e} module of a certain universal $p$-divisible group. We follow \cite{rz}, especially 1.44, 1.45, 3.54 -- 3.75, 5.48, 5.49. Let $B$ be a central simple algebra of dimension $(d+1)^2$ over $K$ and with Brauer invariant $1/(d+1)$. Explicitly it may be presented as follows. Let $\tilde{K}$ be an unramified extension of degree $d+1$ of $K$ contained in $B$ and let $\eta\in{\rm Gal}(\tilde{K}/K)$ be the relative Frobenius automorphism. Then$$B=\tilde{K}[\Pi];\quad \Pi^{d+1}=\pi,\, \Pi x=\eta(x)\Pi\quad (x\in \tilde{K}).$$Let $\overline{\mathbb{F}}_p$ be an algebraic closure of the residue field of ${\mathbb Q}_p$, let $F_0=\quot(W(\overline{\mathbb{F}}_p))$. We view $F_0$ as a subfield of a fixed completed algebraic closure ${\mathbb C}_p$ of ${\mathbb{Q}_p}$. Fix an embedding $\epsilon:{K}\to{\mathbb{C}}_p$ (we will often suppress the symbol $\epsilon$ in our notation) and let $\breve{K}=KF_0$ inside ${\mathbb C}_p$ (using $\epsilon$); abstractly this is $\breve{K}=K\otimes_{K^{t,\epsilon}}F_0$. Denote by $\tau$ the $K$-automorphism of $\breve{K}$ defined by the relative Frobenius of $\overline{\mathbb{F}}_p$ over $k={\mathcal O}_K/(\pi)$. Set$$N=B\otimes_K\breve{K}.$$For $b\in B$ let $[b]\in{\rm End}(B)$ denote the {\it right} multiplication with $b$. Define the $\tau$-linear operator $\Phi'$ on $N$ as $[\pi\Pi^{-1}]\otimes\tau$. Fix an extension $\epsilon':\tilde{K}\to \breve{K}$ of $\epsilon$ (such that $\epsilon'\circ\eta=\tau\circ \epsilon'$). Then we have the decomposition (cf. \cite{rz} p.40)\begin{gather}N=\bigoplus_{i\in{\mathbb{Z}}/d+1}N_i\label{nstrich}\end{gather}$$N_i=\{n\in N\quad|\quad kn=\epsilon'(\eta^i(k))n\,\,{\mbox{for all}}\,\, k\in\tilde{K}\}$$where in the defining equation for $N_i$ the term $kn$ is formed with respect to the action of $\tilde{K}\subset B$ on $B$. The restriction $\Phi'_0$ of $\Phi'\circ[\pi\Pi^{-1}]^{-1}$ to ${N_0}$ is a $\tau$-linear operator $\Phi'_0:N_0\to N_0$. Letting $V_0=(N_0)^{\Phi'_0=1}$ we have $N_0=V_0\otimes_K{\breve{K}}$ and $\id_{V_0}\otimes\tau=\Phi'_0$. Taking intersections with $N_0\otimes_{\breve{K}}{\mathbb C}_p$ is a bijection between the set of $B$-stable ${\mathbb C}_p$-sub vector spaces of $N\otimes_{\breve{K}}{\mathbb C}_p$ of dimension $(d+1)d$, and the set of hyperplanes in the $(d+1)$-dimensional ${\mathbb C}_p$-vector space $N_0\otimes_{\breve{K}}{\mathbb C}_p$ (cf. \cite{rz} p.40/41). Similarly, the group of $B$-linear automorphisms of $N$ (w.r.t. right multiplication on the tensor factor $B$ of $N$) which commute with $\Phi'$ is identified with the group ${\rm GL}(V_0)$. We identify $V_0$ with the standard $K$-rational representation ${\mathcal S}$ of ${\rm GL}_{d+1}$ so that $G={\rm GL}_{d+1}(K)={\rm GL}(V_0)$. The decomposition (\ref{nstrich}) is $G$-stable and for each $i$, right multiplication with $\Pi$ induces a $G$-equivariant isomorphism $[\Pi]:N_i\cong N_{i-1}$. Let $${\bf M}_{{}}=B\otimes_{{\mathbb Q}_p}F_0.$$Let $\sigma$ be the absolute Frobenius automorphism on $F_0$ (extending $\sigma$ on $K^t$ if via $\epsilon$ we view $K^t$ as a subfield of $F_0$). Then (cf. \cite{rz} p.108) \begin{gather}{\bf M}_{{}}=\bigoplus_{j\in{\mathbb{Z}}/f} {\bf M}_{j}\label{fuldec}\end{gather}$${\bf M}_{j}=B\otimes_{K^{t,\sigma^j\epsilon}}F_0.$$We identify ${\bf M}_{0}=N$. We define the $\sigma$-linear Frobenius operator $\Phi_{\bf M}$ on the $F_0$-vector space ${\bf M}_{{}}$ as the sum of the maps $[p\Pi^{-1}]\otimes\sigma:{\bf M}_{0}\to {\bf M}_{1}$ and $p\otimes\sigma:{\bf M}_{j}\to {\bf M}_{j+1}$ (all $j\in({\mathbb{Z}}/f)-\{0\}$). Explicitly, let $${\bf M}=B\otimes_{{\mathbb Q}_p}F_0 \stackrel{\nu_j}{\longrightarrow}{\bf M}_{j}=B\otimes_{K^{t,\sigma^j\epsilon}}F_0$$denote the canonical projection map. If $d\in B$ then $[p\Pi^{-1}]\otimes\sigma:{\bf M}_{0}\to {\bf M}_{1}$ sends $\nu_0(d\otimes 1)$ to $([p\Pi^{-1}]\otimes 1)(\nu_1(d \otimes 1))$ and $p\otimes\sigma:{\bf M}_{j}\to {\bf M}_{j+1}$ sends $\nu_j(d \otimes 1)$ to $p \nu_{j+1}(d \otimes 1)$, cf. \cite{rz} p.39/40. The group of $B$-linear automorphisms of ${\bf M}$ which commute with $\Phi_{\bf M}$ is identified with the group $G={\rm GL}_{d+1}(K)={\rm GL}(V_0)$. This action of $G$ on ${\bf M}$ is $F_0$-linear and respects each ${\bf M}_{j}$; namely, the $G$-action on ${\bf M}_{j}$ is induced from the $G$-action on ${\bf M}_{0}=N$ via the isomorphism ${\bf M}_{j}\cong {\bf M}_{0}\otimes_{F_0^{\sigma^j}}F_0$.

We write $\breve{\mathfrak X}={\mathfrak X}\times_{\spf({\mathcal O}_K)}\spf({\mathcal O}_{\breve{K}})$ and we write $\breve{X}$ for the generic fibre of $\breve{\mathfrak X}$ (as a $\breve{K}$-rigid space). We write ${\bf M}$ and ${\bf M}_{j}$ also for the corresponding constant sheaves on $\breve{X}$. We define a subbundle ${\mathcal F}({\bf M}_{{}}\otimes_{F_0}{\mathcal O}_{\breve{X}})$ of the bundle ${\bf M}_{{}}\otimes_{F_0}{\mathcal O}_{\breve{X}}$ on ${\breve{X}}$ componentwise as follows. For each $j\ne0$ we let ${\mathcal F}({\bf M}_{j}\otimes_{F_0}{\mathcal O}_{\breve{X}})={\bf M}_{j}\otimes_{F_0}{\mathcal O}_{\breve{X}}$. To define ${\mathcal F}({\bf M}_{0}\otimes_{F_0}{\mathcal O}_{\breve{X}})$ consider the decomposition $$N_0\otimes_{F_0}{\mathcal O}_{\breve{X}}=\bigoplus_{\gamma\in{\rm Gal}(\breve{K}/F_0)}N_0\otimes_{\breve{K}^{\gamma}}{\mathcal O}_{\breve{X}}.$$Here $N_0\otimes_{\breve{K}^{\gamma}}{\mathcal O}_{\breve{X}}$ means the tensor product of $N_0$ and ${\mathcal O}_{\breve{X}}$ over $\breve{K}$ using the map $\breve{K}\stackrel{\gamma}{\to}\breve{K}\to{\mathcal O}_{\breve{X}}$; for $\gamma=\id$ we simply write $N_0\otimes_{\breve{K}}{\mathcal O}_{\breve{X}}$. We define ${\mathcal F}(N_0\otimes_{\breve{K}^{\gamma}}{\mathcal O}_{\breve{X}})=N_0\otimes_{\breve{K}^{\gamma}}{\mathcal O}_{\breve{X}}$ for $\gamma\ne\id\in{\rm Gal}(\breve{K}/F_0)={\rm Gal}({K}/K^t)$. We define ${\mathcal F}(N_0\otimes_{\breve{K}}{\mathcal O}_{\breve{X}})$ to be the tautological subbundle of $N_0\otimes_{\breve{K}}{\mathcal O}_{\breve{X}}={{\mathcal S}}\otimes_{{K}}{\mathcal O}_{\breve{X}}$ as defined above. Taking the sum over all $\gamma\in{\rm Gal}(\breve{K}/F_0)$ we obtain a subbundle ${\mathcal F}(N_0\otimes_{F_0}{\mathcal O}_{\breve{X}})$ of $N_0\otimes_{F_0}{\mathcal O}_{\breve{X}}$. For any $i$ we define the subbundle ${\mathcal F}(N_i\otimes_{F_0}{\mathcal O}_{\breve{X}})$ of $N_i\otimes_{F_0}{\mathcal O}_{\breve{X}}$ as the image of ${\mathcal F}(N_0\otimes_{F_0}{\mathcal O}_{\breve{X}})$ under the isomorphism $N_0\otimes_{F_0}{\mathcal O}_{\breve{X}}\cong N_{-i}\otimes_{F_0}{\mathcal O}_{\breve{X}}$ induced by $[\Pi^i]$. These ${\mathcal F}(N_i\otimes_{F_0}{\mathcal O}_{\breve{X}})$ for all $i$ make up a subbundle ${\mathcal F}({\bf M}_{0}\otimes_{F_0}{\mathcal O}_{\breve{X}})$ of ${\bf M}_{0}\otimes_{F_0}{\mathcal O}_{\breve{X}}$, and the ${\mathcal F}({\bf M}_{j}\otimes_{F_0}{\mathcal O}_{\breve{X}})$ for all $j$ make up our ${\mathcal F}({\bf M}\otimes_{F_0}{\mathcal O}_{\breve{X}})$.

In order to make completely explicit the coincidence of ${\mathcal F}({\bf M}\otimes_{F_0}{\mathcal O}_{\breve{X}})$ with the filtration ${\mathcal F}^1$ of ${\bf M}\otimes_{F_0}{\mathcal O}_{\breve{X}}$ as defined in \cite{rz} p.40, consider the decomposition$${\bf M}\otimes_{F_0}{\mathcal O}_{\breve{X}}=B\otimes_{{\mathbb Q}_p}{\mathcal O}_{\breve{X}}=\bigoplus_{\iota\in{\rm Gal}(K/{\mathbb Q}_p)}B\otimes_{K^{\iota}}{\mathcal O}_{\breve{X}}.$$The two term direct sum decomposition of ${\bf M}\otimes_{F_0}{\mathcal O}_{\breve{X}}$ into the summand for $\iota=\id_K$ on the one hand and the sum of the summands for all $\iota\ne\id_K$ on the other hand is used in \cite{rz} p.40 to define a filtration ${\mathcal F}^1$ of ${\bf M}\otimes_{F_0}{\mathcal O}_{\breve{X}}$ as follows: on the ($\iota=\id_K$)-summand, ${\mathcal F}^1$ is the tautological subbundle; on the ($\iota\ne\id_K$)-summands, ${\mathcal F}^1$ is everything. Now this two term decomposition of ${\bf M}\otimes_{F_0}{\mathcal O}_{\breve{X}}$ is the same as the following: the sum of the $N_{-i}\otimes_{\breve{K}}{\mathcal O}_{\breve{X}}=[\Pi^i]N_0\otimes_{\breve{K}}{\mathcal O}_{\breve{X}}$ for all $i\in{\mathbb Z}/(d+1)$ on the one hand, and the sum of the ${\bf M}_{j}\otimes_{F_0}{\mathcal O}_{\breve{X}}$ for all $0\ne j\in{\mathbb Z}/f={\rm Gal}(K^t/{\mathbb Q}_p)$ and of the $N_{-i}\otimes_{\breve{K}^{\gamma}}{\mathcal O}_{\breve{X}}=[\Pi^i]N_0\otimes_{\breve{K}^{\gamma}}{\mathcal O}_{\breve{X}}$ for all $i\in{\mathbb Z}/(d+1)$ and all $\id\ne\gamma\in{\rm Gal}(\breve{K}/F_0)={\rm Gal}({K}/K^t)$ on the other hand. Thus our ${\mathcal F}({\bf M}\otimes_{F_0}{\mathcal O}_{\breve{X}})$ is the ${\mathcal F}^1$ as defined in \cite{rz} p.40. 

Let $\breve{\mathfrak X}_0=\breve{\mathfrak X}\times_{\spf({\mathcal O}_{\breve{K}})}\spec(\overline{\mathbb{F}}_p)$ denote the special fibre of $\breve{\mathfrak X}$. We define a convergent $F$-isocrystal ${\bf M}_{{}}\otimes_{W(\overline{\mathbb{F}}_p)}{\mathcal O}^{conv}_{\breve{\mathfrak {X}}_0}$ on $\breve{\mathfrak X}_0$ as follows. Let ${\mathcal O}^{conv}_{\breve{\mathfrak {X}}_0}$ denote the structure sheaf of the convergent site of $\breve{\mathfrak X}_0$: it assigns to each enlargement $Z$ of $\breve{\mathfrak X}_0$ the structure sheaf ${\mathcal O}_Z$ of $Z$. Now the value of ${\bf M}_{{}}\otimes_{W(\overline{\mathbb{F}}_p)}{\mathcal O}^{conv}_{\breve{\mathfrak {X}}_0}$ at an enlargement $Z$ of $\breve{\mathfrak X}_0$ is given by ${\bf M}_{{}}\otimes_{W(\overline{\mathbb{F}}_p)}{\mathcal O}_Z$. The Frobenius structure on ${\bf M}_{{}}\otimes_{W(\overline{\mathbb{F}}_p)}{\mathcal O}^{conv}_{\breve{\mathfrak {X}}_0}$, which is felt only at enlargements $Z$ of $\breve{\mathfrak X}_0$ with a lifting $\Phi_Z:{\mathcal O}_Z\to{\mathcal O}_Z$ of Frobenius, is given by $\Phi_{\bf M}\otimes \Phi_Z$; in particular, it is 'constant'.

We regard the filtered ${\mathcal O}_{\breve{X}}$-module ${\bf M}\otimes_{F_0}{\mathcal O}_{\breve{X}}$ as the evaluation of ${\bf M}_{{}}\otimes_{W(\overline{\mathbb{F}}_p)}{\mathcal O}^{conv}_{\breve{\mathfrak {X}}_0}$ at the enlargement $\breve{\mathfrak X}$ of $\breve{\mathfrak X}_0$. In this way ${\bf M}\otimes_{W(\overline{\mathbb{F}}_p)}{\mathcal O}^{conv}_{\breve{\mathfrak {X}}_0}$ becomes a convergent filtered $F$-isocrystal on $\breve{\mathfrak X}$. \\

Let ${\mathcal O}_B$ be the maximal order ${\mathcal O}_{\tilde{K}}[\Pi]$ in $B$. A formal $p$-divisible group ${\mathcal G}$ over an ${\mathcal O}_{\breve{K}}$-scheme $S$ (in particular $p$ is locally nilpotent on $S$) with an ${\mathcal O}_B$-action ${\mathcal O}_B\to{\rm End}({\mathcal G})$ is called {\it special} if it satisfies the condition \cite{rz} p. 110. We fix a special formal ${\mathcal O}_B$-module ${\mathcal G}_0$ over $\spec(\overline{\mathbb{F}}_p)$ of $K$-height $(d+1)^2$ (for the concept of $K$-height see \cite{rz} Lemma 3.53 and p.109); it is unique up to quasi-isogeny (\cite{rz} Lemma 3.60). According to Drinfel'd (\cite{rz} Theorem 3.72), the formal ${\mathcal O}_{\breve{K}}$-scheme $\breve{\mathfrak X}$ represents the functor which associates to an ${\mathcal O}_{\breve{K}}$-scheme $S$ the set of $({\mathcal G},\rho)$ consisting of a special formal ${\mathcal O}_B$-module ${\mathcal G}$ over $S$ and a quasiisogeny $\rho:{\mathcal G}_0\times_{\spec(\overline{\mathbb{F}}_p)}\overline{S}\longrightarrow {\mathcal G}\times_S\overline{S}$ of formal ${\mathcal O}_B$-modules over $\overline{S}$ of height $0$. Here we write $\overline{S}=S\otimes_{\mathbb{Z}_p}{\mathbb{F}}_p$.

Hence a universal special formal ${\mathcal O}_B$-module ${\mathcal G}^u$ over $\breve{\mathfrak X}$. (Covariant) Dieudonn\'{e} theory assigns to the special fibre ${\mathcal G}^u_0$ of ${\mathcal G}^u$ a convergent $F$-isocrystal $\breve{\mathfrak E}$ on $\breve{\mathfrak X}_0$. The lifting ${\mathcal G}^u$ of ${\mathcal G}^u_0$ then defines, in addition, a filtration on the evaluation $\breve{\mathfrak E}({\breve{\mathfrak X}})$ of $\breve{\mathfrak E}$ at the enlargement ${\breve{\mathfrak X}}$ of ${\breve{\mathfrak X}}_0$, as follows. Let us denote by $Lie({\mathcal G}^u)$ the Lie-algebra of ${\mathcal G}^u$ and by $H_{1,dR}({\mathcal G}^u)$ the Lie-algebra of the universal vectorial extension of ${\mathcal G}^u$. These give rise to coherent ${\mathcal O}_{\breve{X}}$-modules $Lie({\mathcal G}^u)\otimes_{{\mathcal O}_{\breve{K}}}\breve{K}$ and $H_{1,dR}({\mathcal G}^u)\otimes_{{\mathcal O}_{\breve{K}}}\breve{K}$ and in fact, $H_{1,dR}({\mathcal G}^u)\otimes_{{\mathcal O}_{\breve{K}}}\breve{K}=\breve{\mathfrak E}({\breve{\mathfrak X}})$. There is a canonical surjection\begin{gather}\gamma:\breve{\mathfrak E}({\breve{\mathfrak X}})\longrightarrow Lie({\mathcal G}^u)\otimes_{{\mathcal O}_{\breve{K}}}\breve{K}.\end{gather}The filtration ${\mathcal F}^{\bullet}\breve{\mathfrak E}({\breve{\mathfrak X}})$ of $\breve{\mathfrak E}({\breve{\mathfrak X}})$ is given by ${\mathcal F}^{t}\breve{\mathfrak E}({\breve{\mathfrak X}})=0$ if $t\ge 2$, by ${\mathcal F}^{1}\breve{\mathfrak E}({\breve{\mathfrak X}})=\ker(\gamma)$ and by ${\mathcal F}^{t}\breve{\mathfrak E}({\breve{\mathfrak X}})=\breve{\mathfrak E}({\breve{\mathfrak X}})$ if $t\le 0$. Altogether, $\breve{\mathfrak E}$ becomes a filtered convergent $F$-isocrystal on $\breve{\mathfrak X}$. Let us denote by$$\nabla:\breve{\mathfrak E}({\breve{\mathfrak X}})\longrightarrow\breve{\mathfrak E}({\breve{\mathfrak X}})\otimes_{{\mathcal O}_{\breve{X}}}\Omega^1_{\breve{X}}$$the connection provided by the isocrystal structure (see e.g. \cite{and} 3.6.7 where this connection is discussed for the equivalent contravariant Dieudonn\'{e} theory).

\begin{pro}\label{pdivco} The group of $B$-equivariant quasiisogenies ${\mathcal G}_0\to {\mathcal G}_0$ (of height zero) is naturally identified with $G$. This provides an action of $G$ on $\breve{\mathfrak X}$ and on $\breve{\mathfrak E}$. We have a canonical $G$-equivariant isomorphism of convergent filtered $F$-isocrystals with $G$-action\begin{gather}\breve{\mathfrak E}\cong{\bf M}_{{}}\otimes_{W(\overline{\mathbb{F}}_p)}{\mathcal O}^{conv}_{\breve{\mathfrak {X}}_0}\end{gather}Explicitly, the Frobenius $\Phi$ on $\breve{\mathfrak E}$ corresponds to $\Phi_{\bf M}\otimes \Phi_{{\mathcal O}_{\breve{\mathfrak X}}^{conv}}$ on ${\bf M}_{{}}\otimes_{W(\overline{\mathbb{F}}_p)}{\mathcal O}^{conv}_{\breve{\mathfrak {X}}_0}$, and under the induced isomorphism\begin{gather}\breve{\mathfrak E}({\breve{\mathfrak X}})\cong {\bf M}_{{}}\otimes_{F_0}{\mathcal O}_{\breve{X}}\label{fileval}\end{gather}of evaluations on $\breve{\mathfrak X}$, the filtration step ${\mathcal F}^{1}\breve{\mathfrak E}({\breve{\mathfrak X}})$ corresponds to ${\mathcal F}({\bf M}_{{}}\otimes_{F_0}{\mathcal O}_{\breve{X}})$ and the connection $\nabla$ corresponds to $\id_{{\bf M}_{{}}}\otimes_{F_0}d$.
\end{pro}

{\sc Proof:} This is explained in \cite{rz} 3.58 -- 3.75, 5.48, 5.49. The convergent $F$-isocrystal $\breve{\mathfrak E}$ is constant, with value the $F$-isocrystal associated to ${\mathcal G}_0$ which is ${\bf M}_{{}}$ (cf. also \cite{and} 6.2.5). That the filtration varies as stated is the coincidence of the period map as defined in \cite{rz} with the period map as defined by Drinfel'd (in \cite{rz} 5.48, 5.49, the proof of this coincidence is attributed to Faltings).\hfill$\Box$\\

Let again $M$ be an irreducible $K$-rational representation of $G={\rm GL}\sb {d+1}(K)$ of highest weight $(\lambda_0\ge\lambda_1\ge\ldots\ge\lambda_d=0)$. We define a $G$-equivariant filtered  convergent $F$-isocrystal $\breve{\mathfrak E}_{M}$ on $\breve{\mathfrak X}$ by the same recipe as in section \ref{filconfiso}, as follows. As a $G$-equivariant convergent isocrystal on the reduction $\breve{\mathfrak X}_0$ it is$$\breve{\mathfrak E}_{M}={\bf M}_{M}\otimes_{{\mathcal O}_{K^t}}{\mathcal O}^{conv}_{\breve{\mathfrak {X}}_0}$$with diagonal $G$-action. If $Z$ is an enlargement of $\breve{\mathfrak X}_0$ endowed with a lifting $\Phi_Z:{\mathcal O}_Z\to{\mathcal O}_Z$ of the $p$-power Frobenius endomorphism of the reduction then the Frobenius on $\breve{\mathfrak E}_{M}(Z)={\bf M}_{M}\otimes_{{\mathcal O}_{K^t}}{\mathcal O}_Z$ is $\Phi_{{\bf M}_{M}}\otimes\Phi_Z$. The evaluation $\breve{\mathfrak E}_{M}(\breve{\mathfrak X})$ of $\breve{\mathfrak E}_{M}$ at $\breve{\mathfrak X}$ decomposes as\begin{align}\breve{\mathfrak E}_{M}(\breve{\mathfrak X})&=\bigoplus_{j\in{\mathbb Z}/f}\bigoplus_{i\in{\mathbb Z}/(d+1)}M_{\sigma^{j}}\otimes_{{\mathcal O}_{K^t}}{\mathcal O}_{\breve{\mathfrak X}}\notag \\{}&=\bigoplus_{j\in{\mathbb Z}/f}\bigoplus_{i\in{\mathbb Z}/(d+1)}\bigoplus_{\delta\in{\rm Gal}(K/{\mathbb Q}_p)_{\sigma^{j}}}M_{\delta}\otimes_{{\mathcal O}_{K}}{\mathcal O}_{\breve{\mathfrak X}}\notag\end{align}where we write ${\rm Gal}(K/{\mathbb Q}_p)_{\sigma^{j}}=\{\gamma \in {\rm Gal}(K/{\mathbb Q}_p);\,\,\gamma|_{K^t}=\sigma^{j}\}$. We filter the direct summands $M_{\delta}\otimes_{{\mathcal O}_{K}}{\mathcal O}_{\breve{\mathfrak X}}$ as follows. If $\delta=\id_K$ we let\begin{gather}{f}^s(M\otimes_{{\mathcal O}_{K}}{\mathcal O}_{\breve{\mathfrak X}})={\mathcal O}_{\breve{\mathfrak X}}.\overline{u}(z)({f}^sM)\label{glaube}\end{gather}for $s\in{\mathbb Z}$. If $\delta\ne\id_K$ we let$${f}^{r_M}(M_{\delta}\otimes_{{\mathcal O}_{K}}{\mathcal O}_{\breve{\mathfrak X}})=M_{\delta}\otimes_{{\mathcal O}_{K}}{\mathcal O}_{\breve{\mathfrak X}}$$and ${f}^{r_M+1}(M_{\delta}\otimes_{{\mathcal O}_{K}}{\mathcal O}_{\breve{\mathfrak X}})=0$. (Recall that we defined $r_M=\sum_{i=0}^d\lambda_i=\sum_{i=1}^d\lambda_i$.) Together we get a filtration ${f}^{\bullet}\breve{\mathfrak E}_{M}(\breve{\mathfrak X})$ of $\breve{\mathfrak E}_{M}(\breve{\mathfrak X})$.

\begin{lem}\label{filglei} Take $M={\mathcal S}$.\\(a)  The filtration step ${f}^1({{\mathcal S}}\otimes_{{\mathcal O}_K}{\mathcal O}_{\breve{\mathfrak X}})$ of ${{\mathcal S}}\otimes_K{\mathcal O}_{\breve{X}}={{\mathcal S}}\otimes_{{\mathcal O}_K}{\mathcal O}_{\breve{\mathfrak X}}$ as defined in formula (\ref{glaube}) is the tautological subbundle.\\
(b) $\breve{\mathfrak E}$ and $\breve{\mathfrak E}_{\mathcal S}$ are isomorphic as $G$-equivariant filtered convergent $F$-isocrystals.
\end{lem}

{\sc Proof:} (a) We have ${f}^{2}{\mathcal S}=0$, ${f}^1{\mathcal S}=<e_1,\ldots,e_d>$ and ${f}^0{\mathcal S}={\mathcal S}$ in the standard basis $e_0,\ldots,e_d$ of ${\mathcal S}$. From this an easy computation gives the result.\\(b) By Proposition \ref{pdivco} we know $\breve{\mathfrak E}\cong {\bf M}_{{}}\otimes_{W(\overline{\mathbb{F}}_p)}{\mathcal O}^{conv}_{\breve{\mathfrak {X}}_0}$. Next we construct an isomorphism \begin{gather}{\bf M}\cong{\bf M}_{\mathcal S}\otimes_{K^t}F_0\label{cois}\end{gather}of $G$-equivariant $F_0$-vector spaces with $\sigma$-linear endomorphism. In the decomposition ${\bf M}_{{}}=\oplus_{j\in{\mathbb{Z}}/f} {\bf M}_{j}$ (see (\ref{fuldec})) we may regard the $F_0$-vector space ${\bf M}_{-j}=B\otimes_{K^{t,\sigma^{-j}\epsilon}}F_0$ as $({\bf M}_{0})_{\sigma^{j}}$ (send $d\otimes y\in {\bf M}_{-j}$ to $d\otimes\sigma^j(y)\in({\bf M}_{0})_{\sigma^{j}}$). Recall the decomposition ${\bf M}_{0}=N=\oplus_{i\in{\mathbb Z}/(d+1)}N_i$ and that right multiplication with $\Pi$ induces isomorphisms $[\Pi]:N_i\to N_{i-1}$. We obtain decompositions ${\bf M}_{-j}=({\bf M}_{0})_{\sigma^{j}}=\oplus_{i\in{\mathbb Z}/(d+1)}(N_i)_{\sigma^{j}}$ and isomorphisms $[\Pi]:(N_i)_{\sigma^{j}}\to(N_{i-1})_{\sigma^{j}}$ for any $j$. The identification $N_0\cong{\mathcal S}\otimes_K{\breve K}={\mathcal S}\otimes_{K^t}F_0$ as $G$-representations gives rise, for any $0\le i\le d$ (we insist on this particular representative modulo $(d+1)$ in order to get a welldefined power of $\Pi$ in the following isomorphism), to the isomorphism of $G$-representations $$(N_i)_{\sigma^{j}}\stackrel{[\Pi^i]}{\to}(N_0)_{\sigma^{j}}\cong({\mathcal S}\otimes_K{\breve K})_{\sigma^{j}}\cong({\mathcal S}\otimes_{K^t}F_0)_{\sigma^{j}}\cong{\mathcal S}_{\sigma^{j}}\otimes_{K^t}F_0$$ where the last isomorphism sends $x\otimes y\in ({\mathcal S}\otimes_{K^t}F_0)_{\sigma^{j}}$ to $x\otimes \sigma^{-1}(y)\in {\mathcal S}_{\sigma^{j}}\otimes_{K^t}F_0$. We define the searched for isomorphism (\ref{cois}) as the direct sum of these isomorphisms, letting the direct summand $(N_i)_{\sigma^{j}}$ of ${\bf M}$ correspond to the $(j,i)$ direct summand of ${\bf M}_{\mathcal S}\otimes_{K^t}F_0$. The Frobenius and the filtration on $\breve{\mathfrak E}_{\mathcal S}$ are defined in terms of the number $r_{\mathcal S}$. As the highest weight of ${\mathcal S}$ is given by $(\lambda_0,\ldots,\lambda_d)=(1,0,\ldots,0)$, this number is $r_{\mathcal S}=1$. Now tracing back our definitions of the respective Frobenii we see that indeed $\Phi_{\bf M}$ on ${\bf M}$ corresponds to $\Phi_{{\bf M}_{\mathcal S}}\otimes{\sigma}$ on ${\bf M}_{\mathcal S}\otimes_{K^t}F_0$ under (\ref{cois}). Hence (\ref{cois}) defines an isomorphism of $G$-equivariant convergent $F$-isocrystals $\breve{\mathfrak E}\cong\breve{\mathfrak E}_{\mathcal S}$. Moreover, it respects the respective filtrations on evaluations at $\breve{\mathfrak X}$. Indeed, (a) describes the filtration on the $\delta={\rm id}_K$-direct summand of $\breve{\mathfrak E}_{\mathcal S}(\breve{\mathfrak X})$, and this coincides with our definition of ${\mathcal F}(N_0\otimes_{{\mathcal O}_{\breve{K}}}{\mathcal O}_{\breve{\mathfrak X}})$ and hence with the filtration on the $\delta={\rm id}_K$-direct summand of $\breve{\mathfrak E}(\breve{\mathfrak X})$. On the $\delta\ne{\rm id}_K$-direct summands, in either case the filtrations jump from zero to everything (at the same filtration step).\hfill$\Box$\\

The general theory tells us that $\breve{\mathfrak E}$ is isomorphic with the {\it contra}variant Dieudonne module of the $p$-divisible group $\check{\mathcal G}^u$ over ${\breve{\mathfrak X}}$ {\it dual} to ${\mathcal G}^u$, and it is this interpretation of $\breve{\mathfrak E}$ which will be of interest to us. Let
$\breve{\mathfrak E}^*$ denote the {\it contra}variant Dieudonne module of ${\mathcal G}^u$. Let $${\mathfrak D}=\breve{\mathfrak E}\oplus\breve{\mathfrak E}^*.$$ Let again $M$ be arbitrary and recall that we defined $r_M=\sum_{h=0}^d\lambda_h=\sum_{h=0}^{d-1}\lambda_h$. Let the symmetric group $S_{r_M}$ act on the $r_M$-fold direct sum ${\mathfrak D}^{r_M}$ through permutations of the summands. By functoriality of taking exterior powers this induces an action of $S_{r_M}$ on the $r_M$-th exterior power $\bigwedge^{r_M}({\mathfrak D}^{r_M})$ of ${\mathfrak D}^{r_M}$, hence an action of the group algebra ${\mathbb Q}[S_{r_M}]$ on the filtered $F$-isocrystal $\bigwedge^{r_M}({\mathfrak D}^{r_M})$.

\begin{pro}\label{iotas} There exist an embedding of $G$-equivariant convergent filtered $F$-isocrystals$$\breve{\mathfrak E}_{M}\stackrel{\iota}{\longrightarrow}\bigwedge^{r_M}({\mathfrak D}^{r_M})$$whose image is respected by the action of $S_{r_M}$ on $\bigwedge^{r_M}({\mathfrak D}^{r_M})$ and split with respect to all structure elements, i.e. it is a direct summand of $\bigwedge^{r_M}({\mathfrak D}^{r_M})$ as a convergent filtered $F$-isocrystal.\end{pro} 

{\sc Proof:} We will construct $\iota$ as a composition of $G$-equivariant embeddings$$\breve{\mathfrak E}_{M}\stackrel{\iota_1}{\longrightarrow}(\breve{\mathfrak E}_{{\mathcal S}})^{\otimes r_M}\stackrel{\iota_2}{\longrightarrow}{\mathfrak D}^{\otimes r_M}\stackrel{\iota_3}{\longrightarrow}\bigwedge^{r_M}({\mathfrak D}^{r_M})$$all of which are respected by $S_{r_M}$ and have the stated splitting property. For $1\le i\le r_M$ let $\lambda_i:{\mathfrak D}\to {\mathfrak D}^{r_M}$ denote the inclusion which identifies ${\mathfrak D}$ with the $i$-th direct summand of ${\mathfrak D}^{r_M}$. We define $\iota_3$ as the map which sends $x_1\otimes\ldots\otimes x_{r_M}$ to $\lambda_1(x_1)\wedge\ldots\wedge\lambda_{r_M}(x_{r_M})$. By Lemma \ref{filglei}(b) we have ${\mathfrak D}\cong\breve{\mathfrak E}_{{\mathcal S}}\oplus\breve{\mathfrak E}^*$ and this gives an obvious definition of $\iota_2$. That $\iota_2$ and $\iota_3$ are split as stated is clear by elementary linear algebra of direct sums and exterior products. To define $\iota_1$ we describe a tensor power variant $\breve{\mathfrak E}_{{\mathcal S},r_M}$ of $\breve{\mathfrak E}_{{\mathcal S}}$. Consider the $K^t$-vector space with $G$-action$${\bf M}_{{\mathcal S},r_M}=\bigoplus_{j\in{\mathbb Z}/f}\bigoplus_{i\in{\mathbb Z}/(d+1)}({\mathcal S}^{\otimes r_M})_{\sigma^{j}}.$$We endow ${\bf M}_{{\mathcal S},r_M}$ with a $\sigma$-linear operator $\Phi_{{\bf M}_{{\mathcal S},r_M}}$ as follows: it sends an element $x$ from the $(j,i)$-direct summand to the element $p^{r_M}x$ in the $(j-1,i)$ direct summand if $j\ne0$, to the element $p^{r_M}x$ in the $(j-1,i-1)$ direct summand if $j=0$ and $i\ne 0$, and to the element $p^{r_M}\pi^{-r_M}x$ in the $(j-1,i-1)$ direct summand if $j=0$ and $i=0$. We then let$$\breve{\mathfrak E}_{{{\mathcal S},r_M}}={\bf M}_{{{\mathcal S},r_M}}\otimes_{{\mathcal O}_{K^t}}{\mathcal O}^{conv}_{\breve{\mathfrak {X}}_0}$$with diagonal $G$-action and Frobenius $\Phi_{{\bf M}_{{\mathcal S},r_M}}\otimes\Phi$. The evaluation $\breve{\mathfrak E}_{{\mathcal S},r_M}(\breve{\mathfrak X})$ decomposes as\begin{align}\breve{\mathfrak E}_{{\mathcal S},r_M}(\breve{\mathfrak X})&=\bigoplus_{j\in{\mathbb Z}/f}\bigoplus_{i\in{\mathbb Z}/(d+1)}({\mathcal S}^{\otimes r_M})_{\sigma^{j}}\otimes_{{\mathcal O}_{K^t}}{\mathcal O}_{\breve{\mathfrak X}}\notag \\{}&=\bigoplus_{j\in{\mathbb Z}/f}\bigoplus_{i\in{\mathbb Z}/(d+1)}\bigoplus_{\delta\in{\rm Gal}(K/{\mathbb Q}_p)_{\sigma^{j}}}({\mathcal S}^{\otimes r_M})_{\delta}\otimes_{{\mathcal O}_{K}}{\mathcal O}_{\breve{\mathfrak X}}.\notag\end{align} We filter the direct summands $({\mathcal S}^{\otimes r_M})_{\delta}\otimes_{{\mathcal O}_{K}}{\mathcal O}_{\breve{\mathfrak X}}$ as follows. If $\delta=\id_K$ we filter$$({\mathcal S}^{\otimes r_M})_{\delta}\otimes_{{\mathcal O}_{K}}{\mathcal O}_{\breve{\mathfrak X}}=({\mathcal S}\otimes_{{\mathcal O}_{K}}{\mathcal O}_{\breve{\mathfrak X}})\otimes_{{\mathcal O}_{\breve{\mathfrak X}}}\ldots\otimes_{{\mathcal O}_{\breve{\mathfrak X}}}({\mathcal S}\otimes_{{\mathcal O}_{K}}{\mathcal O}_{\breve{\mathfrak X}})$$by the tensor product filtration ${f}^{\bullet}(({\mathcal S}^{\otimes r_M})_{\delta}\otimes_{{\mathcal O}_{K}}{\mathcal O}_{\breve{\mathfrak X}})$ of the filtrations on the tensor factors ${\mathcal S}\otimes_{{\mathcal O}_{K}}{\mathcal O}_{\breve{\mathfrak X}}$ given by formula (\ref{glaube}) (for $M={\mathcal S})$. On the other hand, if $\delta\ne\id_K$ and $s\le r_M$ we let$${f}^{s}(({\mathcal S}^{\otimes r_M})_{\delta}\otimes_{{\mathcal O}_{K}}{\mathcal O}_{\breve{\mathfrak X}})=({\mathcal S}^{\otimes r_M})_{\delta}\otimes_{{\mathcal O}_{K}}{\mathcal O}_{\breve{\mathfrak X}}$$and if $s\ge r_M+1$ we let ${f}^{s}(({\mathcal S}^{\otimes r_M})_{\delta}\otimes_{{\mathcal O}_{K}}{\mathcal O}_{\breve{\mathfrak X}})=0$. Together we get a filtration ${f}^{\bullet}\breve{\mathfrak E}_{{\mathcal S},r_M}(\breve{\mathfrak X})$ of $\breve{\mathfrak E}_{{\mathcal S},r_M}(\breve{\mathfrak X})$. Now we will construct the searched for embedding $\iota_1$ as a composition of embeddings$$\breve{\mathfrak E}_{M}\stackrel{\iota'_1}{\longrightarrow}\breve{\mathfrak E}_{{\mathcal S},r_M}\stackrel{\iota''_1}{\longrightarrow}(\breve{\mathfrak E}_{{\mathcal S}})^{\otimes r_M}.$$Consider the decomposition as an isocrystal$$(\breve{\mathfrak E}_{{\mathcal S}})^{\otimes r_M}=\bigoplus_{j_{\bullet}\in({\mathbb Z}/f)^{r_M}}\bigoplus_{i_{\bullet}\in({\mathbb Z}/(d+1))^{r_M}}({\mathcal S}_{\sigma^{j_{1}}}\otimes_{K^t}\ldots\otimes_{K^t}{\mathcal S}_{\sigma^{j_{r_M}}})\otimes_{{\mathcal O}_{K^t}}{\mathcal O}^{conv}_{\breve{\mathfrak {X}}_0}.$$Taking the sum over all direct summands with $j_1=\ldots=j_{r_M}$ and $i_1=\ldots=i_{r_M}$ defines a direct summand isomorphic with\begin{gather}\bigoplus_{j\in{\mathbb Z}/f}\bigoplus_{i\in{\mathbb Z}/(d+1)}(({\mathcal S}^{\otimes_{K^t}r_M})_{\sigma^{j}})\otimes_{{\mathcal O}_{K^t}}{\mathcal O}^{conv}_{\breve{\mathfrak {X}}_0}.\label{zwie}\end{gather}Here ${\mathcal S}^{\otimes_{K^t}r_M}$ denotes the $r_M$-fold tensor power of ${\mathcal S}$ over $K^t$. It decomposes as$${\mathcal S}^{\otimes_{K^t}r_M}=\bigoplus_{\delta_{\bullet}\in{\rm Gal}(K/K^t)^{r_M-1}}{\mathcal S}_{\delta_1}\otimes_K{\mathcal S}_{\delta_2}\otimes_K\ldots\otimes_K{\mathcal S}_{\delta_{r_M-1}}\otimes_K{\mathcal S}.$$In this way the $r_M$-fold tensor power ${\mathcal S}^{\otimes r_M}$ of ${\mathcal S}$ over $K$ becomes a direct summand of ${\mathcal S}^{\otimes_{K^t}r_M}$, namely the one corresponding to $\delta_1=\ldots=\delta_{r_M-1}={\bf 1}$. Hence$$\breve{\mathfrak E}_{{\mathcal S},r_M}=\bigoplus_{j\in{\mathbb Z}/f}\bigoplus_{i\in{\mathbb Z}/(d+1)}({\mathcal S}^{\otimes r_M})_{\sigma^{j}}\otimes_{{\mathcal O}_{K^t}}{\mathcal O}^{conv}_{\breve{\mathfrak {X}}_0}$$is a direct summand of (\ref{zwie}). Our definitions are just made in such a way that this direct summand is stable under Frobenius and strict with respect to filtrations. We have thus defined ${\iota''_1}$. The induced action of $S_{r_M}$ on ${\mathfrak D}^{\otimes r_M}$ and on $(\breve{\mathfrak E}_{{\mathcal S}})^{\otimes r_M}$ is given by permutations of the tensor factors, hence the same is true for the induced action of $S_{r_M}$ on each direct summand of $\breve{\mathfrak E}_{{\mathcal S},r_M}$. By \cite{calu} p.211 or \cite{bchlls} there is a $G$-equivariant embedding $\iota:M\to{\mathcal S}^{\otimes r_M}$ whose image admits a characterization of the following type. There exists a subset $T\subset S_{r_M}$ consisting of transpositions, and there exists a collection $\{P_n\}_{n}$ of subsets of $S_{r_M}$ such that for all $x\in{\mathcal S}^{\otimes r_M}$ we have $x\in\bi(\iota)$ if and only if $\sigma(x)=\epsilon(\sigma)x$ for each $\sigma\in T$ and $(\sum_{\sigma\in P_n}\epsilon(\sigma)\sigma)(x)=0$ for each $P_n$ (here $\epsilon(\sigma)$ denotes the signum of $\sigma$, and $\sigma\in S_{r_M}$ acts on ${\mathcal S}^{\otimes r_M}$ through permutations of the tensor factors). It follows that there exists an idempotent $U_{M}\in {\mathbb Q}[S_{r_M}]$ which acts on ${\mathcal S}^{\otimes r_M}$ as a projector onto $\iota(M)$. Namely, take $U_M$ as the product of all $${\rm id}-\frac{1}{2}(\sigma-\epsilon(\sigma){\rm id})\quad\quad(\sigma\in T),\quad\quad\quad\quad {\rm id}-|P_n|^{-1}\sum_{\sigma\in P_n}\epsilon(\sigma)\sigma\quad\quad (P_n\in\{P_n\}_n).$$ From $\iota$ we derive similar embeddings $$\iota^{(j)}:M_{\sigma^{j}}\to({\mathcal S}^{\otimes r_M})_{\sigma^{j}}$$for any $j\in{\mathbb Z}/f$. Source and target of the searched for embedding $\iota'_1$ both are graded over the same index set $\{(j,i)\in ({\mathbb Z}/f)\times({\mathbb Z}/(d+1))\}$. We define $\iota'_1$ as the graded embedding which is $\iota^{(j)}$ on the $(i,j)$-component. It is clear that the image of $\iota'_1$ is the image of $U_M$ acting on $\breve{\mathfrak E}_{{\mathcal S},r_M}$. In particular, the complement $\bi(1-U_M)$ is stable under Frobenius since the actions of elements of $S_{r_M}$ on $\bigwedge^{r_M}({\mathfrak D}^{r_M})$ respect the Frobenius structure (because this is true for the permutation action on ${\mathfrak D}^{r_M}$). Finally, since $\iota:M\to{\mathcal S}^{\otimes r_M}$ is $G$-equivariant it follows that the filtration on ${\mathcal S}^{\otimes r_M}$, which is the tensor product of the filtrations ${f}^{\bullet}{\mathcal S}$ on its tensor factors, induces the filtration ${f}^{\bullet}M$ on $M$. It then follows that the induced embedding $M\otimes_{{\mathcal O}_{K}}{\mathcal O}_{\breve{\mathfrak X}}\to{\mathcal S}^{\otimes r_M}\otimes_{{\mathcal O}_{K}}{\mathcal O}_{\breve{\mathfrak X}}=({\mathcal S}\otimes_{{\mathcal O}_{K}}{\mathcal O}_{\breve{\mathfrak X}})^{\otimes r_M}$ is strict for filtrations, hence $\iota'_1(\breve{\mathfrak X})$ is strict for filtrations at the direct summands of $\breve{\mathfrak E}_{M}(\breve{\mathfrak X})$ with $j=0\in{\mathbb Z}/f$. For the other direct summands this is clear anyway.\hfill$\Box$\\

{\it Remark:} We can now give a conceptual explanation for our definition of the filtration on $\breve{\mathfrak E}_{M}(\breve{\mathfrak X})$. Firstly, for $M={\mathcal S}$ it follows from Proposition \ref{pdivco} and Lemma \ref{filglei} that the filtration on $\breve{\mathfrak E}_{{\mathcal S}}(\breve{\mathfrak X})$ is exactly the one imposed by geometry, i.e. under the isomorphism (\ref{fileval}) it becomes the natural filtration on the covariant Dieudonne module $\breve{\mathfrak E}$ of ${\mathcal G}^u$. Secondly, by Proposition \ref{iotas} a general $M$ admits a split embedding $M\to{\mathcal S}^{\otimes r_M}$ and the tensor product filtration on $\breve{\mathfrak E}_{{\mathcal S}^{\otimes r_M}}(\breve{\mathfrak X})$ restricts to the one on $\breve{\mathfrak E}_{M}(\breve{\mathfrak X})$ as we defined it.\\

Let $\Gamma<{\rm SL}_{d+1}(K)$ be a cocompact discrete torsionfree subgroup as before such that ${{\mathfrak X}}_{\Gamma}$ has strictly semistable reduction. Then ${\breve{\mathfrak X}}_{\Gamma}=\Gamma\backslash{\breve{\mathfrak X}}$ is just ${{\mathfrak X}}_{\Gamma}\times_{\spf({\mathcal O}_K)}\spf({\mathcal O}_{\breve{K}})$ and this is a projective ${\mathcal O}_{\breve{K}}$-scheme with strictly semistable reduction. Its generic fibre is ${\breve{X}}_{\Gamma}=\Gamma\backslash \breve{X}=X_{\Gamma}\otimes_K{\breve{K}}$, a smooth projective ${\breve{K}}$-variety.

By its modular characterization the universal $p$-divisible group ${\mathcal G}^u$ over ${\breve{\mathfrak X}}$ descends to a $p$-divisible group ${\mathcal G}^u_{\Gamma}$ over ${\breve{\mathfrak X}}_{\Gamma}$. Similarly $\check{\mathcal G}^u$ descends to a $p$-divisible group $\check{\mathcal G}^u_{\Gamma}$ over ${\breve{\mathfrak X}}_{\Gamma}$ (the dual of ${\mathcal G}_{\Gamma}^u$). The contravariant Dieudonne module ${\mathfrak D}^{\Gamma}$ of ${\mathcal G}^u_{\Gamma}\times\check{\mathcal G}^u_{\Gamma}\to{\breve{\mathfrak X}}$, as a convergent filtered $F$-isocrystal on ${\breve{\mathfrak X}}_{\Gamma}$, is a descent of ${\mathfrak D}$.

From now on we assume that $\Gamma$ is of arithmetic type in the sense that \cite{rz} 6.51 applies to ${\breve{\mathfrak X}}_{\Gamma}$ (cf. section \ref{unifsect}). Then there is an abelian scheme (with ${\mathcal O}_B$-action) $f:{\mathcal A}\to{\breve{\mathfrak X}}_{\Gamma}$ over ${\breve{\mathfrak X}}_{\Gamma}$ whose associated $p$-divisible group is ${\mathcal G}^u_{\Gamma}\times\check{\mathcal G}^u_{\Gamma}$. For $r\ge0$ let$$f^{r}:{\mathcal A}^{r}={\mathcal A}\times_{{\breve{\mathfrak X}}_{\Gamma}}\ldots\times_{{\breve{\mathfrak X}}_{\Gamma}}{\mathcal A}\longrightarrow {\breve{\mathfrak X}}_{\Gamma}$$ denote the $r$-fold fibre power of ${\mathcal A}$ over ${\breve{\mathfrak X}}_{\Gamma}$, again an abelian scheme over ${\breve{\mathfrak X}}_{\Gamma}$.

The relative crystalline cohomology ${\bf R}f^r_{0,crys}({\mathcal A}_0^r/{\breve{\mathfrak X}}_{\Gamma,0})$ of the reduction $f_0^r:{\mathcal A}_0^r\to{\breve{\mathfrak X}}_{\Gamma,0}$ of $f^r$ is a convergent $F$-isocrystal on ${\breve{\mathfrak X}}_{\Gamma,0}$. Its evaluation at the lifting ${\breve{\mathfrak X}}_{\Gamma}$ can be identified with the relative de Rham cohomology ${\bf R}f^r_{\mathbb Q,*}(\Omega^{\bullet}_{{\mathcal A}^r_{\mathbb Q}/{\breve{X}}_{\Gamma}})$ of the generic fibre $f^r_{\mathbb Q}:{\mathcal A}^r_{\mathbb Q}\to{\breve{X}}_{\Gamma}$ of $f^r$ and as such it is provided with a Hodge filtration. Together this defines a filtered convergent $F$-isocrystal (with ring structure) on ${\breve{\mathfrak X}}_{\Gamma}$ which we denote by ${\mathbb R}f_*({\mathcal A}^r/{\breve{\mathfrak{X}}}_{\Gamma})$.

\begin{lem}\label{drelco} There is a canonical isomorphism of filtered convergent $F$-isocrystals (with ring structure) $${\mathbb R}f^r_*({\mathcal A}^r/{\breve{\mathfrak{X}}}_{\Gamma})\cong \bigwedge(({\mathfrak D}^{\Gamma})^{r})$$ on ${\breve{\mathfrak X}}_{\Gamma}$. In particular, ${\mathbb R}^rf_*({\mathcal A}^r/{\breve{\mathfrak{X}}}_{\Gamma})\cong \bigwedge^r(({\mathfrak D}^{\Gamma})^{r})$
\end{lem}

{\sc Proof:} For $r=1$ this follows from our definitions. For $r>1$ it then follows from general facts on the (relative) cohomology of abelian schemes.\hfill$\Box$\\

Let again $M$ be as before. The $G$-equivariant filtered convergent $F$-isocrystals $\breve{\mathfrak E}_{M}$ descend to filtered convergent $F$-isocrystals $\breve{\mathfrak E}_{M}^{\Gamma}$ on ${\breve{\mathfrak X}}_{\Gamma}$, and similarly do their tensor and exterior powers. Denoting the filtration of $\breve{\mathfrak E}_{M}^{\Gamma}({\breve{\mathfrak X}}_{\Gamma})$ by ${f}^{\bullet}\breve{\mathfrak E}_{M}^{\Gamma}({\breve{\mathfrak X}}_{\Gamma})$ we filter the de Rham complex $\breve{\mathfrak E}_{M}^{\Gamma}({\breve{\mathfrak X}}_{\Gamma})\otimes\Omega^{\bullet}_{\breve{X}_{\Gamma}}$ by setting\begin{gather}\breve{\mathcal F}^j_M=[({f}^{j}\breve{\mathfrak E}_{M}^{\Gamma}({\breve{\mathfrak X}}_{\Gamma}))\otimes\Omega^{0}_{\breve{X}_{\Gamma}}\to({f}^{j-1}\breve{\mathfrak E}_{M}^{\Gamma}({\breve{\mathfrak X}}_{\Gamma}))\otimes\Omega^{1}_{\breve{X}_{\Gamma}}\to({f}^{j-2}\breve{\mathfrak E}_{M}^{\Gamma}({\breve{\mathfrak X}}_{\Gamma}))\otimes\Omega^{2}_{\breve{X}_{\Gamma}}\to\ldots].\label{drcofi}\end{gather}This filtration gives rise to the spectral sequence\begin{gather}E_1^{st}=H^{s+t}({\breve X}_{\Gamma},\breve{\mathcal F}_M^{s}/\breve{\mathcal F}_M^{s+1})\Rightarrow H^{s+t}({\breve X}_{\Gamma},\breve{\mathfrak E}_{M}^{\Gamma}({\breve{\mathfrak X}}_{\Gamma})\otimes\Omega^{\bullet}_{\breve{X}_{\Gamma}}).\label{tilho}\end{gather}

\begin{satz}\label{degen} The spectral sequence (\ref{tilho}) degenerates in $E_1$.
\end{satz}

{\sc Proof:} We have ${\bf R}^{r_M}(f_{\mathbb Q}^{r_M})_*(\Omega^{\bullet}_{{\mathcal A}_{\mathbb Q}^{r_M}/{\breve{X}}_{\Gamma}})=\bigwedge^{r_M}(({\mathfrak D}^{\Gamma})^{r_M})({\breve{\mathfrak X}})$ as filtered ${\mathcal O}_{\breve{X}_{\Gamma}}$-modules with connection. Let $$\breve{\mathcal K}^{\bullet}=(\bigwedge^{r_M}(({\mathfrak D}^{\Gamma})^{r_M})({\breve{\mathfrak X}}))\otimes_{{\mathcal O}_{\breve{X}_{\Gamma}}}{\Omega}^{\bullet}_{\breve{X}_{\Gamma}}$$denote its de Rham complex, let $(\breve{\mathcal F}^{j}{\breve{\mathcal K}^{\bullet}})_j$ denote its associated filtration. By Proposition \ref{kkue} the spectral sequence $$E_1^{st}=H^{s+t}(\breve{X}_{\Gamma},\breve{\mathcal F}^s{\breve{\mathcal K}^{\bullet}}/\breve{\mathcal F}^{s+1}{\breve{\mathcal K}^{\bullet}})\Longrightarrow H^{s+t}(\breve{X}_{\Gamma},\breve{\mathcal K}^{\bullet})$$degenerates in $E_1$. The degeneration of (\ref{tilho}) now follows once we recognize $\breve{\mathfrak E}_{M}^{\Gamma}({\breve{\mathfrak X}}_{\Gamma})$ as a direct summand (as a filtered ${\mathcal O}_{\breve{X}}$-module with connection) of $\bigwedge^{r_M}(({\mathfrak D}^{\Gamma})^{r_M})({\breve{\mathfrak X}_{\Gamma}})$. But for that we may use Proposition \ref{iotas}.\hfill$\Box$\\

{\it Remark:} The inclusions $(\iota'_1)^{\Gamma}(\breve{\mathfrak X}):\breve{\mathfrak E}_{M}^{\Gamma}\to\breve{\mathfrak E}_{{\mathcal S},r_M}^{\Gamma}$ and ${\iota_3}^{\Gamma}:({\mathfrak D}^{\otimes r_M})^{\Gamma}{\to}\bigwedge^{r_M}({\mathfrak D}^{r_M})^{\Gamma}$ derived from Proposition \ref{iotas} (and its proof) can be cut out geometrically by ${\breve{\mathfrak X}}_{\Gamma}$-endomorphisms of ${\mathcal A}^{r_M}$. Indeed, for $(\iota'_1)^{\Gamma}$ this follows from the fact that the action of ${\mathbb Q}[S_{r_M}]$ on $\bigwedge^{r_M}({\mathfrak D}^{r_M})^{\Gamma}$ is induced by permutations of the factors of ${\mathcal A}^{r_M}$. To see it for $\iota_3^{\Gamma}$ let $\lambda_i$ (for $1\le i\le r_M$) be as in the proof of Proposition \ref{iotas} and let $p_i:{\mathfrak D}^{r_M}\to{\mathfrak D}^{r_M}$ denote the projector onto $\oplus_{j\ne i}\lambda_j({\mathfrak D})$ with kernel $\lambda_i({\mathfrak D})$. Then the image of $\iota_3$ is the image of the idempotent $\prod_{1\le i\le r_M}(\id-\wedge^{r_M}(p_i))$ acting on $\bigwedge^{r_M}({\mathfrak D}^{r_M})$. Clearly the $(p_i)^{\Gamma}$ are induced by endomorphisms of ${\mathcal A}^{r_M}$.\\

Let $$\breve{D}=H_{conv}^{d}(\breve{\mathfrak X}_{{\Gamma},0},\breve{\mathfrak E}_{M}^{\Gamma})$$denote the $d$-th logarithmic convergent cohomology group of the convergent $F$-isocrystal $\breve{\mathfrak E}_{M}^{\Gamma}$ on the reduction $\breve{\mathfrak X}_{{\Gamma},0}$ of $\breve{\mathfrak X}_{{\Gamma}}$ relative to $\breve{T}=(\spf(W(\overline{\mathbb F})),1\mapsto 0)$. Thus $D$ is a finite dimensional $F_0$-vector space endowed with a $\sigma$-linear Frobenius endomorphism $\varphi$ and a nilpotent $F_0$-linear endomorphism $N$ such that $N\varphi=p\varphi N$. Alternatively, $D$ is the $d$-th (rational) log crystalline cohomology group of an integral structure (as an $F$-crystal on $\breve{\mathfrak X}_{{\Gamma},0}$) inside $\breve{\mathfrak E}_{M}^{\Gamma}$ with respect to the logarithmic divided power thickening $\breve{T}$ of $(\spec(\overline{\mathbb F}),1\mapsto 0)$. Let$$\breve{D}_{\breve{K}}=H_{dR}^{d}(\breve{X}_{{\Gamma}},\breve{\mathfrak E}_{M}^{\Gamma}(\breve{\mathfrak X}_{{\Gamma}}))$$denote the $d$-th de Rham cohomology group of the filtered ${\mathcal O}_{\breve{X}_{{\Gamma}}}$-module with connection $\breve{\mathfrak E}_{M}^{\Gamma}(\breve{\mathfrak X}_{{\Gamma}})$. We can view the pair $(\breve{D},\breve{D}_{\breve{K}})$ via the isomorphism $\breve{D}_{\breve{K}}\cong \breve{D}\otimes_{F_0}\breve{K}$ as a filtered $(\varphi,N)$-module (over $(F_0,\breve{K})$). 

\begin{satz}\label{weakad} $(\breve{D},\breve{D}_{\breve{K}})$ is a weakly admissible filtered $(\varphi,N)$-module (over $(F_0,\breve{K})$). The monodromy filtration coincides with the weight filtration.  
\end{satz}

{\sc Proof:} ${\mathcal A}^{r_M}$ is a projective ${\mathcal O}_{\breve{K}}$-scheme with semistable reduction. Let $H_{conv}^{d+r_M}({\mathcal A}^{r_M}_0)$ denote the $(d+r_M)$-th logarithmic convergent cohomology group with constant coefficients, with respect to the basis $\breve{T}$, of its special fibre. Equivalently one can take (rational) logarithmic crystalline cohomology. We have the isomorphism\begin{gather}H_{dR}^{d+r_M}({\mathcal A}^{r_M}_{\mathbb Q}/\breve{K})\cong H_{conv}^{d+r_M}({\mathcal A}^{r_M}_0)\otimes_{F_0}\breve{K}\label{hkcon}\end{gather}with the $(d+r_M)$-th de Rham cohomology group $H_{dR}^{d+r_M}({\mathcal A}^{r_M}_{\mathbb Q}/\breve{K})$, and by Tsuji's theorem \cite{tsu}, the pair $(H_{conv}^{d+r_M}({\mathcal A}^{r_M}_0),H_{dR}^{d+r_M}({\mathcal A}^{r_M}_{\mathbb Q}/\breve{K}))$ (glued together via the isomorphism (\ref{hkcon})) is a weakly admissible filtered $(\varphi,N)$-module. We intend to show that $(\breve{D},\breve{D}_{\breve{K}})$ is a direct summand (with respect to all structure elements). That $\breve{D}_{\breve{K}}$ is a direct summand of $H_{dR}^{d+r_M}({\mathcal A}^{r_M}_{\mathbb Q}/\breve{K})$ (in the filtered sense) follows from Proposition \ref{kkue} by the proof of Theorem \ref{degen}. Consider the spectral sequence corresponding to the composition of derived functors $${\bf R}\Gamma_{conv}({\mathcal A}^{r_M}_0,.)={\bf R}\Gamma_{conv}(\breve{\mathfrak X}_{\Gamma,0},.)\circ {\bf R}f^{r_M}_{0,{conv}}(.)$$where ${\bf R}f^{r_M}_{0,{conv}}(.)$ denotes relative convergent cohomology (as $f^{r_M}_0$ is (classically) smooth it may be computed as relative crystalline cohomology, cf.  \cite{faalet}, \cite{ogusII}, \cite{ogusc}, \cite{shihoArX}). It tells us that $H_{conv}^{d+r_M}({\mathcal A}^{r_M}_0)={\bf R}^{d+r_M}\Gamma_{conv}({\mathcal A}^{r_M}_0,{\mathcal O}^{conv}_{{\mathcal A}^{r_M}_0})$ (with ${\mathcal O}^{conv}_{{\mathcal A}^{r_M}_0}$ denoting the trivial convergent $F$-isocrystal on ${\mathcal A}^{r_M}_0$) has an exhausting and separated filtration, indexed by $q\in{\mathbb Z}$, such that the $q$-th graded piece is a subquotient of $H_{conv}^{d+r_M-q}(\breve{\mathfrak X}_{\Gamma,0},{\bf R}^qf^{r_M}_{0,{conv}}({\mathcal O}^{conv}_{{\mathcal A}^{r_M}_0}))$. But in fact this subquotient is all of $H_{conv}^{d+r_M-q}(\breve{\mathfrak X}_{\Gamma,0},{\bf R}^qf^{r_M}_{0,{conv}}({\mathcal O}^{conv}_{{\mathcal A}^{r_M}_0}))$: this follows from the corresponding fact for de Rham cohomology of the generic fibres (which we already proved, using Proposition \ref{kkue}) by a dimension count argument in view of the isomorphisms (\ref{hkcon}) and $$H^s(\breve{X}_{\Gamma},{\bf R}^{t}(f_{\mathbb Q}^{r_M})_*(\Omega^{\bullet}_{{\mathcal A}_{\mathbb Q}^{r_M}/{\breve{X}}_{\Gamma}}))\cong  H_{conv}^{s}(\breve{\mathfrak X}_{\Gamma,0},{\bf R}^tf^{r_M}_{0,{conv}}({\mathcal O}^{conv}_{{\mathcal A}^{r_M}_0}))\otimes_{F_0}\breve{K}$$for all $s,t$. Moreover, this $q$-th graded piece is in fact a direct summand of $H_{conv}^{d+r_M}({\mathcal A}^{r_M}_0)$ (with respect to all structure elements): namely, it is characterized as the subspace on which the $n$-multiplication maps $n:{\mathcal A}^{r_M}_0\to {\mathcal A}^{r_M}_0$ on the abelian ${\breve{\mathfrak X}_{\Gamma,0}}$-scheme ${\mathcal A}^{r_M}_0$ act for any $n\in{\mathbb N}$ as multiplication with $n^q$ (so that this graded piece is the image of an idempotent in the group ring ${\mathbb Q}[{\rm End}_{\breve{\mathfrak X}_{\Gamma,0}}({\mathcal A}^{r_M}_0)]$; the action of this ring respects all structure elements). Now taking $q=r_M$ we obtain that $H_{conv}^{d}(\breve{\mathfrak X}_{\Gamma,0},{\bf R}^{r_M}f^{r_M}_{0,{conv}}({\mathcal O}^{conv}_{{\mathcal A}^{r_M}_0}))=H_{conv}^{d}(\breve{\mathfrak X}_{\Gamma,0},\bigwedge^{r_M}(({\mathfrak D}^{\Gamma})^{r_M}))$ is a direct summand of $H_{conv}^{d+r_M}({\mathcal A}^{r_M}_0)$. But $\breve{D}$ is a direct summand of $H_{conv}^{d}(\breve{\mathfrak X}_{\Gamma,0},\bigwedge^{r_M}(({\mathfrak D}^{\Gamma})^{r_M}))$ (with respect to all structure elements) as one can see by the same arguments as in the proof of Theorem \ref{degen}. Hence $\breve{D}$ is a direct summand of $H_{conv}^{d+r_M}({\mathcal A}^{r_M}_0)$ (with respect to all structure elements). It is clear that our direct sum embeddings on the crystalline and on the de Rham side agree (via the isomorphisms $\breve{D}_{\breve{K}}\cong \breve{D}\otimes_{F_0}\breve{K}$ and (\ref{hkcon})), by the functoriality of all constructions.

Finally, the coincidence of the monodromy filtration with the weight filtration follows from Theorem \ref{moweeas} by base field extension.\hfill$\Box$\\

{\sc Proof of Theorem \ref{wead}:} The spectral sequence (\ref{tilho}) arises from the spectral sequence (\ref{tilhoa}) by the base field extension $K\to\breve{K}$ (or $K^t\to F_0$), hence Theorem \ref{degen} implies the degeneration of (\ref{tilhoa}). Similarly, $(\breve{D},\breve{D}_{\breve{K}})$ arises from $(D,D_{{K}})$ by the base field extension $K^t\to F_0$ (with $\sigma$-linear extension of the Frobenius from $D$ to $\breve{D}$) as follows by well known cohomological principals. Hence theorem \ref{weakad} implies the weak admissibility of $(D,D_{{K}})$.\hfill$\Box$\\

\begin{satz} Let $r\ge0$ and $i\ge0$. For $H_{crys}^{i}({\mathcal A}^{r}_0)=H_{conv}^{i}({\mathcal A}^{r}_0/T)$ (cohomology with constant coefficients), the monodromy-weight conjecture (as formulated in \cite{mokr}) holds true: the monodromy filtration coincides with the weight filtration.
\end{satz}

{\sc Proof:} In the proof of Theorem \ref{weakad} we learned that the log crystalline (or log convergent) cohomology groups of ${\mathcal A}^{r}_0$ (with constant coefficients) are just direct sums of the log convergent cohomology groups of $\breve{\mathfrak X}_{{\Gamma},0}$ with coefficients in ${\mathbb R}f^r_*({\mathcal A}^r/{\breve{\mathfrak{X}}}_{\Gamma})$. Consequently, in view of Lemma \ref{drelco} it is enough to show that $\bigwedge(({\mathfrak D}^{\Gamma})^{r})$ is a direct sum of convergent $F$-isocrystals on $\breve{\mathfrak X}_{\Gamma}$ such that the monodromy filtration on their convergent cohomology groups is pure. By a base field extension argument it is therefore enough to show that $\bigwedge(({\mathfrak D}^{\Gamma})^{r})$ decomposes, possibly after raising Frobenius endomorphisms to some power, into direct summands each of which is isomorphic to some Frobenius twisting of a convergent $F$-isocrystal as treated in Theorem \ref{moweeas} (via base field extension $K\to\breve{K}$). 

In Lemma \ref{filglei} we saw $\breve{\mathfrak E}\cong \breve{\mathfrak E}_{\mathcal S}$. A similar description can be given for $\breve{\mathfrak E}^*$, as follows. Dual to the definition of ${\mathbb M}$ resp. of ${\mathbb M}_{\mathcal S}$ (recall Lemma \ref{filglei}, in particular the isomorphism (\ref{cois})) we now define the $G$-representation$${\bf M}^*=\bigoplus_{j\in{\mathbb Z}/f}\bigoplus_{i\in{\mathbb Z}/(d+1)}{\mathcal S}^*_{\sigma^{-j}}\otimes_{K^t}F_0$$where ${\mathcal S}^*$ denotes the dual of the standard representation ${\mathcal S}$ of $G$. We endow ${\bf M}^*$ with a $\sigma$-linear operator $\Phi_{{\bf M}^*}$ as follows: it sends an element $m\otimes 1$ from the $(j,i)$ direct summand to the element $m\otimes 1$ in the $(j-1,i)$ direct summand if $j\ne0\in{\mathbb Z}/f$, to the element $m\otimes 1$ in the $(j-1,i+1)$ direct summand if $j=0\in{\mathbb Z}/f$ and $i\ne d\in{\mathbb Z}/(d+1)$, and to the element $\pi m\otimes 1$ in the $(j-1,i+1)$ direct summand if $j=0\in{\mathbb Z}/f$ and $i=d\in{\mathbb Z}/(d+1)$. Then $\Phi_{{\bf M}^{*}}$ on ${\bf M}^{*}$ is isoclinic of slope $\frac{1}{n(d+1)}$. It turns out that $\breve{\mathfrak E}^*\cong {\bf M}^*\otimes_{W(\overline{\mathbb F}_p)}{\mathcal O}_{\breve{\mathfrak X}}^{conv}$ as convergent $F$-isocrystals, where the Frobenius structure on ${\bf M}^*\otimes_{W(\overline{\mathbb F}_p)}{\mathcal O}_{\breve{\mathfrak X}}^{conv}$ is given by $\Phi_{{\bf M}^{*}}\otimes \Phi_{{\mathcal O}_{\breve{\mathfrak X}}^{conv}}$. From this it follows that ${\mathfrak D}^{\Gamma}=\breve{\mathfrak E}^{\Gamma}\oplus (\breve{\mathfrak E}^*)^{\Gamma}$ is, after passing to a power of the Frobenius structure, the direct sum of twists of convergent $F$-isocrystals as treated in Theorem \ref{moweeas} (via base field extension $K\to\breve{K}$). But then it is not hard to see that $\bigwedge(({\mathfrak D}^{\Gamma})^{r})$ admits a direct sum decomposition with direct summands with the same property.\hfill$\Box$\\

\section{Appendix: Uniformization of Shimura varieties}

\label{unifsect}

For the following compare with \cite{har}, \cite{rz}. Let $E$ be a totally real field with a single place $v$ dividing $p$ such that there exists an isomorphism of the $v$-adic completion $E_{v}$ with $K$ (we fix this isomorphism). We choose an imaginary quadratic field $L_0$ in which $p$ splits as a product $p=p_1p_2$ and let $L=L_0E$. We let $v_i$ (for $i=1,2$) be the prime of $L$ dividing $p_i$ and $v$. We choose a central division algebra ${\mathcal B}$ of dimension $(d+1)^2$ over $L$ which stays a division algebra at $v_1$ and $v_2$ such that (for $B$ as before) there are isomorphisms$${\mathcal B}_{v_1}\cong B,\quad\quad\quad {\mathcal B}_{v_2}\cong B^{opp}$$(extending $E_{v}\cong K$). Let ${\mathcal B}'$ denote the central simple algebra of dimension $(d+1)^2$ over $L$ which splits at $v_1$ and $v_2$ but whose invariants at other finite places are the same as those of ${\mathcal B}$. 
We assume ${\mathcal B}$ endowed with an involution $\alpha$ of the second kind relative to the non-trivial $L_0$-automorphism of $L$; then also ${\mathcal B}'$ is endowed with such an involution $\alpha$. We define$$G{\mathcal G}=GU({\mathcal B},\alpha)=\{g\in R_{L/{\mathbb Q}}{\mathcal B}^{opp,\times}\,\,|\,\,g\cdot g^{\alpha}\in{\mathbb G}_{m,\mathbb Q}\},$$$$G{\mathcal J}=GU({\mathcal B}',\alpha)=\{g\in R_{L/{\mathbb Q}}({\mathcal B}')^{opp,\times}\,\,|\,\,g\cdot g^{\alpha}\in{\mathbb G}_{m,\mathbb Q}\}$$which we view as algebraic groups over ${\mathbb Q}$. Our group $G$ considered earlier comes in via $G={\rm GL}_{d+1}(K)\to{\rm PGL}_{d+1}(K)\cong G{\mathcal J}_{ad}({\mathbb Q}_p)$. Finally, letting$${\mathcal G}=U({\mathcal B},\alpha)=\{g\in{\mathcal B}^{opp,\times}\,\,|\,\,g\cdot g^{\alpha}=1\}$$we assume that for one real place $s$ of $E$ we have ${\mathcal G}(E_s)=U(d,1)$ while for all the other real places $s'\ne s$ of $E$ we have ${\mathcal G}(E_{s'})=U(d+1)$ (compact unitary group).

We fix an embedding of algebraic closures $\phi:\overline{\mathbb Q}\to \overline{\mathbb Q}_p$ and let $\Phi$ be the set of embeddings $\eta:L\to\overline{\mathbb Q}$ such that $\phi\circ\eta$ induces the prime $p_1$ of $L_0$. Then $\Phi$ is a CM-type of $L$ and its reflex field is $L_0$. We use $\Phi$ to identify ${\mathcal B}\otimes_{\mathbb Q}{\mathbb C}$ with the product $M(d+1,{\mathbb C})^{d+1}$ of $d+1$ copies of the matrix algebra and let the first copy of $M(d+1,{\mathbb C})$ correspond to the complex place $t\in\Phi$ above $s$. We define$$h_{\Phi}:{\mathbb C}^{\times}\cong R_{{\mathbb C}/{\mathbb R}}{\mathbb G}_{m,{\mathbb C}}({\mathcal B}^{opp}\otimes_{{\mathbb Q}}{\mathbb C})^{\times},$$$$h_{\phi}(z)=(\diag(z\cdot I_{d},\bar{z}), zI_{d+1},\ldots, zI_{d+1}).$$Let $G{\mathcal G}_{\infty}=\prod_{s'}G{\mathcal G}(E_{s'})$, the product running over all real places of $E$. The $G{\mathcal G}_{\infty}$-conjugacy class $X_d$ of homomorphisms from $R_{{\mathbb C}/{\mathbb R}}{\mathbb G}_{m,{\mathbb C}}$ to $G{\mathcal G}({\mathbb R})$ containing $h_{\Phi}$ may be $G{\mathcal G}_{\infty}$-equivariantly identified with the unit ball in ${\mathbb C}^d$. Let $Sh=Sh(G{\mathcal G},X_d)$ denote the associated Shimura variety. It is defined over $L$, where $L$ is viewed as a subfield of ${\overline{\mathbb Q}}$ via $t$. Let $C_p\subset G{\mathcal G}({\mathbb Q}_p)$ be the unique maximal compact subgroup and let $C^p\subset G{\mathcal G}({\mathbb A}_f^p)$ be a sufficiently small open compact subgroup. Extending scalars $L\to L_{v_1}\cong E_v\cong K$ we may view the Shimura variety $Sh$ as being defined over $K$. At level $C=C_pC^p$ it then has a model $Sh_C$ over ${\mathcal O}_{K}$. If ${\mathcal H}(G{\mathcal G}({\mathbb A}_f)//C)$ denotes the Hecke-algebra with respect to $C$, there is a ${\mathcal H}(G{\mathcal G}({\mathbb A}_f)//C)$ equivariant isomorphism of $p$-adic formal schemes \cite{rz} 6.51\begin{gather}G{\mathcal J}({\mathbb Q})\backslash \breve{{\mathfrak X}}\times G{\mathcal G}({\mathbb A}_f)/C\cong \widehat{Sh_C}\otimes_{{\mathcal O}_K}{\mathcal O}_{\breve{K}}.\label{rzunifo}\end{gather}Here $G{\mathcal J}({\mathbb Q})$ acts diagonally on $\breve{{\mathfrak X}}\times G{\mathcal G}({\mathbb A}_f)$, on the first factor $\breve{{\mathfrak X}}$ through the isomorphism ${\rm PGL}_{d+1}(K)\cong G{\mathcal J}_{ad}({\mathbb Q}_p)$ and the natural action of ${\rm PGL}_{d+1}(K)$ on $\breve{{\mathfrak X}}$, on the second factor through the isomorphism $G{\mathcal G}({\mathbb A}_f)\cong G{\mathcal J}({\mathbb A}_f)$. The left hand side of (\ref{rzunifo}) is a disjoint sum of formal schemes $\breve{{\mathfrak X}}_{\Gamma}$ as considered before.

\section{Appendix: Residue maps and log convergent cohomology}
\label{resisec}

Let ${\mathcal O}_K$ be a complete discrete valuation ring with fixed uniformizer $\pi\in{\mathcal O}_K$, with fraction field $K$ of characteristic zero and perfect residue field $k$ of characteristic $p>0$. Then $K$ is a finite field extension of $K^t$, the fraction field of the ring $W(k)$ of Witt vectors with coefficients in $k$. Let $\sigma$ denote the Frobenius endomorphism of $K^t$. Let ${\mathfrak X}$ be a strictly semistable proper formal ${\mathcal O}_K$-scheme of pure relative dimension $d$. We endow ${\mathfrak X}$ with the log structure defined by its special fibre $Y={\mathfrak X}\otimes_{{\mathcal O}_K}k$ which is a normal crossings divisor on ${\mathfrak X}$. We endow $Y$ with the pull back log structure. Let $\{Y_s\}_{s\in R}$ denote the set of irreducible components of $Y$, all of which are smooth proper $k$-schemes of pure dimension $d$. We fix a total ordering of $R$. For $i\in{\mathbb N}$ let $S_i$ be the set of subsets of $R$ consisting of precisely $i$ elements. For $\sigma\in S_i$ let $Y_{\sigma}=\cap_{s\in\sigma}Y_s$ and then let $Y^i=\coprod_{\sigma\in S_i}Y_{\sigma}$.

The following results from \cite{hk} are presented there mainly in terms of log rigid cohomology, but here we are going to present them in terms of log convergent cohomology instead --- as $Y$ and hence all its components and component intersections are proper this leads to the same cohomology groups, by \cite{hk} Proposition 5.6. 

Let $T$ denote the formal log scheme with underlying formal scheme ${\rm Spf}(W(k))$ and with log structure defined by the chart ${\mathbb N}\to W(k),$ $1\mapsto 0$. It is endowed with a Frobenius endomorphism $\sigma:T\to T$. Let $T_1=T\otimes_{W(k)}k$ with its induced log structure. We have a natural morphism of log schemes $Y\to T_1$ which is log smooth. Let $E$ be a locally constant sheaf (for the Zariski topology) of finite dimensional $K^t$-vector spaces on ${\mathfrak X}$ (or on $Y$). We can then consider the log convergent cohomology groups $H_{conv}^i(Y/T,E)$ of $Y$ relative to $T$, with coefficients in $E$, see \cite{hk}. These are finite dimensional $K^t$-vector spaces endowed with a $K^t$-linear nilpotent monodromy operator $N$. In fact, $N^i=0$ if $Y^{i+1}=\emptyset$ (but not necessarily vice versa) and in particular $N^{d+1}=0$ at any rate. Moreover, if $E$ is endowed with a $\sigma^r$-linear bijective Frobenius endomorphism (for some $r\in{\mathbb N}$), then also $H_{conv}^i(Y/T,E)$ is endowed with a $\sigma^r$-linear bijective Frobenius endomorphism $\Phi$; it satisfies $N\Phi=p^r\Phi N$. For example, if $E=K^t$, the trivial constant sheaf, then there is an isomorphism of $H_{conv}^i(Y/T,K^t)=H_{conv}^i(Y/T)$ with the log crystalline cohomology of $Y$ relative to $T$ (with trivial coefficients) as defined by Hyodo and Kato, such that the respective Frobenius- and monodromy operators coincide. 

Let now $F$ be a locally constant sheaf (for the Zariski topology) of finite dimensional $K$-vector spaces on ${\mathfrak X}$ (or on $Y$); by scalar restriction to $K^t$ it is in particular a sheaf in finite dimensional $K^t$-vector spaces. We can consider the de Rham complex $\Omega_{{\mathfrak X}_{\mathbb Q}}^{\bullet}\otimes_{K}F$ of the generic fibre ${\mathfrak X}_{\mathbb Q}$ (as a $K$-rigid space) of ${\mathfrak X}$ with coefficients in $F$ (we mean: the pull back of $F$ via the specialization map ${\mathfrak X}_{\mathbb Q}\to{\mathfrak X}$). By \cite{hk} we then have for any $i\in{\mathbb Z}$ an isomorphism of $K$-vector spaces \begin{gather}H^i_{dR}({\mathfrak X}_{\mathbb Q},F)=H^i({\mathfrak X}_{\mathbb Q},\Omega_{{\mathfrak X}_{\mathbb Q}}^{\bullet}\otimes_{K}F)\cong H_{conv}^i(Y/T,F)\label{hyokat}\end{gather}(depending on our chosen uniformizer $\pi$). By transport of structures, this isomorphism endows $H^i_{dR}({\mathfrak X}_{\mathbb Q},F)$ with the following additional structures. First of all, the monodromy operator $N$ on $H_{conv}^i(Y/T,F)$ induces an $N$ on $H^i_{dR}({\mathfrak X}_{\mathbb Q},F)$. Secondly, if $F$ carries a $\sigma^r$-linear bijective Frobenius endomorphism, then $H_{conv}^i(Y/T,F)$ and hence $H^i_{dR}({\mathfrak X}_{\mathbb Q},F)$ gets one. Finally, if $F=E\otimes_{K^t}K$ for some locally constant sheaf $E$ of $K^t$-vector spaces (e.g. if $F=K$ take $E=K^t$), then $H_{conv}^i(Y/T,F)= H_{conv}^i(Y/T,E)\otimes_{K^t}K$ and hence $H^i_{dR}({\mathfrak X}_{\mathbb Q},F)$ gets a $K^t$-structure. [In this connection we remark that the first sentence following the statement of Theorem 0.1 in the introduction of \cite{hk} might be misleading, since at that point the need for such an $E$ (to get a $K^t$-structure on $H^i_{dR}({\mathfrak X}_{\mathbb Q},F)$) was not emphasized.] In practice it often happens that there is a Frobenius structure only on such an $E$, not on $F$, hence a Frobenius structure only on the corresponding $K^t$-lattice $H_{conv}^i(Y/T,E)$.

We have a Cech spectral sequence\begin{gather}E_1^{rs}=H_{dR}^s(]Y^{r+1}[,F)\Rightarrow H_{dR}^{s+r}({\mathfrak X}_{\mathbb Q},F).\label{drspec}\end{gather}

Here, by definition,$$H_{dR}^s(]Y^{r+1}[,F)=\coprod_{W}F(W)\otimes_K H_{dR}^s(]W[)$$ where $W$ runs through the connected components of $Y^{r+1}$ and where $F(W)$ denotes the value at $W$ of the restriction of $F$ to $W$; note that this restriction is constant. (In \cite{hk} we used the notation $F(Y^{r+1})\otimes_KH_{dR}^s(]Y^{r+1}[)$ instead.) As $d=\dim({\mathfrak X}_{\mathbb Q})$ we get a boundary morphism$$H_{dR}^0(]Y^{d+1}[,F)\stackrel{\alpha}\longrightarrow H_{dR}^{d}({\mathfrak X}_{\mathbb Q},F).$$On the other hand we have the sum of restriction maps$$H_{dR}^{d}({\mathfrak X}_{\mathbb Q},F)\stackrel{\beta}\longrightarrow H_{dR}^d(]Y^{d+1}[,F).$$Now we define a residue map$$Res:H_{dR}^d(]Y^{d+1}[,F) \longrightarrow H_{dR}^0(]Y^{d+1}[,F).$$For $\sigma\in S_{d+1}$ consider the decomposition $Y_{\sigma}=\coprod_{j=1}^{t(\sigma)}Y_{\sigma,j}$ with $t(\sigma)\in{\mathbb Z}_{\ge0}$ such that for each $Y_{\sigma,j}$ the underlying scheme is just $\spec(k)$. We let $Res$ be the sum of residue maps $F(Y_{\sigma,j})\otimes_K res_{\sigma,j}$ where $$res_{\sigma,j}:H_{dR}^d(]Y_{\sigma,j}[) \longrightarrow H_{dR}^0(]Y_{\sigma,j}[)$$ is defined as follows. Write $\sigma=\{s_0,\ldots,s_d\}$ with $s_0<\ldots<s_d$ in our fixed total ordering of $R$. We find $t_0,\ldots,t_d\in{\mathcal O}_{\mathfrak{X}}(\mathfrak{U})$ for some open affine neighbourhood $\mathfrak{U}\subset \mathfrak{X}$ of $Y_{\sigma,j}$ such that $Y_{s_i}\cap \mathfrak{U}\subset \mathfrak{U}$ is defined by $t_i$ (any $i$) and such that $T_i\mapsto t_i$ defines an \'{e}tale morphism $$\mathfrak{U}\longrightarrow{\rm Spf}({\mathcal O}_K\{T_0,\ldots,T_d\}/(T_0\cdots T_d-\pi))$$(with ${\mathcal O}_K\{T_0,\ldots,T_d\}$ denoting the $\pi$-adic completion of ${\mathcal O}_K[T_0,\ldots,T_d]$). Let ${\mathbb C}_p$ denote the completion of an algebriac closure of $K$. The set of ${\mathbb C}_p$-valued points of the tube $]Y_{\sigma,j}[$ is then given (since $Y_{\sigma,j}$ is the common vanishing locus of all $t_i$) as\begin{gather}]Y_{\sigma,j}[({\mathbb C}_p)=\{y\in \mathfrak{U}_{\mathbb Q}({\mathbb C}_p)\quad|\quad |t_i(y)|<1\mbox{ for }i=1,\ldots, d\mbox{ and }|t_1(y)\cdots t_d(y)|>|\pi|\}.\label{gannulli}\end{gather}In fact, $]Y_{\sigma,j}[$ maps isomorphically to the tube of the common vanishing locus (this vanishing locus is just a $k$-rational point) of all $T_i$ in the generic fibre of ${\rm Spf}({\mathcal O}_K\{T_0,\ldots,T_d\}/(T_0\cdots T_d-\pi))$, and this latter tube is characterized by the analogous inequalities in terms of the $T_i$. Therefore the local analysis of \cite{schtei} p.404 -- 407 carries over verbatim: $H_{dR}^d(]Y_{\sigma,j}[)$ is one-dimensional, with basis the class of the $d$-form$$\dlog(t_1)\wedge\ldots\wedge\dlog(t_d).$$We now define $res_{\sigma,j}$ as the $K$-linear map which sends the class of $\dlog(t_1)\wedge\ldots\wedge\dlog(t_d)$ to $1\in K=H_{dR}^0(]Y_{\sigma,j}[)$. The definition is independent on the choice of $t_0,\ldots,t_d$. Indeed, for a unit $\epsilon\in {\mathcal O}^{\times}_{\mathfrak{X}}(\mathfrak{U})$ the restriction of\begin{gather}\dlog(t_1)\wedge\ldots\wedge\dlog(t_{i-1})\wedge\dlog(\epsilon)\wedge\dlog(t_{i+1})\wedge\ldots\wedge\dlog(t_d)\label{closed}\end{gather}to $]Y_{\sigma,j}[$ is closed. To see this observe that $\epsilon|_{]Y_{\sigma,j}[}$ may be expanded into a convergent {\it power} (not just Laurent) series in $t_1,\ldots t_d$, we therefore find another convergent (on $]Y_{\sigma,j}[$) power series $\tilde{\epsilon}$ whose derivative with respect to the variable $t_i$ is $\epsilon$. Thus$$\tilde{\epsilon}\dlog(t_1)\wedge\ldots\wedge\dlog(t_{i-1})\wedge\dlog(t_{i+1})\wedge\ldots\wedge\dlog(t_d)$$is (up to sign) a primitive for (\ref{closed}). It follows that $\dlog(t_1)\wedge\ldots\wedge\dlog(t_d)$ and $\dlog(t_1)\wedge\ldots\wedge\dlog(\epsilon t_i)\wedge\ldots\wedge\dlog(t_d)$, whose difference is the form (\ref{closed}), define the same class in $H_{dR}^d(]Y_{\sigma,j}[)$.  

\begin{satz}\label{nviares} The composition $$\alpha\circ Res\circ\beta: H_{dR}^{d}({\mathfrak X}_{\mathbb Q},F)\longrightarrow H_{dR}^{d}({\mathfrak X}_{\mathbb Q},F)$$coincides (up to sign) with $N^d$, the $d$-fold iterated monodromy operator.
\end{satz}

{\sc Proof:} (For $d=1$ and $F=K$ this has been proved in \cite{coliov}.) First we trace the way of $\alpha\circ Res\circ\beta$ through the isomorphism (\ref{hyokat}). Besides the logarithmic convergent cohomology $H_{conv}^d(Y/T,F)$ relative to the logarithmic basis $T=(\spf(W(k)),[{\mathbb N}\to W(k), 1\mapsto0])$ one may consider the logarithmic convergent cohomology $H_{conv}^d(Y/S,F)$ relative to the logarithmic basis $S=(\spf({\mathcal O}_K),[{\mathbb N}\to {\mathcal O}_K, 1\mapsto\pi])$. As ${\mathfrak X}$ is naturally a log smooth formal $S$-log scheme lifting $Y$, its de Rham cohomology $H_{dR}^{d}({\mathfrak X}_{\mathbb Q},F)$ identifies (practically by the definitions) with $H_{conv}^d(Y/S,F)$. For the same reason, the Cech spectral sequence\begin{gather}E_1^{rs}=H_{conv}^s(Y^{r+1}/S,F)\Rightarrow H_{conv}^{s+r}(Y/S,F)\label{sspec}\end{gather}gets naturally identified with (\ref{drspec}). Consider the spectral sequence\begin{gather}E_1^{rs}=H_{conv}^s(Y^{r+1}/T,F)\Rightarrow H_{conv}^{s+r}(Y/T,F)\label{tspec}.\end{gather}By the main result of \cite{hk} the isomorphism (\ref{hyokat}) comes in fact with an isomorphism of spectral sequences (\ref{sspec})$\cong$ (\ref{tspec}) (depending on $\pi$). In particular, this isomorphism provides us with the vertical isomorphisms in the following diagram:$$\xymatrix{H_{dR}^{d}({\mathfrak X}_{\mathbb Q},F)\ar[r]^{\beta}\ar[d]^{=}&H_{dR}^d(]Y^{d+1}[,F)\ar[d]^{=}\ar[r]^{Res}&H_{dR}^0(]Y^{d+1}[,F)\ar[d]^{=}\ar[r]^{\quad\alpha}& H_{dR}^{d}({\mathfrak X}_{\mathbb Q},F)\ar[d]^{=}\\       H_{conv}^d(Y/S,F)\ar[d]^{\cong}\ar[r]^{\beta}& H_{conv}^d(Y^{d+1}/S,F)\ar[d]^{\cong}\ar[r]^{Res} & H_{conv}^0(Y^{d+1}/S,F)\ar[d]^{\cong}\ar[r]^{\quad\alpha}& H_{conv}^d(Y/S,F)\ar[d]^{\cong}\\H_{conv}^d(Y/T,F)\ar[r]^{\beta'}&H_{conv}^d(Y^{d+1}/T,F)\ar[r]^{Res'} & H_{conv}^0(Y^{d+1}/T,F)\ar[r]^{\quad\alpha'}& H_{conv}^d(Y/T,F)}$$Here $\alpha'$ is the boundary map in (\ref{tspec}) and $\beta'$ is the natural restriction map, thus the outer squares commute. Let us compute the map $Res'$ which makes the inner square commute. Let again $\sigma=\{s_0,\ldots,s_d\}\in S_{d+1}$ with $s_0<\ldots <s_d$, consider the above decomposition $Y_{\sigma}=\coprod_{j=1}^{t(\sigma)}Y_{\sigma,j}$ and fix some $1\le j\le t(\sigma)$. Let $V=\spec(W(k)[X])$, endowed with the log structure defined by $X$. Then we have exact closed embeddings $T\to V$ and $S\to V$ defined by $X\mapsto 0$, resp. by $X\mapsto \pi$. We endow $B=\spec(W(k)[T_0,\ldots,T_d])$ with the log structure defined by the chart ${\mathbb N}^d\to B$, $1_i\mapsto T_i$, and view it as a $V$-log scheme via $X\mapsto T_0\cdots T_d$. We then have an exact closed immersion of $V$-log schemes $Y_{\sigma,j}\to B$ as follows: on underlying schemes it is just sending $Y_{\sigma,j}=\spec(k)$ to the common vanishing locus of $p,T_0,\ldots,T_d$. On the level of log structures it sends $T_i$ to a local equation $t_i$ for $Y_{s_i}$ in $Y$ around $Y_{\sigma,j}$ (note that a chart for the log structure of $Y_{\sigma,j}$ is given be the free abelian group on such a system $t_0,\ldots,t_d$). Of course, $Y_{\sigma,j}\to B$ factors through exact closed embeddings $Y_{\sigma,j}\to B\times_VS=B_S$ and $Y_{\sigma,j}\to B\times_VT=B_T$. Since $B_S$ is log smooth over $S$, the log convergent cohomology $H_{conv}^i(Y_{\sigma,j}/S)$ may be computed as the log de Rham cohomology of the tube $]Y_{\sigma,j}[_{B_S}$ of $Y_{\sigma}$ in the generic fibre of the $\pi$-adic completion of $B_S$. Similarly, since $B_T$ over $T$, the log convergent cohomology $H_{conv}^i(Y_{\sigma,j}/T)$ may be computed as the log de Rham cohomology of the tube $]Y_{\sigma,j}[_{B_T}$ of $Y_{\sigma,j}$ in the generic fibre of the $p$-adic completion of $B_T$. We saw that $]Y_{\sigma,j}[_{B_S}$ is given by (\ref{gannulli}) and we defined the residue map $res_{\sigma,j}:H_{dR}^d(]Y_{\sigma,j}[) \to H_{dR}^0(]Y_{\sigma,j}[)$ by sending the class of $\dlog(T_1)\wedge\ldots\wedge\dlog(T_d)$ to $1$. On the other hand let $B_T^0$ be the exact closed formal $T$-subscheme of $B_T$ defined by $T_0=\ldots=T_d=0$; the underlying formal scheme of $B_T^0$ is $\spf(W(k))$. Its generic fibre $B_{T,\mathbb Q}^0$ is an exact Zariski closed rigid subspace of  $]Y_{\sigma,j}[_{B_T}$. Let $\Omega^{\bullet}_{]Y_{\sigma,j}[_{B_T}}$ resp. $\Omega^{\bullet}_{{B_{T,\mathbb Q}^0}}$ denote the log de Rham complexes relative to $T\otimes{\mathbb Q}$. By \cite{hk} Proposition 4.2 or \cite{mono} Proposition 1.8 the restriction maps$$H^i_{conv}(Y_{\sigma,j}/T)=H^i(]Y_{\sigma,j}[_{B_T},\Omega^{\bullet}_{]Y_{\sigma,j}[_{B_T}})\longrightarrow H^i( B_{T,\mathbb Q}^0, \Omega^{\bullet}_{{B_{T,\mathbb Q}^0}})=\Omega^{i}_{{B_{T,\mathbb Q}^0}}$$are isomorphisms. Now $H^d( B_{T,\mathbb Q}^0, \Omega^{\bullet}_{{B_{T,\mathbb Q}^0}})=\Omega^{d}_{{B_{T,\mathbb Q}^0}}$ is the one-dimensional $K^t$-vector space with basis $\dlog(T_1)\wedge\ldots\wedge\dlog(T_d)$, while $H^0( B_{T,\mathbb Q}^0, \Omega^{\bullet}_{{B_{T,\mathbb Q}^0}})=K^t$. Using these isomorphisms, the $K^t$-linear map $$\Omega^{d}_{{B_{T,\mathbb Q}^0}}\longrightarrow K^t,\quad\quad\quad\dlog(T_1)\wedge\ldots\wedge\dlog(T_d)\mapsto1$$defines a $K^t$-linear map $$Res'_{\sigma,j}:H_{conv}^d(Y_{\sigma,j}/T)\longrightarrow H_{conv}^0(Y_{\sigma,j}/T)$$for any $\sigma$, $j$. We claim that the sum$$H_{conv}^d(Y^{d+1}/T,F)=\coprod_{\sigma\in S_{d+1}}\coprod_{j\in\{1,\ldots,t_{\sigma}\}}F(Y_{\sigma,j})\otimes H_{conv}^d(Y_{\sigma,j}/T)$$$$\stackrel{\coprod_{\sigma,j} {\rm id}\otimes Res'_{\sigma,j}}{\longrightarrow} \coprod_{\sigma\in S_{d+1}}\coprod_{j\in\{1,\ldots,t_{\sigma}\}}F(Y_{\sigma,j})\otimes H_{conv}^0(Y_{\sigma,j}/T)=H_{conv}^0(Y^{d+1}/T,F)$$is the map $Res'$ in question. Following the definitions of the residue maps involved, all we need to show is that for any $\sigma$, $j$, the isomorphism $H_{conv}^d(Y_{{\sigma},j}/S)\cong H_{conv}^d(Y_{{\sigma},j}/T)  \otimes_{K^t}K$ identifies the respective classes which the global section $\dlog(T_1)\wedge\ldots\wedge\dlog(T_d)\in\Omega_{B/V}^d$ induces on either side. But this is immediately clear from the construction of this isomorphism. Indeed, in \cite{hk} we construct a certain (formal) 'log scheme with boundary' $\overline{B}$ together with an open dense immersion $B\to\overline{B}$ and an extension of the relative logarithmic de Rham complex $\Omega_{B/V}^{\bullet}$ to a complex $\Omega_{\overline{B}}^{\bullet}$ on $\overline{B}$, in such a way that $\dlog(T_1)\wedge\ldots\wedge\dlog(T_d)$ extends to a (closed) form in $\Omega_{\overline{B}}^{d}$. The restriction maps $$H_{conv}^i(Y_{{\sigma},j}/S)\longleftarrow H^i(\overline{B},\Omega_{\overline{B}}^{\bullet})\otimes_{K^t}K \longrightarrow H_{conv}^i(Y_{{\sigma},j}/T)  \otimes_{K^t}K$$ induced by $$Y_{\sigma,j}\to B_S\longrightarrow B\longrightarrow \overline{B}\longleftarrow B\longleftarrow B_T\longleftarrow Y_{\sigma,j}$$are bijective and this provides the said isomorphism.
        
To prove the theorem it remains to check $\alpha'\circ Res'\circ\beta'=N^d$ (up to sign) as endomorphisms of $H_{conv}^d(Y/T,F)$. This holds true for any locally constant sheaf of $K^t$-vector spaces $E$ on $Y$ (like $F$, viewed as a sheaf of $K^t$-vector spaces). In the constant coefficient case, $E=K^t$, this is shown in \cite{mono} Theorem 4.2 (where, in fact, more generally $N^i$ for any $i\ge0$ is described in a similar way). But the reasoning leading to \cite{mono} Theorem 4.2 relies entirely on (Zariski-)local constructions and is therefore valid for general $E$ as above. For the convenience of the reader, we indicate the crucial points in section \ref{steenbrink} below. \hfill$\Box$\\

\section{Appendix: The monodromy operator}
\label{steenbrink}

We recall constructions from \cite{hk} p.420--422. We keep the notations used in the proof of Theorem \ref{nviares}. Recall that we endow $V=\spec(W(k)[X])$ with the log structure defined by $X$. By (severe) abuse of notation we denote its $p$-adic formal completion again by $V$. An {\it admissible lift} of the semistable $k$-log scheme $Y$ is a $p$-adic formal $V$-log scheme ${\mathfrak Z}$ together with an isomorphism of $T_1$-log schemes $Y\cong{\mathfrak Z}\times_{V}T_1$ satisfying the following conditions: on underlying formal schemes ${\mathfrak Z}$ is smooth over $W(k)$, flat over $V$ and its reduction modulo $(p)$ is generically smooth over $k[X]$; the fibre ${\mathfrak Y}$ above $X=0$ is a divisor with normal crossings on ${\mathfrak Z}$, and the log structure on ${\mathfrak Z}$ is defined by this divisor. (Thus, \'{e}tale locally, ${\mathfrak Z}$ looks like our $B$ as considered in the proof of \ref{nviares}.) We denote an admissible lift by $({\mathfrak Z},{\mathfrak Y})$. Locally on $Y$, admissible lifts exist.

Choose an open covering $Y=\cup_{h\in H}U_h$ of $Y$, together with admissible liftings $({\mathfrak Z}_h,{\mathfrak Y}_h)$ of the $U_h$ (so $U_h$ is the reduction of ${\mathfrak Y}_h$). For a subset $G\subset H$ let $U_G=\cap_{h\in G}U_h$. Let $\{U_{G,\beta}\}_{\beta\in \Upsilon_G}$ be the set of irreducible components of $U_G$. For $h\in G$ and $\beta\in \Upsilon_G$ let ${\mathfrak Y}_{h,\beta}$ be the unique $W(k)$-flat irreducible component of ${\mathfrak Y}_h$ with $U_{G,\beta}={\mathfrak Y}_{h,\beta}\cap U_G$. Let ${\mathfrak K}'_G$ be the blowing up of $\times_{{\spec}(W(k))}({\mathfrak Z}_h)_{h\in G}$ along $\sum_{\beta\in \Upsilon_G}(\times_{{\spec}(W(k))}({\mathfrak Y}_{h,\beta})_{h\in G})$, let ${\mathfrak K}_G$ be the complement of the strict transforms in ${\mathfrak K}'_G$ of all ${\mathfrak Y}_{h_0,\beta}\times (\times({\mathfrak Z}_h)_{h\in G-\{h_0\}})$ (i.e. all $h_0\in G$, all $\beta\in \Upsilon_G$), and let ${\mathfrak Y}_G$ be the exceptional divisor in ${\mathfrak K}_G$. It is a normal crossings divisor, and its $W(k)$-flat irreducible components are indexed by $\Upsilon_G$: they are the inverse images of the $W(k)$-flat irreducible components of ${\mathfrak Y}_h$, for any $h\in G$. We denote the $p$-adic formal completions of ${\mathfrak Y}_G$ and ${\mathfrak K}_G$ again by ${\mathfrak Y}_G$ and ${\mathfrak K}_G$. By construction, the diagonal embedding $U_G\to\times_{{\spec}(W(k))}({\mathfrak Y}_h)_{h\in G}$ lifts canonically to an embedding $$U_G\longrightarrow{\mathfrak Y}_G\longrightarrow {\mathfrak K}_G.$$ Viewing ${\mathfrak K}_G$ as a formal $V$-log scheme (with log structure defined by ${\mathfrak Y}_G$), this is an exact closed embedding of (formal) $V$-log schemes. 
Denote by $\tilde{\omega}_{{\mathfrak K}_G}^{\bullet}$ the logarithmic de Rham complex of ${\mathfrak K}_G$ (relative to ${\rm Spf}(W(k))$ endowed with its trivial log structure). Write $\theta=\dlog(X)$ and let $$\tilde{\omega}_{{\mathfrak Y}_G}^{\bullet}=\tilde{\omega}_{{\mathfrak K}_G}^{\bullet}\otimes{\cal O}_{{\mathfrak Y}_G}\quad\quad\quad\quad\quad\omega_{{\mathfrak Y}_G}^{\bullet}=\frac{\tilde{\omega}_{{\mathfrak Y}_G}^{\bullet}}{\tilde{\omega}_{{\mathfrak Y}_G}^{\bullet-1}\wedge\theta}.$$So $\omega_{{\mathfrak Y}_G}^{\bullet}$ is the logarithmic de Rham complex of the morphism of formal log schemes ${\mathfrak Y}_G\to T$. Let ${\mathfrak Y}_{G,\mathbb{Q}}$ be the generic fibre of ${\mathfrak Y}_G$, a $K^t$-dagger space, and let $\tilde{\omega}_{{\mathfrak Y}_{G,\mathbb{Q}}}^{\bullet}$ resp. ${\omega}_{{\mathfrak Y}_{G,\mathbb{Q}}}^{\bullet}$ denote the sheaf complexes on ${\mathfrak Y}_{G,\mathbb{Q}}$ obtained from $\tilde{\omega}_{{\mathfrak Y}_G}^{\bullet}$ resp. ${\omega}_{{\mathfrak Y}_G}^{\bullet}$ by tensoring with $\mathbb{Q}$. Let $E$ be a locally constant sheaf of $K^t$-vector spaces on $Y$. On the admissible open subspace $]U_G[_{{\mathfrak Y}_G}$ of ${\mathfrak Y}_{G,\mathbb{Q}}$ we define the sheaf complexes $$E\otimes_{K^t}\tilde{\omega}_{{\mathfrak Y}_{G,\mathbb{Q}}}^{\bullet}=sp^{-1}(E|_{U_G})\otimes_{K^t}\tilde{\omega}_{{\mathfrak Y}_{G,\mathbb{Q}}}^{\bullet}|_{]U_G[_{{\mathfrak Y}_G}}$$$$E\otimes_{K^t}{\omega}_{{\mathfrak Y}_{G,\mathbb{Q}}}^{\bullet}=sp^{-1}(E|_{U_G})\otimes_{K^t}{\omega}_{{\mathfrak Y}_{G,\mathbb{Q}}}^{\bullet}|_{]U_G[_{{\mathfrak Y}_G}}$$where $sp:]U_G[_{{\mathfrak Y}_G}\to U_G\subset Y$ is the specialization map. For $G_1\subset G_2$ we have natural transition maps $]U_{G_2}[_{{\mathfrak Y}_{G_2}}\to]U_{G_1}[_{{\mathfrak Y}_{G_1}}$. Hence a site $(]U_{G}[_{{\mathfrak Y}_{G}})_{G\subset H}=]U_{\bullet}[_{{\mathfrak Y}_{\bullet}}$ with sheaf complexes $E\otimes_{K^t}\tilde{\omega}_{{\mathfrak Y}_{\bullet}}^{\bullet}$ and $E\otimes_{K^t}{\omega}_{{\mathfrak Y}_{\bullet}}^{\bullet}$ on it. As the ${\mathfrak Y}_{G}$ are log smooth over $T$ and the $U_G\to{\mathfrak Y}_G$ are exact we have$${\mathbb{R}}\Gamma_{conv}(Y/T,E)={\mathbb{R}}\Gamma(]U_{\bullet}[_{{\mathfrak Y}_{\bullet}},E\otimes_{K^t}{\omega}_{{\mathfrak Y}_{\bullet}}^{\bullet}).$$By construction, we have a short exact sequence\begin{gather}0\longrightarrow E\otimes_{K^t}{\omega}_{{\mathfrak Y}_{\bullet}}^{\bullet}[-1]\stackrel{\wedge\theta}{\longrightarrow} E\otimes_{K^t}\tilde{\omega}_{{\mathfrak Y}_{\bullet}}^{\bullet}\longrightarrow E\otimes_{K^t}{\omega}_{{\mathfrak Y}_{\bullet}}^{\bullet}\longrightarrow 0.\label{monodef}\end{gather}By definition, the {\it monodromy operator}$$N:H^i_{conv}(Y/T,E)\longrightarrow H^i_{conv}(Y/T,E)$$(any $i$) is the connecting homomorphism in cohomology associated with (\ref{monodef}).

Now we draw on a construction of Steenbrink (which in a crystalline setting was adapted by Mokrane \cite{mokr}). For $j\ge 0$ let $$P_j\tilde{\omega}_{{\mathfrak K}_{\bullet}}^{k}=\bi(\tilde{\omega}_{{\mathfrak K}_{\bullet}}^{j}\otimes\Omega^{k-j}_{{\mathfrak K}_{\bullet}}\longrightarrow \tilde{\omega}_{{\mathfrak K}_{\bullet}}^{k})$$where $\Omega^{\bullet}_{{\mathfrak K}_{\bullet}}$ denotes the non-logarithmic de Rham complex on the simplicial formal scheme ${\mathfrak K}_{\bullet}$. Then let $$P_j\tilde{\omega}_{{\mathfrak Y}_{\bullet}}^{\bullet}=\frac{P_j\tilde{\omega}_{{\mathfrak K}_{\bullet}}^{\bullet}}{\tilde{\omega}_{{\mathfrak K}_{\bullet}}^{\bullet}\otimes{\mathfrak I}_{{\mathfrak Y}_{\bullet}}}$$where ${\mathfrak I}_{{\mathfrak Y}_{\bullet}}$ is the ideal of ${{\mathfrak Y}_{\bullet}}$ in ${\mathfrak K}_{\bullet}$. On ${\mathfrak Y}_{\bullet}$ these complexes give rise to a filtration $P_{\bullet}\tilde{\omega}_{{\mathfrak Y}_{\bullet}}^{\bullet}$ of $\tilde{\omega}_{{\mathfrak Y}_{\bullet}}^{\bullet}$. On ${\mathfrak Y}_{G,\mathbb{Q}}$ define the double complex $A_G^{\bullet\bullet}$ as follows: let$$A_G^{ij}=sp^{-1}(E|_{U_G})\otimes_{K^t}\frac{\tilde{\omega}_{{\mathfrak Y}_{G,\mathbb{Q}}}^{i+j+1}}{P_j(\tilde{\omega}_{{\mathfrak Y}_{G,\mathbb{Q}}}^{i+j+1})},$$as differentials $A_G^{ij}\to A_G^{(i+1)j}$ take those induced by ${\rm id}\otimes(-1)^jd$, and as differentials $A_G^{ij}\to A_G^{i(j+1)}$ take those induced by $e\otimes\omega\mapsto e\otimes(\omega\wedge\theta)$. Let $A_G^{\bullet}$ be the associated total complex. The augmentation $E\otimes_{K^t}\tilde{\omega}_{{\mathfrak Y}_{G,\mathbb{Q}}}^{\bullet}\to A_G^{\bullet 0}$ defined by $e\otimes\omega\mapsto e\otimes(\omega\wedge\theta)$ induces a quasi-isomorphism $E\otimes_{K^t}{\omega}_{{\mathfrak Y}_{G,\mathbb{Q}}}^{\bullet}\to A_G^{\bullet}$ (see \cite{hk}, p. 422). Combining for varying $G$ we get a quasi-isomorphism $E\otimes_{K^t}{\omega}_{{\mathfrak Y}_{\bullet}}^{\bullet}\to A_{\bullet}^{\bullet}$ of sheaf complexes on ${\mathfrak Y}_{\bullet}$. In particular we see $${\mathbb{R}}\Gamma_{conv}(Y/T,E)={\mathbb{R}}\Gamma(]U_{\bullet}[_{{\mathfrak Y}_{\bullet}},E\otimes_{K^t}{\omega}_{{\mathfrak Y}_{\bullet}}^{\bullet})={\mathbb{R}}\Gamma(]U_{\bullet}[_{{\mathfrak Y}_{\bullet}}, A_{\bullet}^{\bullet}).$$Let $\nu$ be the bihomogeneous endomorphism of bidegree $(-1,1)$ on $A_{\bullet}^{\bullet\bullet}$ such that $(-1)^{j+1}\nu$ is the natural projection $A_{\bullet}^{i,j}\to A_{\bullet}^{i-1,j+1}$. Denote again by $\nu$ the induced endomorphism of $A_{\bullet}^{\bullet}$.

\begin{lem}\label{nviastee} For any $i\in {\mathbb Z}$, the endomorphism $\nu$ of $A_{\bullet}^{\bullet}$ induces the monodromy operator $N$ on $H^i_{conv}(Y/T,E)=H^i(]U_{\bullet}[_{{\mathfrak Y}_{\bullet}}, A_{\bullet}^{\bullet}).$
\end{lem}

{\sc Proof:} The easy proof is the same as in \cite{mokr} 3.18 (and as in the original work of Steenbrink).\hfill$\Box$\\

\begin{lem} In notations of the proof of Theorem \ref{nviares} we have $\alpha'\circ Res'\circ\beta'=N^d$ (up to sign) as endomorphisms of $H_{conv}^d(Y/T,E)$.
\end{lem}

{\sc Proof:} According to Lemma \ref{nviastee} the composite$$E\otimes_{K^t}{\omega}_{{\mathfrak Y}_{\bullet}}^{\bullet}\stackrel{\wedge\theta}{\longrightarrow}A_{\bullet}^{\bullet}\stackrel{\nu^d}{\longrightarrow} A_{\bullet}^{\bullet}\stackrel{\wedge\theta}{\longleftarrow} E\otimes_{K^t}{\omega}_{{\mathfrak Y}_{\bullet}}^{\bullet}$$(the last (like the first) arrow is a quasiisomorphisms, so we invert it) induces $N^d$ in cohomology. By construction, $\nu^d\circ(\wedge\theta):E\otimes_{K^t}{\omega}_{{\mathfrak Y}_{\bullet}}^{\bullet}\to A_{\bullet}^{\bullet}$ factors through the subcomplex $A^{\bullet d}_{\bullet}[d]$ of $A_{\bullet}^{\bullet}$. Consider the isomorphism (in the derived category) $${\rm res'}:A^{\bullet d}_{\bullet}[d]= \frac{ \tilde{\omega}_{{\mathfrak Y}_{\bullet,\mathbb{Q}}}^{\bullet} }{ P_d \tilde{\omega}_{{\mathfrak Y}_{\bullet,\mathbb{Q}}}^{\bullet} }[-1] \cong \bigoplus_{x\in Y^{d+1}}(E|_x)[d]$$obtained by taking residues (cf. \cite{hk} p.422) and applying the Poincar\'{e} Lemma (which states that open disks are de Rham cohomologically trivial, the starting point for the definition of convergent cohomology). Moreover, identifying a sheaf on the discrete set $Y^{d+1}$ with its global sections we get an identification $$\bigoplus_{x\in Y^{d+1}}(E|_x)[d]= H^0_{conv}(Y^{d+1}/T,E).$$The composite ${\rm res'}\circ\nu^d\circ(\wedge\theta)$ induces the map $Res'\circ\beta'$ in cohomology. On the other hand, the composite $$H^0_{conv}(Y^{d+1}/T,E)= \bigoplus_{x\in Y^{d+1}}(E|_x)[d]\stackrel{ {\rm res'}^{-1}}{\cong} A^{\bullet d}_{\bullet}[d]\to A_{\bullet}^{\bullet}\stackrel{\wedge\theta}{\longleftarrow} E\otimes_{K^t}{\omega}_{{\mathfrak Y}_{\bullet}}^{\bullet}$$ is easily seen to induce the connecting map $\alpha'$ in cohomology.\hfill$\Box$\\

\section{Appendix: Cohomology of Abelian schemes}
\label{klaussec}

Let $F$ be a field of characteristic $0$, let $S$ be a smooth $F$-variety and let $f:A\to S$ be an abelian $S$-scheme. Denote by $\Omega^{\bullet}_A$ and $\Omega^{\bullet}_S$ the respective de Rham complexes relative to $\spec(F)$ and by $\Omega^{\bullet}_{A/S}$ the relative de Rham complex for $f$. On the $q$-th relative de Rham cohomology ${\bf R}^qf_*\Omega^{\bullet}_{A/S}$ (any $q\ge0$) we have the Gauss-Manin connection$$\nabla_{GM}:{\bf R}^qf_*\Omega^{\bullet}_{A/S}\longrightarrow{\bf R}^qf_*\Omega^{\bullet}_{A/S}\otimes_{{\mathcal O}_S}\Omega_S^1$$which in the usual way extends to a sheaf complex $({\bf R}^qf_*\Omega^{\bullet}_{A/S}\otimes_{{\mathcal O}_S}\Omega_S^{\bullet},\nabla_{GM})$. We filter this complex by setting$$F^s({\bf R}^qf_*\Omega^{\bullet}_{A/S}\otimes_{{\mathcal O}_S}\Omega_S^{\bullet})=\sum_{t_1+t_2=s}{\bf R}^qf_*\Omega^{\bullet}_{A/S,\ge t_1}\otimes_{{\mathcal O}_S}\Omega_{S,\ge t_2}^{\bullet}$$for $s\ge 0$. (Strictly speaking, we better wrote $\bi[{\bf R}^qf_*\Omega^{\bullet}_{A/S,\ge t_1}\to {\bf R}^qf_*\Omega^{\bullet}_{A/S}]$ instead of ${\bf R}^qf_*\Omega^{\bullet}_{A/S,\ge t_1}$, but by the degeneration of the relative Hodge-to-de Rham spectral sequence for $f$ this is the same.) We write$$gr_F^s({\bf R}^qf_*\Omega^{\bullet}_{A/S}\otimes_{{\mathcal O}_S}\Omega_S^{\bullet})=\frac{F^s({\bf R}^qf_*\Omega^{\bullet}_{A/S}\otimes_{{\mathcal O}_S}\Omega_S^{\bullet})}{F^{s+1}({\bf R}^qf_*\Omega^{\bullet}_{A/S}\otimes_{{\mathcal O}_S}\Omega_S^{\bullet})}.$$

\begin{pro}\label{kkue} Let $q\ge 0$. If $S$ is projective, the spectral sequence\begin{gather}E_1^{st}=H^{s+t}(S,gr_F^s({\bf R}^qf_*\Omega^{\bullet}_{A/S}\otimes_{{\mathcal O}_S}\Omega_S^{\bullet}))\Longrightarrow H^{s+t}(S,{\bf R}^qf_*\Omega^{\bullet}_{A/S}\otimes_{{\mathcal O}_S}\Omega_S^{\bullet})\label{gmho}\end{gather}degenerates in $E_1$. The cohomology group $H^{s}(S,{\bf R}^qf_*\Omega^{\bullet}_{A/S}\otimes_{{\mathcal O}_S}\Omega_S^{\bullet})$ is a filtered direct summand of $H_{dR}^{s+q}(A/F)$, for any $s$.
\end{pro}

{\sc Proof:} By complex or $p$-adic Hodge theory we know that the spectral sequence\begin{gather}E_1^{st}=H^t(A,\Omega^{s}_{A})\Longrightarrow H^{s+t}(A,\Omega^{\bullet}_{A})=H_{dR}^{s+t}(A/F)\label{puho}\end{gather}degenerates in $E_1$. Consider the filtration by subcomplexes$$G^p\Omega^{\bullet}_A=\bi[\Omega_A^{\bullet-p}\otimes_{{\mathcal O}_A}f^*(\Omega_S^p)\longrightarrow\Omega^{\bullet}_A]$$($p\ge0$) of $\Omega^{\bullet}_A$. Its graded pieces are $$gr_G^p\Omega^{\bullet}_A=\frac{G^p\Omega^{\bullet}_A}{G^{p+1}\Omega^{\bullet}_A}=\Omega_{A/S}^{\bullet-p}\otimes_{{\mathcal O}_A}f^*(\Omega_S^p).$$(We put $\Omega_A^i=\Omega_{A/S}^i=0$ for $i<0$.) By the projection formula we get $${\bf R}^{p+q}f_*(gr_G^p\Omega^{\bullet}_A)={\bf R}^qf_*\Omega^{\bullet}_{A/S}\otimes_{{\mathcal O}_S}\Omega_S^p.$$Under this identification, the connecting homomorphisms$$d_1^{\bullet,q}:{\bf R}^{p+q}f_*(gr_G^p\Omega^{\bullet}_A)\longrightarrow{\bf R}^{p+q+1}f_*(gr_G^{p+1}\Omega^{\bullet}_A)$$are those provided by the Gauss-Manin connection on ${\bf R}^qf_*\Omega^{\bullet}_{A/S}$. Thus, for fixed $q$ we have an identification of complexes\begin{gather}({\bf R}^{p+q}f_*(gr_G^p\Omega^{\bullet}_A),d_1^{\bullet,q})_{p\ge 0}= ({\bf R}^qf_*\Omega^{\bullet}_{A/S}\otimes_{{\mathcal O}_S}\Omega_S^{\bullet},\nabla_{GM}).\label{grddr}\end{gather}The stupid filtration $(\Omega^{\bullet}_{A,\ge s})_s$ of $\Omega^{\bullet}_{A}$ giving rise to the spectral sequence (\ref{puho}) induces a filtration on all $gr_G^p\Omega^{\bullet}_A$ and hence on the left hand side of (\ref{grddr}). Under (\ref{grddr}) it corresponds to the filtration $F^{\bullet}$ considered above, giving rise to the spectral sequence (\ref{gmho}).

Using pull back along the zero section of the abelian scheme $A/S$ we split the canonical exact sequence$$0\longrightarrow f^*\Omega_S^{1}\longrightarrow \Omega_A^{1}\longrightarrow \Omega_{A/S}^{1}\longrightarrow0$$to get decompositions\begin{gather}\Omega_A^s=\bigoplus_{s_1+s_2=s}\Omega_{A/S}^{s_2}\otimes_{{\mathcal O}_A}f^*\Omega_S^{s_1}.\label{splitsecdif}\end{gather}

Now we compute \begin{align}\sum_j{\rm dim}H^j(A,\Omega^{\bullet}_A)&=\sum_{i,j}{\rm dim}H^j(A,\Omega^{i}_A)\notag \\
{} & = \sum_{j,p,q}{\rm dim}H^j(A,\Omega^{q}_{A/S}\otimes f^*\Omega_S^p)\notag \\
{} & \ge  \sum_{j,p,q,t}{\rm dim}H^j(S,{\bf R}^qf_*\Omega_{A/S}^t\otimes \Omega_S^p)\notag \\
{} & \ge  \sum_{j,p,q}{\rm dim}H^j(S,{\bf R}^qf_*\Omega_{A/S}^{\bullet}\otimes \Omega_S^p)\notag \\
{} & =  \sum_{j,p,q}{\rm dim}H^j(S,{\bf R}^qf_*(gr_G^p \Omega_{A}^{\bullet} ))\notag \\
{} & \ge  \sum_{j,q}{\rm dim}H^j(S,{\bf R}^qf_*(\Omega_{A}^{\bullet} ))\notag \\
{} & \ge  \sum_{j}{\rm dim}H^j(A,\Omega_{A}^{\bullet}).\notag\end{align}

Here the first equality follows from the degeneration of (\ref{puho}) in $E_1$, and the second equality follows from the decomposition (\ref{splitsecdif}), while all the inequalities result from the existence of various spectral sequences provided by general cohomology theory. These inequalities must then be equalities and these spectral sequences degenerate, in particular the maps (induced by the differential on $\Omega_A^{\bullet}$)$$ H^j(S,{\bf R}^qf_*\Omega_{A/S}^t\otimes \Omega_S^p)\longrightarrow  H^j(S,{\bf R}^qf_*\Omega_{A/S}^{t'}\otimes \Omega_S^{p'})$$for $(t',p')=(t+1,p)$ or $(t',p')=(t,p+1)$ are zero. Now$$gr_F^s({\bf R}^qf_*\Omega^{\bullet}_{A/S}\otimes_{{\mathcal O}_S}\Omega_S^{\bullet})=\bigoplus_{t+p=s}{\bf R}^qf_*\Omega_{A/S}^t\otimes \Omega_S^p$$by the degeneration of the relative Hodge-to-de Rham spectral sequence for $f$. Thus the spectral sequence (\ref{gmho}) degenerates in $E_1$ (and (\ref{gmho}) is a subquotient of (\ref{puho})). 

Moreover, the subquotient $H^{s}(S,{\bf R}^qf_*\Omega^{\bullet}_{A/S}\otimes_{{\mathcal O}_S}\Omega_S^{\bullet})$ of $H_{dR}^{s+q}(A/F)$ is in fact a filtered direct summand: it can be characterized as the subspace of all elements on which for all $m\in{\mathbb Z}$ the map $[m]:A\to A$, multiplication with $m$, induces the multiplication with $m^q$ (so that this graded piece is the image of an idempotent in the group ring ${\mathbb Q}[{\rm End}_{S}(A)]$; the action of this ring respects filtrations; cf. also \cite{demu}).\hfill$\Box$\\

  %@@

\end{document}